%% file: Epstein_Wvolume_arxiv.tex
\documentclass{article}
\usepackage{amsmath,amsfonts,epsfig,xcolor}
\usepackage{mathptmx}

\usepackage{amsmath,amsfonts,latexsym,amssymb,hyperref}
\usepackage{mathrsfs}
\usepackage{verbatim}
\usepackage[utf8]{inputenc}
\usepackage[T1]{fontenc}
\usepackage{ae,aecompl}
\usepackage{braket}

\input{./macros}

\newtheorem{corollary}[theorem]{Corollary}

\renewcommand{\hbar}{\bar{{\mathbb H}}^3}

\newcommand{\CC}{\mathbb C}
\newcommand{\Sph}{\mathbb S}
\newcommand{\R}{\mathbb R}

\newcommand{\Hp}{{\mathbb H}^2}
\newcommand{\Hs}{{{\mathbb H}^3}}

\newcommand{\psl}{{\rm PSL}(2,\CC)}

\renewcommand{\htwo}{{\mathbb H}^2}
\def\eproof{$\Box$ \medskip}
\renewcommand{\area}{{\operatorname{{\bf area}}}}

\newcommand{\ddx}{{\partial_x}}
\newcommand{\ddy}{{\partial_y}}
\renewcommand{\ddz}{{\partial_z}}
\newcommand{\ddzbar}{{\partial_\zbar}}
\newcommand{\hess}{{\operatorname{Hess}}}
\newcommand{\os}{{\operatorname{\bf{OS}}}}
\newcommand{\geu}{{g_{\rm euc}}}

\newcommand{\gfl}{{\mathfrak g}}
\newcommand{\Id}{{\operatorname{Id}}}

\newcommand{\hh}{{\mathfrak h}}
\newcommand{\cb}{{\mathcal A}}
\newcommand{\II}{{\rm II}}

\newcommand{\dt}{\dot}
\newcommand{\cotanh}{{\operatorname{cotanh}}}

\newcommand{\ext}{{\operatorname{ext}}}
\newcommand{\ep}{{\operatorname{ep}}}
\newcommand{\la}{{\langle}}
\newcommand{\ra}{{\rangle}}
\def\a{\alpha}
\def\b{\beta}

\begin{document}

\title{\bf Epstein Surfaces, W-Volume, and the Osgood-Stowe Differential}
 \author{Martin Bridgeman\thanks{Research supported by NSF grant DMS-2005498 and the Simons Fellowship 675497} \   and Kenneth Bromberg\thanks{Research supported by NSF grant DMS-1906095.} }

\date{\today}

\maketitle

\begin{abstract}
In a seminal paper, Epstein introduced the theory of what are now called Epstein surfaces, which construct surfaces in $\Hs$ associated to a conformal metric on a domain in $\chat$. More recently, these surfaces have been used by Krasnov-Schlenker to define the W-volume and renormalized volume associated with a convex co-compact hyperbolic 3-manifold.  In this paper we consider Epstein surfaces, W-volume and renormalized volume in two main parts. In the first, we develop an alternate construction of Epstein surfaces using the Osgood-Stowe differential, a generalization of the Schwarzian derivative. Krasnov-Schlenker showed that the metric and shape operator of a surface in hyperbolic space is naturally  dual to a conformal metric and shape operator on a  projective structure via the hyperbolic Gauss map. We show that this projective shape operator can be derived from the Osgood-Stowe differential. This  approach allows us to give a comprehensive and self-contained development of Epstein surfaces, W-volume, and renormalized volume.  In the second part, we use the theory developed in the first  to prove a number of new results including a generalization of Epstein's univalence criterion, a variational formula for W-volume in terms of the Osgood-Stowe differential, and a use of W-volume to relate the length of the bending lamination of the convex core to the norm of the Schwarzian derivative of the associated univalent map on the conformal boundary.
\end{abstract}

\tableofcontents

\section{Introduction}

In \cite{epstein-envelopes} Epstein gave a beautiful construction of a surface in $\Hs$ associated to a conformal metric on a domain $\Omega \subset \chat = \del \Hs$. This construction has been very influential. For example, in a pair of follow-up papers \cite{epstein:gaussmap,epstein:univalent} Epstein used these surface to give a new univalence condition for conformal maps of the disk. More recently, these surfaces have been used to define the {\em renormalized volume} of a convex co-compact hyperbolic 3-manifold.

In Epstein's original construction the surface was described as an {\em envelope of horospheres}. Epstein gave parametric equations for these {\em Epstein surfaces} and then from this parameterization was able to calculate the first and second fundamental forms of the surface. Here will we give an alternate construction of Epstein's surfaces. In our approach rather than calculate the first and second fundamental forms from a parameterization of the surface, motivated by work of Krasnov-Schlenker, we describe them as ``dual'' to certain forms on the domain $\Omega$ (or more generally on a complex projective structure).  These forms satisfy the classical surface equations of Gauss and Codazzi and therefore are realized as surfaces in $\Hs$. 

The pair of forms on $\Omega$ are a conformal metric and the metric's {\em Osgood-Stowe differential} (see \cite{Osgood-Stowe}), which plays the role of the second fundamental form. The Osgood-Stowe differential generalizes the Schwarzian quadratic differential - if the conformal metric is the hyperbolic metric then the  Osgood-Stowe differential is exactly the Schwarzian differential of the uniformizing map.

After developing the theory of Epstein surfaces we then revisit Epstein's univalence criteria, giving new proofs that also generalize Epstein's results.

In recent work motivated from physics, Krasnov-Shlenker \cite{KS08} introduced the notion of W-volume for convex, compact hyperbolic 3-manifolds. Via Epstein surfaces, W-volume can be interpreted as a smooth functional on the space of conformal metrics on the boundary of conformally compact but infinite volume hyperbolic 3-manifold. In \cite{KS08} and \cite{schlenker-qfvolume}  they prove many fundamental properties of  W-volume and the related functional, renormalized volume $V_R$. These include a variational formula for W-volume and a new proof that $dV_R$ is the Schwarzian derivative of the uniformization map for the conformal boundary. Our approach to Epstein surfaces naturally gives a new variational formula for W-volume in terms of the Osgood-Stowe differential from which the properties of W-volume developed by Krasnov-Schlenker can easily be extended. In particular, we give a new inequality bounding the $L^2$-norm of the Schwarzian derivative of a projective structure by  the length of its associated bending lamination.

The paper has two parts, in the first  we develop an alternative approach to the subject, and in the second we apply this new approach to prove a number of new results. In the first, we introduce a new framework for examining the interconnected areas of Epstein surfaces, W-volume, and renormalized volume. This framework facilitates a straightforward and comprehensive approach to the topic, enabling us to provide  proofs of the main results found in \cite{epstein:gaussmap,epstein:univalent,epstein-envelopes} and \cite{KS08}. For completeness, we include a number of proofs of classical theorems such as Bonnet's theorem on immersed surfaces. In the second part, we use the framework developed in the first, to establish several new results. These include a generalization of Epstein's univalence criterion, a variational formula for W-volume in relation to the Osgood-Stowe differential, and a bound on the norm of the Schwarzian derivative in terms of the length of the bending lamination. We also mention recent work of Quinn  who has extended many of these ideas in higher dimensions (see \cite{Quinn_Conf_Flat}).

{\bf Acknowledgements:} We would like to thank Mladen Bestvina, Francois Labourie, Curt McMullen, Keaton Quinn, Franco Vargas-Pallete, Jean-Marc Schlenker and Rich Schwartz  for useful comments and suggestions.

\section{Results}
We now describe the main results of the paper. As well as these, we will also use our approach to give elementary proofs of many of the fundamental results  on Epstein surfaces and W-volume. Our hope is that this will give a self-contained description to both subjects as well as a unifying approach.

\subsection{Duality, Osgood-Stowe Differential}
If $S$ is a surface embedded in a hyperbolic 3-manifold then the induced metric on $g$ and shape operator $B\colon TS\to TS$ satisfy the {\em Gauss-Codazzi} equations. Conversely, if $g$ is metric and $B$ a bundle endomorpshism that satisfy the Gauss-Codazzi equations then, by a classical theorem of Bonnet, they are the metric and shape operator of a surface in a hyperbolic 3-manifold. Here we will develop a parallel theory for conformal metrics on {\em complex projective structures} where our replacement for the shape operator will be a tensor developed by Osgood-Stowe  in \cite{Osgood-Stowe}. Moreover we will see that there is a duality between these two theories.

If $g$ is Riemannian metric then any conformally equivalent metric can be written as $e^{2\phi} g$ where $\phi$ is a real valued function. The {\em Osgood-Stowe tensor} is the traceless part of the Hessian of $\phi$ and can be considered a measure of the difference between the two metrics.
While this is defined in all dimensions here we will restrict to the case of conformal metrics on a Riemann surface where the tensor can be naturally identified with a quadratic differential on the surface. One very concrete example is when $\Omega$ is a hyperbolic domain in $\CC$. Then the Osgood-Stowe tensor between the hyperbolic metric on $\Omega$ and the Euclidean metric on $\CC$ gives the Schwarzian derivative of the conformal map uniformizing $\Omega$. A slight variation of these ideas gives a version of the Osgood-Stowe tensor for any conformal metric on a {complex projective structure} on a surface.

We will see that this tensor plays a role similar to that of the second fundamental form of a surface embedded in a space. Namely, if $g$ is a conformal metric on a projective structure $\Sigma$ then we can define a {\em projective shape operator} $B(\Sigma, g)$ via the Osgood-Stowe tensor. Then the {\em projective Gauss-Codazzi equations} are 
$$\begin{array}{ll}
d^{\nabla} B = 0 & \qquad \mbox{(Codazzi)}\\
\Tr( B) = -2K( g) & \qquad \mbox{(Gauss)}
\end{array}$$
where $K(g)$ is the Gaussian curvature of $g$, $\nabla$ the Riemannian connection for $g$, and $d^\nabla$ is defined on an endomorphism $A:TS\rightarrow TS$ by
$$d^\nabla A(X,Y) = \nabla_X (A(Y))-\nabla_Y (A(X)) - A([X,Y]).$$
We prove:
\begin{theorem}\label{proj_GC}
A pair $(g,B)$ on a surface $S$ satisfies the projective Gauss-Codazzi equations if and only if there exists a projective structure $\Sigma$ on $S$ such that $g$ is a conformal metric on $\Sigma$ and $B = B(\Sigma,g)$. Furthermore $\Sigma$ is unique.
\end{theorem}

Before describing the duality we review the classical theory of immersed surfaces. Let $M$ be a hyperbolic 3-manifold with Levi-Civita connection $\overline\nabla$ and $S$  an immersed surface in $M$ with  induced metric $g$. Then the {\em shape operator} of $S$ is the endomorphism $B:TS \rightarrow TS$ of the tangent bundle $TS$ of $S$  defined by $B(v) = \overline{\nabla}_v n$ where $n$ is a choice of unit normal vector field to $S$. For $S$ the boundary of a manifold $N$ our convention will be to let $n$  be the outward normal\footnote{The standard convention   is to  define $B(v) = -\overline{\nabla}_v n$ where $n$ is the inward normal, which is equivalent.}.

The shape operator $B$ is symmetric with respect to the metric $g$ and its eigenvalues at a point are  the {\em principal curvatures}.  The pair $(g,B)$ satisfy the {\em Gauss-Codazzi equations in $\Hs$} given by
$$\begin{array}{ll}
d^\nabla B = 0 & \qquad \mbox{(Codazzi)}\\
\det(B) = K(g)+1. & \qquad \mbox{(Gauss)}
\end{array}$$
A classical result due to Bonnet,  referred to as {\em the fundamental theorem for immersed surfaces}, gives the following converse.

\begin{theorem}[Bonnet]
If $S$ is simply connected and the pair $(g, B)$ satisfy the Gauss-Codazzi equations in $\Hs$ then there is an immersion
$$f\colon S\to\Hs$$
with $g= f^* g_\Hs$ and shape operator $B$. Furthermore $f$ is unique up to post-composition with isometries of $\Hs$.
\end{theorem}

In \cite{KS08} Krasnov-Schlenker show that the image of an immersed surface in $\Hs$ under the hyperbolic Gauss map to $\hat\CC$ defines  a pair $(\hat g, \hat B)$ which satisfy   the  projective Gauss-Codazzi equations. There is a natural duality of pairs $(g,B) \longleftrightarrow (\hat g,\hat B)$ which can be defined by 
$$\hat g(X,(\Id+\hat B)Y) = 2g((\Id+B)X,Y) \qquad (\Id+B)(\Id+\hat B) =2\Id.$$
Equivalently we can write the duality by
$$\hat g(X,Y) = g((\Id+B)X,(\Id+B)Y) \qquad \hat B = (\Id+B)^{-1}(\Id-B)$$
and 
$$ g(X,Y) = \frac{1}{4}\hat g((\Id+\hat B)X,(\Id+\hat B)Y) \qquad  B = (\Id+\hat B)^{-1}(\Id-\hat B).$$
We have:
\begin{theorem}[see also {Schlenker, \cite{schlenker:dual}}]
Let  $(g,B)$ and $\left(\hat g,\hat B\right)$ be dual. Then $(g,B)$ satisfies the Gauss-Codazzi equations in $\Hs$ if and only if $\left(\hat g,\hat B\right)$ satisfy the projective Gauss-Codazzi equations.
\end{theorem}

Observe that a local diffeomorpshism $f\colon S \to \chat$ defines a projective structure $\Sigma = (S, f)$ on $S$. If $\hat g$ is a conformal metric on $\Sigma$ with $\hat B = B(\Sigma, \hat g)$ the projective shape operator then by the theorem of Bonnet there is an immersion $f_0\colon S \to \Hs$ with metric and shape operator $(g, B)$ that is dual to $(\hat g, \hat B)$. On the other hand with the exact same setup Epstein (\cite{epstein-envelopes}) described a map $f_\ep \colon S \to \Hs$ as an {\em envelope of horospheres}. We show that these surfaces are the same:
\begin{theorem}\label{epstein_equal}
Let $\Sigma = (S,f)$ be a simply connected projective structure and $\hat g$ a conformal metric on $\Sigma$. Then $f_0 = f_{\ep}$.
\end{theorem}

\subsection{Univalence criteria}
In a pair of papers  (see \cite{epstein:gaussmap,epstein:univalent}) used (what are now called) Epstein surfaces to give  new univalence criteria for conformal maps of the disk to $\chat$ as well as a criterion for when the conformal map extended to a quasiconformal homeomorphism of the Riemann sphere. 
These criteria generalized Nehari's univalence condition (see \cite{Nehari:schwarzian}) as well as the Ahlfors-Weill Extension Theorem  (see \cite{ahlforsweill}). Osgood-Stowe in \cite{OS_nehari} and subsequently  Chuaqui in \cite{Chuaqui_nehari}, reinterpreted Epstein's univalence criteria in terms of the Osgood-Stowe differential and proved generalizations. 
Using a different approach, we  give the following further generalization of Epstein and those of Osgood-Stowe and Chuaqui. Our results are most conveniently stated in terms of the {\em Schwarzian tensor} $Q(\Sigma, g)$ which is a complexification of the Osgood-Stowe tensor.
\begin{theorem}\label{nehari}
Let $\hat g$ be a complete metric on $\Delta$, $f\colon \Delta\to \chat$ a locally univalent map and $\Sigma = (\Delta, f)$ the associated projective structure on $\Delta$. 
\begin{enumerate}
\item If 
$$\|Q(\Sigma,\hat g)\| \leq -\frac{1}{4}K(\hat g)$$
then $f$ is univalent and extends continuously to $\overline\Delta$. 
\item If 
$$\|Q(\Sigma,\hat g)\| < -\frac{1}{4}K(\hat g)$$
then $f$ extends to a homeomorphism of $\hat\CC$.
\item If $K(\hat g) < 0$ and
$$-\frac{4\|Q(\Sigma,\hat g)\|}{K(\hat g)} \leq k < 1.$$
Then $f$ extends to a quasiconformal homeomorphism  with beltrami differential $\|\mu\|_\infty \leq k$. 
\end{enumerate}
\end{theorem}

Setting $\hat g= g_{\Hp} $ we observe that the univalence criterion above gives Nehari's univalence criterion that a locally univalent map $f:\Delta \rightarrow \hat\CC$ is univalent if the Schwarzian $\|Sf\|_\infty \leq \frac{1}{2}$ (see \cite{Nehari:schwarzian}). The quasiconformal criterion similarly generalizes the  Ahlfors-Weill Extension Theorem which proves that if $\|Sf\|_\infty < \frac{k}{2} < \frac{1}{2}$ then $f$ extends to a $(1+k)/(1-k)$-quasiconformal homeomorphism (see \cite{ahlforsweill}).

\subsection{W-volume}
If $\bar X$ is a compact, hyperbolizable 3-manifold with boundary then given any metric $\hat g$ on the $\del \bar X$ there is a complete hyperbolic structure on $X$ (the interior of $\bar X$) that is conformally compactified by a conformal structure in the conformal class of $\hat g$. The conformal boundary of this structure also determines a projective structure on $\del \bar X$ so there is a projective shape operator $\hat B = B(\del \bar X, \hat g)$. The dual pair $(g, B)$ will determine a surface in $X$. For now we will assume that the surface is embedded and therefore bounds a compact submanifold $N \subset X$. Then 
the W-volume of $N$ is defined by
$$W(N) = V(N) - \frac{1}{2}\int_{\partial N} H dA$$
where $H$ is the mean curvature of the boundary $\partial N$.
The W-volume can be considered as a functional on the space of smooth metric on $\del \bar X$ (again here we are ignoring some technicalities that  arise when the surface determined by $(g, B)$ is not embedded). If $\nu$ is a variation of $\hat g$ the we have a formula for the variation $dW(\nu)$ which will be in terms of Schwarzian tensor $Q(\del \bar X, \hat g)$ for the projective boundary, the Beltrami differential $\mu(\nu)$ which describes the infinitesimal change of conformal structure induced by $\nu$, and $dK(\nu)$, the infinitesimal change in curvature induced by $\nu$. 
\begin{theorem}
The variational formula for $W$  is given by
$$dW(v) = \frac{1}{4}\int_{\partial \bar X} dK(v)dA_{\hat g} - 2\Re \int_{\partial \bar X} Q(\del\bar X, \hat g)\cdot\mu(v) .$$
\end{theorem}

For two metrics in the same conformal class the underlying hyperbolic 3-manifolds will be isometric. In each conformal class there is a unique hyperbolic metric so by restricting to these metrics we get a function on $\mathcal{CC}(X)$ the space of convex co-compact hyperbolic metrics on $X$. This is the {\em renormalized volume}, denoted $V_R$. Then our variational formula reduces to the well known variational formula for $V_R$  first proved by Takhtajan-Zograf and Takhtajan-Teo (see \cite{TZ}, \cite{TT}) in their work on the Liouville action.  A proof of the variational formula for $V_R$ using methods similar to those in this paper was given by Krasnov-Schlenker in \cite{KS08}.

\begin{theorem}[\cite{TZ}, \cite{TT},  \cite{KS08}]
The variation of $V_R$ on $\mathcal{CC}(X)$ is given by
$$(dV_R)(v) = \Re \int_{\partial \bar X} Q(X, \hat g)\cdot \mu(v).$$ 
\end{theorem}

\subsection{W-volume for pairs, projective structures}
Given two conformal metrics $g,h$ on $\partial \bar X$ it is natural to define the W-volume of the pair by
$$W(g,h) = W(h) -W(g).$$
One motivation for defining this is that it can be extended to projective structures. In particular if $\cG(\Sigma)$ is the space of smooth conformal metrics on a projective structure $\Sigma$ on a closed surface $S$ of genus $g\geq 2$ then for $g,h \in \cG(\Sigma)$ we define $W(g,h)$. This W-volume for pairs agrees with the W-volume for pairs on $\partial\bar X$ by taking $\Sigma = \Sigma_X$.

In \cite{KS08}, Krasnov-Schlenker proved a monotonicity property of W-volume, that if $g\leq h$ pointwise and both have non-positive curvature, then $W(g) \leq W(h)$ or equivalently that $W(g,h) \geq 0$. We generalize this monotonicity property as well as the scaling property of W-volume in the below theorem. If $g$  is a smooth metric we define the curvature form by $\Omega_g = K(g)dA_g$.

\begin{theorem}\label{Wdiff}
 Let $g,h \in \cG(\Sigma)$ with $h = e^{2u}g$. Then
$$W(g,h) = -\frac{1}{4} \int_{\Sigma} u\left( \Omega_g+ \Omega_h\right) = 
 - \frac{1}{8}(\Omega_g+\Omega_h)\left(\log\left(h/g\right)\right).$$
In particular  for $g \leq h$ pointwise and both  non-positively curved, then $W(g,h) \geq 0$.
\end{theorem}

In \cite{KS08}, Krasnov-Schlenker proved that in the quasifuchsian case, the hyperbolic metric had maximal W-volume over all conformal metrics with the same area. We  prove the following generalization which gives a new proof of the original maximality property of Krasnov-Schlenker. 

\begin{theorem}{}\label{wmax}
 Let $g,h \in \cG(\Sigma)$ with $h = e^{2u}g$ such that  $g,h$ have the same area and $g$ has constant curvature. Then
$$W(g,h) \leq -\frac{1}{4}\|\nabla u\|_2^2$$
where $\nabla$ is the gradient of $ g$. In particular $ h =  g$ is the unique maximum.
\end{theorem}

\subsection{Schwarzian derivatives and bending laminations}
Thurston showed that complex projective structures on a surface are parameterized by pairs $(g, \lambda)$ where $g$ is a hyperbolic metric and $\lambda$ is a measured lamination. When $\Sigma$ is the conformal boundary of a hyperbolic 3-manifold then $g$ is the metric on the boundary convex core facing $\Sigma$ and $\lambda$ is the bending lamination (see \cite{Kamishima:Tan} for details).  Also associated to the projective structure is the Schwarzian quadratic differential $\phi_\Sigma$.
Using W-volume, given a projective structure $\Sigma$ we show the following surprising relation between the norms $\|\phi_\Sigma\|_2, \|\phi_\Sigma\|_\infty,$ and the length $L(\lambda_\Sigma)$  of the bending lamination $\lambda_\Sigma$  of $\Sigma$ given by its Thurston parametrization.

\begin{theorem}\label{newbound}
 Let $\Sigma$ be a projective structure with Thurston parameterization $\left(Y_\Sigma, \lambda_\Sigma\right)$. Then
$$\|\phi_\Sigma\|_2 \leq (1+\|\phi_\Sigma\|_\infty)\sqrt{L(\lambda_\Sigma)}.$$
\end{theorem}

By  Nehari (see \cite{Nehari:schwarzian}) if $\Sigma$ is the quotient of a domain in $\chat$ then $\|\phi_\Sigma\|_\infty \leq 3/2$. Therefore we obtain the following corollary.
\begin{corollary}
 Let $\Sigma$ be the quotient of a domain in $\chat$.  Then 
$$\|\phi_\Sigma\|_2 \leq  \frac{5}{2}\sqrt{L(\lambda_\Sigma)} .$$
\end{corollary}
For example, this situation arises when $\Sigma$ is an incompressible boundary component of a convex co-compact hyperbolic 3-manifold.

\section{Linear algebra} 
As we will be working with conformal metrics on  Riemann surfaces, it will be natural  to describe tensors and forms in terms of the underlying complex structure even when the objects are $\R$-linear. We now describe the elementary linear algebra needed to  do this.

 Let $J$ be the standard (almost) complex structure on $\R^2$. If $\ddx$ and $\ddy$ are the standard basis of $\R^2$ then 
 $$\ddz = \frac12\left(\ddx - i\ddy\right) \mbox{ and } \ddzbar = \frac12\left(\ddx+i\ddy\right)$$
are a basis of $\R^2\otimes \CC$. If $dx$ and $dy$ are the dual basis (with respect to $\R$) of $\R^2$ then
$$dz = dx+idy \mbox{ and } d\zbar = dx -idy$$
are the dual basis (with respect to $\CC$) to $\ddz$ and $\ddzbar$. Note that both $dz$ and $d\zbar$ are $\CC$-linear on $\R^2\otimes \CC$ but when restricted to $\R^2$, $dz$ is $\CC$-linear with respect to the complex structure $\left(\R^2, J\right)$ and $d\zbar$ is $\CC$-anti-linear. It will be useful to write objects in complex coordinates (even if they are not $\CC$-linear with respect to $\left(\R^2, J\right)$).

The following ``dot'' notation will be useful. Let $T$ be an $\R$-linear 2-tensor and $B$ an $\R$-linear map. Then define the 2-tensor $T\cdot B$ by
$$T\cdot B(X,Y) = T(BX,Y).$$

If we fix a metric on $g$ on $\R^2$ we can use this to {\em change the type} of a tensor: 
If we are given either $T$ or $B$ the formula
$$T = g\cdot B$$
determines the other.

We can then define the {\em trace} of $T$ by $\Tr_g(T) = \Tr(B)$. As $B$ depends on the metric so will $\Tr_g(T)$. The metric $g$ is conformal (with respect to the complex structure $(\R^2,J)$) if $J$ is an isometry. This is equivalent to being a multiple of the usual Euclidean metric $\geu$. We note that $\Tr_\geu(T) = 0$ if and only $\Tr_g(T) =0$ for all conformal metrics.

Both 2-tensors on $\R^2$ and linear maps of $\R^2$ to itself have unique $\CC$-linear extensions to $\R^2\otimes \CC$. If the $\R$-linear versions of $T$, $B$ and $g$ satisfy the equation $T = g\cdot B$ then so will their $\CC$-linear extensions.
The change of basis from complex coordinates to real coordinates is given by the matrix
$$\cb = \left(\begin{array}{cc} 1/2& 1/2\\-i/2 & i/2\end{array}\right)$$
so if $T$ is written as a matrix in the basis $\{\ddx, \ddy\}$ then its matrix representation in the complex basis $\{\ddz, \ddzbar\}$ is 
$$\cb^t T \cb.$$
For example, the usual Euclidean inner product $\geu$ on $\R^2$ is given by the identity matrix in the real basis so its $\CC$-linear extension in the complex basis is given by
$$\cb^t \Id \cb = \left(\begin{array}{cc} 0& 1/2\\1/2 & 0 \end{array}\right).$$
In particular $\geu(\ddz, \ddz) = \geu(\ddzbar, \ddzbar) =0$ while $\geu(\ddz, \ddzbar) = \geu(\ddzbar, \ddz) = 1/2$.

More generally if
$$T= \left( \begin{array}{cc} x+y & r +s \\ r - s & x - y \end{array} \right)$$
then
$$\cb^t T \cb  = \frac12\cdot\left(\begin{array}{cc} y-ir & x+is\\x-is& y+ir \end{array}\right).$$
This gives us the following lemma.

\begin{lemma}\label{2tensor}
Let $T$ be a $\CC$-bilinear tensor on $\R^2\otimes\CC$ that is real on $\R^2\otimes \R$ and $g$ a conformal metric. If $T$ is symmetric and traceless over $\R^2\otimes\R$ then
$$T = Qdz^2 + \bar Q d\zbar^2$$
where
$$Q = \frac12\cdot(T(\ddx, \ddx) -iT(\ddx, \ddy)).$$
\end{lemma}

We also get formulas for the trace.
\begin{lemma}\label{traceformula}
Let $T$ be a $\CC$-bilinear tensor on $\R^2\otimes\CC$ that is real on $\R^2\otimes \R$ and $g$ a conformal metric. Then
$$\Tr_g(T) = \frac{T(\ddz, \ddzbar) + T(\ddzbar, \ddz)}{g(\ddz, \ddzbar)}$$
and if $T$ is symmetric
$$\Tr_g(T) = 2\Re \frac{T(\ddz, \ddzbar)}{g(\ddz, \ddzbar)}.$$
\end{lemma}

Given a conformal metric $g = \lambda\cdot \geu$ the area form is the alternating 2-tensor given by the matrix
$$\left(\begin{array}{cc} 0& \lambda\\ -\lambda & 0 \end{array}\right)$$
in the basis $\{\ddx, \ddy\}$ for $\R^2$. In complex coordinates we then have:
\begin{lemma}\label{areaformula}
Let $g$ be a conformal metric. Then the area form for  $g$ is
$$dA_g = 2g(\ddz,\ddzbar) dx\wedge dy = -i g(\ddz, \ddzbar) dz \wedge d\zbar.$$
\end{lemma}

\subsection{Pairing}\label{pair_notation}
We define a pairing between $T$ an $B$ by
$$\langle T, B \rangle_g = \frac{1}{2}\Tr_g(T\cdot B)dA_g$$
where $g$ is any  conformal metric  with respect to $J$. It is clear that the definition is independent of the conformal metric $g$ chosen and the pairing is $\R$-linear in both $T$ and $B$. We will drop the subscript when the conformal class is obvious. Choosing $g = g_{euc}$ by Lemma \ref{traceformula}
$$\langle T, B \rangle = \left(T\left( \ddz,B\ddzbar\right) +T\left(\ddzbar, B\ddz\right)\right)dx\wedge dy.$$
Our normalization is chosen such that
$$\left\langle dz^2, \ddz\otimes d\zbar\right\rangle = dx\wedge dy.$$
This pairing is the usual pairing between quadratic differentials and Beltrami differentials from Teichm\"uller theory. The following is an immediate consequence of the definition.

\begin{lemma}\label{pairing}
Let  $g$ be a conformal metric and $A,B$ linear maps. Then
$$\langle g\cdot A,  B\rangle = \langle g\cdot B,  A\rangle = \frac{1}{2}\Tr(AB)dA_g.$$
Furthermore if $T$ is symmetric and traceless then
$$\langle T,  B\rangle = \langle T,  B_0\rangle.$$
\end{lemma}

{\bf Proof:} As $(g\cdot A)\cdot B = g\cdot (AB)$ and $\Tr_g(g\cdot M) = \Tr(M)$ we have
$$\langle g\cdot A,  B\rangle = \frac{1}{2} \Tr_g((g\cdot A)\cdot B)dA_g = \frac{1}{2}\Tr_g(g\cdot(AB))dA_g = \frac{1}{2}\Tr(AB)dA_g .$$
The first result follows by  symmetry and the second by orthogonality.\eproof

\subsection{Variations}
If $X_t$ is a one-parameter family of tensors then its time zero derivative $\delta X$ is a tensor of the same type. Then $\delta X$ is a {\em variation} of $X$. Often a variation of $X$ will induce a variation of another tensor and we will use the $\delta$-notation to refer to the variation of this other tensor. For example, if $P$ is a linear isomorphism and $\hat P$ is its inverse the variation $\delta P$ of $P$ will induce a variation $\delta \hat P$ of $\hat P$ and by the product rule we have
$$\delta P\hat P  + P \delta \hat P = 0.$$

We now show that the pairing for two metrics satisfy  a variational property. 

\begin{lemma}\label{invert}
Let $g,\hat g$ be metrics on $\R^2$ related by  $\hat g = P^*g$ where $P$ is symmetric with respect to $g$. Equivalently $g = \hat P^*\hat g$ where $\hat P$ is symmetric with respect to $\hat g$ and $P\hat P = \Id$. 
Then  under a smooth variation of $(g,P)$
$$\langle \hat g, \delta \hat P\rangle_{\hat g}+\langle\delta \hat g,\hat P_0\rangle_{\hat g} =-\sgn(\det(P))\left( \langle g, \delta P\rangle_g+\langle \delta g,P_0\rangle_g\right)  $$
where $P_0, \hat P_0$ are the traceless part of $P,\hat P$ respectively.
\end{lemma}

{\bf Proof:}
We differentiate  $P\hat P = \Id$. This gives
$$\delta \hat P   = - P^{-1}\delta P  P^{-1}.$$
By the Cayley-Hamilton Theorem for linear maps on $\R^2$
$$P^2-\Tr(P)P+\det(P)\Id = 0.$$ Multiplying by $ P^{-1}$ we get $$P = 2P_0+\det(P) P^{-1} \qquad \det(P)(P^{-1})_0 = -P_0 \qquad \Id = 2P_0P^{-1}+\det(P)P^{-2}.$$
As $\hat P = P^{-1}$  then
$$\hat P_0 = -\frac{ P_0}{\det(P)}.$$
As $\hat g =  P^* g$ and $P$ is symmetric with respect to $g$, then  $\hat g =  g  \cdot P^2$. We define linear maps $\eta, \hat\eta$ by $\delta g = g \cdot \eta, \delta\hat g = \hat g\cdot \hat\eta$. Then differentiating $\hat g = g\cdot P^2$ gives
$$ \delta \hat g =  \delta g \cdot P^{2}+  g\cdot(\delta P  P +  P\delta  P) = g\cdot (\eta P^2+\delta P  P +  P\delta  P).$$
As $\delta \hat g =  \hat g\cdot\hat \eta  = g\cdot (P^2\hat \eta)$ we get 
$$\hat \eta = P^{-2}\eta P^2+P^{-2}\delta P  P +  P^{-1}\delta  P.$$
By definition of the pairing $\langle\ ,\ \rangle_{\hat g}$ 
$$\langle \hat g, \delta \hat P\rangle_{\hat g}+\langle \delta \hat g,\hat P_0\rangle_{\hat g} =\frac{1}{2} \Tr_{\hat g}(\hat g\cdot\delta \hat P+ \delta\hat g\cdot\hat P_0)dA_{\hat g}=\frac{1}{2} \Tr(\delta \hat P+ \hat \eta\hat P_0)dA_{\hat g}.$$
By definition of the area form we have
$$dA_{\hat g} = |\det(P)|dA_{ g} = \sgn(\det(P))\det(P) dA_g.$$  
Applying $\Tr(XY) = \Tr(YX)$ and the Cayley-Hamilton Theorem
\begin{eqnarray*}
\Tr(\delta \hat P+\hat \eta \hat P_0)dA_{\hat g} &=&-\sgn(\det(P))\Tr\left( \delta P\left(2 P_0 P^{-1}+ \det(P)P^{-2}\right)+\eta P_0\right)dA_{ g} \\
&=& -\sgn(\det(P))\Tr(\delta  P+ \eta  P_0)dA_{ g}.
\end{eqnarray*}
 The result then follows.\eproof

\subsection{Metric variations and the strain tensor}\label{metric-var}
If $g$ is a metric in the conformal class of the usual complex structure on $\R^2$ and $\delta g$ is a variation of $g$ then we define the linear map $\eta\colon \R^2 \to \R^2$ by
$$\delta g = 2g \cdot \eta .$$
A simple calculation shows that the variation of the area form $dA_g$ satisfies 
$$\delta(dA_g) = \Tr(\eta)dA_g.$$
The traceless part $\eta_0$ of $\eta$ is the  {\em strain}  of the variation and it measures the  change of conformal structure. If we complexify, by Lemma \ref{2tensor} we see that
$$\eta_0 = \nu \ddz\otimes d\zbar + \overline{\nu} \ddzbar\otimes dz$$
for some $\nu \in \CC$. We then define the {\em Betrami differential} for the metric variation by $\mu(\delta g) = \nu \ddz\otimes d\zbar$. 

\subsection{Linear algebra on the complexified tangent bundle}
We can extend this discussion to the tangent bundle $TX$ of a Riemann surface $X$. Namely, the one complex dimensional subspaces spanned by $\ddz$ and $\ddzbar$ are invariant under automorphisms of the almost complex structure $\left(\R^2, J\right)$. In particular, the bundle $TX\otimes\CC$ has a decomposition into complex line bundles $T^{1,0} X$ and $T^{0,1} X$ where the former is spanned by $\ddz$ and the later by $\ddzbar$. We will only use when $X$ is a domain in $\CC$ where using complex coordinates will greatly simplify some computations. 

For example if $f\colon X\to \R$ is a smooth function then we can take the usual differential $df$ which is $\R$-valued 1-form (a section of the co-tangent bundle) and extend it to a section of the dual bundle (over $\CC$) of $TX\otimes \CC$ by declaring it to be $\CC$-linear. With this extension we have that
$$df = f_z dz + f_\zbar d\zbar$$
is a section of the dual bundle of $TX \otimes \CC$ that when restricted to $TX \cong TX \otimes \R$ this section is $\R$-valued and is the usual differential of $f$.

Let $g_\phi(,) = e^{2\phi}\langle, \rangle$ be a conformal metric on $X$ where $\langle, \rangle$ is the usual Euclidean inner product on $\R^2$. Let $\nabla^\phi$ be the Levi-Civita connection for $g_\phi$. These both extend to $TX\otimes \CC$ by declaring them to be $\CC$-linear. Then in real coordinates the gradient is
$$\nabla^\phi f = e^{-2\phi}\left(f_x\ddx + f_y\ddy\right).$$
Using the matrix $\cb^{-1}$ to change to complex coordinates this becomes
$$\nabla^\phi f = 2e^{-2\phi}\left(f_\zbar \ddz + f_z \ddzbar\right)$$
which naturally is a section of $TX\otimes \CC$. We have the following elementary calculations.

\begin{prop}\label{hesscalc}
Let $g_\phi = e^{2\phi}\geu$ with Levi-Civita connection $\nabla^\phi$. Then
\begin{itemize}
\item $\nabla^\phi_X Y = \nabla_X Y + X(\phi) Y + Y(\phi) X - \geu(X,Y) \nabla \phi$;

\item $\nabla^\phi_{\ddz} \ddz = 2\phi_z\ddz \qquad \nabla^\phi_{\ddzbar} \ddzbar=2\phi_\zbar\ddzbar \qquad \nabla^\phi_{\ddz} \ddzbar = \nabla^\phi_{\ddzbar} \ddz = 0$;

\item $\hess f\left(\ddz, \ddz\right) = f_{zz} -2\phi_z f_z$.
\end{itemize}
\end{prop}

{\bf Proof:} We take the equation relating $\nabla^\phi$ and $\nabla$ as the definition of $\nabla^\phi$. One then checks that (1) it is a connection, (2) it is torsion free and (3) it is compatible with the metric $g_\phi$. Therefore it must be the Levi-Civita connection. Once we have verified this over $\R$ then the formula extends over $\CC$.

The formulas for the covariant derivatives of $\ddz$ and $\ddzbar$ then follow where we are using that $\geu\left(\ddz, \ddz\right) = \geu\left(\ddzbar, \ddzbar\right) = 0$ and $\geu\left(\ddz, \ddzbar\right) = \geu\left(\ddzbar, \ddz\right) = 1/2$.

To calculate the Hessian we use the formula
$$\hess f\left(X, Y\right) = g_\phi\left(\nabla^\phi_X \nabla^\phi f, Y\right).$$
To evaluate this when $X=Y=\ddz$ we only need to know the $\ddzbar$ term of $\nabla^\phi_\ddz \nabla^\phi f$ and by the previous formulas this is
$$2e^{-2\phi}\left(-2\phi_z  f_z +  f_{zz}\right) \ddzbar$$
and therefore
$$\hess f\left(\ddz, \ddz\right) = 2e^{-2\phi}\left(-2\phi_z  f_z +  f_{zz}\right)g_\phi\left(\ddzbar, \ddz\right) = f_{zz} -2\phi_z f_z.$$
\eproof

\section{The Osgood-Stowe differential}
If $(M,g)$ is a Riemannian manifold and $f\colon M\to \R$ is a smooth function then the {\em Osgood-Stowe differential} $\os(f;g)$ is the traceless part of the 2-tensor
$$\hess f - df\otimes df.$$
If $M$ is  Riemann surface and $g$ is a conformal metric then there is nice formula for $\os(f;g)$ in local coordinates. 
\begin{prop}\label{OSlocal}
Let $\Omega$ be a domain in $\CC$ and $g_\phi = e^{2\phi}\geu$ a conformal metric on $\Omega$. If $f\colon \Omega \to \R$ is a smooth function then
$$\os(f;g_\phi) = Q(f;g_\phi) dz^2 + \overline{Q(f;g_\phi)}d\zbar^2$$
where
$$Q(f;g_\phi) = f_{zz} - 2\phi_z f_z - (f_z)^2.$$
\end{prop}

{\bf Proof:} By Lemma \ref{2tensor}, a symmetric traceless 2-tensor is of the form $Q dz^2 + \overline{Q} d\zbar^2$ so to find $\os(f;g_\phi)$ we only need to find $\os(f;g_\phi)(\ddz, \ddz) = Q$. By Proposition \ref{hesscalc}
$$\hess f(\ddz, \ddz) = f_{zz} - 2\phi_z f_z.$$
As
$$df\otimes df (\ddz,\ddz) = (f_z)^2$$
the proposition follows. \eproof

Let $g_0$ and $g_1$ be conformal metrics on a Riemann surface $X$. Then $g_1 = e^{2\phi}g_0$ for some real valued function $\phi$. We then define
$$\os(g_0,g_1) = \os(\phi; g_0) \qquad \mbox{and} \qquad Q(g_0, g_1) = Q(\phi;g_0).$$
The tensor $Q$ is the {\em Schwarzian tensor} and has been studied in related  papers by Dumas and Schlenker among others (see \cite{dumas2007} and \cite{schlenker:dual}). The Schwarzian tensor is a quadratic differential, i.e.  a section of the tensor product of the holomorphic cotangent bundle of $X$ with itself. {\em Holomorphic quadratic differentials} appear frequently in Teichm\"uller theory, and will play a special role here, but in general the Schwarzian tensor will not be holomorphic. We now state some basic properties of the Schwarzian tensor and include the proofs for completeness.

\begin{prop}\label{OS_properties} The Osgood-Stowe differential satisfies
\begin{itemize}
\item $\os(g_0, g_1) + \os(g_1,g_2) = \os(g_0, g_2)$ \hspace{1in} (cocycle property)
\item $Q(\geu, f^*\geu) = \frac12 Sf$
\item $\os(\geu, g)=0$ if $g$ is a Euclidean metric on $\chat \setminus \{z\}$, a complete metric of constant curvature on a round disk, or a metric of constant positive curvature on $\chat$.
\end{itemize}
\end{prop}

{\bf Proof:} We first prove the cocycle property when $g_0 = \geu$, $g_1 = e^{2\phi_1}\geu$ and $g_2 = e^{2\phi_2}g_1 = e^{2(\phi_1 + \phi_2)}\geu$. By Proposition \ref{OSlocal}
\begin{eqnarray*}
Q(\geu, g_1) & = & (\phi_1)_{zz} -(\phi_1)_z^2 \\
Q(g_1, g_2) &=& (\phi_2)_{zz} -2(\phi_1)_z (\phi_2)_z - (\phi_2)_z^2\\
Q(\geu, g_2) & =& (\phi_1 + \phi_2)_{zz} - (\phi_1 + \phi_2)_z^2
\end{eqnarray*}
and the sum of the first two lines equals the third. This proves the cocycle property in this case.

Note that if $g= g_1 = g_2$ we see that $\os(\geu, g) = -\os(g, \geu)$. Therefore
\begin{eqnarray*}
\os(g_0, g_1) + \os(g_1, g_2) & = & \os(g_0, g_1) + \os(g_1, \geu) + \os(\geu, g_1) + \os(g_1, g_2)\\
& =& \os(g_0, \geu) + \os(\geu, g_2) \\
&= & \os(g_0, g_2)
\end{eqnarray*}
and we have the cocycle property in general.

For the second property we observe that
$$f^*\geu = e^{\log |f_z|^2} \geu = e^{2\phi}\geu$$
where $\phi = \frac12 \cdot(\log f_z + \log \overline{f_z})$. As $\log \overline{f_z}$ is anti-holomorphic our formula for $Q$ gives
$$Q(\geu, f^*\geu)= \left(\frac12\cdot\frac{f_{zz}}{f_z} \right)_z - \left(\frac12 \cdot\frac{f_{zz}}{f_z} \right)^2 = \frac12 Sf.$$

For the last property, in all cases the isometries of $g$ are restrictions of M\"obius transformations and the stabilizer of any point is the full rotation group of the circle.  By the cocycle property
$$\os(\geu, f^*\geu) + \os(f^*\geu, f^*g) = \os(\geu, f^*g).$$
If $f$ is an isometry (and therefore a M\"obius transformation) then $Sf = 0$ so the first term is zero. As $g = f^* g$ this gives
$$\os(f^*\geu, g)  = \os(\geu, g)$$
so $\os(\geu, g)$ is invariant under the action of the isometry group of $g$.

$$\os(\geu, g)(f_* X, f_* Y) =\os(f^*\geu, f^*g)(X,Y) = \os(\geu,g)(X,Y)$$
so at each point $z$, $\os(\geu, g)$ is invariant under a faithful $S^1$-action. As $\os(\geu, g)$ is traceless this implies that it is zero. \eproof

\subsection{The Osgood-Stowe differential for projective structures}
We now show how the Osgood-Stowe differential can be used to define a differential associated to complex projective structures on a Riemann surface. 

A {\em complex projective structure} $\Sigma$ on a surface $S$ is an atlas of charts to $\chat$ with transition maps restrictions of M\"obius transformations. A map between projective structures is an isomorphism if in every chart it is the restriction of a M\"obius transformation. One way to define a complex projective structure is as a pair $(S,f)$ where $S$ is a surface and $f$ is an immersion of $S$ into $\chat$. Then two pairs $(S_0, f_0)$ and $(S_1, f_1)$ are equivalent projective structures if there is a homeomorphism $\phi\colon S_0\to S_1$ and an element $\gamma \in \psl$ such that $\gamma \circ f_0 = f_1\circ \phi$.

Not every projective structure $\Sigma$ can be written as such a pair, however, one sufficient condition is for $S$ to be simply connected. In particular, $\Sigma$ is a projective structure on $S$ then it lifts to a projective structure $\tilde \Sigma$ on $\tilde S$ and this projective structure can be represented by a pair $(\tilde S, f)$. As $\tilde\Sigma$ is the lift of a projective structure on $\Sigma$ the map $f$ will have a certain equivariance. That is there will be a homomorphism $\rho\colon \pi_1(S) \to \psl$ such that for all $\beta\in\pi_1(S)$ we have $f\circ \beta = \rho(\beta) \circ f$. Then $f$ is the {\em developing map} for $\Sigma$ and $\rho$ the holonomy representation. The developing map will be unique up to post composition of maps in $\psl$. Post composing $f$ with $\gamma \in\psl$ will have the effect of conjugating the holonomy representation $\rho$ by $\gamma$.

A projective structure $\Sigma$ also determines a conformal structure $X$ on $S$. Let $g$ be a conformal metric on $\Sigma$. Let $(U,\psi)$ be a projective chart for $\Sigma$ and assume that the image of $U$ is in $\CC$. Then $\os(\psi^*\geu, g)$ defines a 2-tensor on $U\subset \Sigma$. By Proposition \ref{OSlocal} this 2-tensor doesn't depend on the choice of chart and determines a 2-tensor $\os(\Sigma,g)$ on all of $\Sigma$ with associated quadratic differential $Q(\Sigma, g)$. We note that if $\Sigma$ is a projective structure given by locally univalent map $f:\Delta \rightarrow \hat\CC$ then
$$Q(\Sigma, g_{\Hp}) = Q(f^*\geu,g_{\Hp}) = -Q(g_{\Hp}, f^*\geu) = -\frac{1}{2}Sf.$$

If $\Sigma_0$ and $\Sigma_1$ are two different projective structures and $(U, \psi_0)$ and $(U, \psi_1)$ are charts for the same neighborhood $U$ then  $\os(\psi^*_0 \geu, \psi^*_1 \geu)$ defines a 2-tensor on $U$. This doesn't depend on the choice of chart and gives a global tensor $\os(\Sigma_0, \Sigma_1)$ on $X$ with associated quadratic differential $Q(\Sigma_0, \Sigma_1)$.

\begin{theorem}\label{holomorphic}
If $\Sigma_0$ and $\Sigma_1$ are projective structures on a conformal structure $X$ then $Q(\Sigma_0, \Sigma_1)$ is a holomorphic quadratic differential. Furthermore for any projective structure $\Sigma_0$ and any holomorphic quadratic differential $Q$ on $X$ there is a unique projective structure $\Sigma_1$ with $Q(\Sigma_0, \Sigma_1) = Q$.
\end{theorem}

{\bf Proof:} In a neighborhood $U$ of $X$ with charts $(U, \psi_0)$ and $(U, \psi_1)$, $Q(\Sigma_0, \Sigma_1)$ is the Schwarzian of $\psi_1\circ \psi_0^{-1}$ and hence holomorphic.

For the second statement we use that for any chart $(U, \psi_0)$ for $\Sigma$ there is a locally univalent map
$$f\colon \psi_0(U) \to \chat$$
such that the Schwarzian $Sf$ is $2 Q$ on $U$.  Furthermore, $f$ is unique up to post-composition by a M\"obius transformation so if we let $\psi_1 = f\circ \psi_0$ then $(U, \psi_1)$ defines a projective atlas on $X$. \eproof
\subsection{The projective second fundamental form}

We define a {\em projective second fundamental form} of $(\Sigma, g)$ by
$$\II(\Sigma,g) = 2\os(\Sigma,g) - K(g)\cdot g$$
where $K(g)$ is the curvature. We also define the {\em projective shape operator}
 $$B(\Sigma,g) \colon T\Sigma\to T\Sigma$$
by
$$\II(\Sigma,g)(X,Y) = g(B(\Sigma,g) X, Y).$$
We will see that the the second fundamental form for a projective structure plays a similar role as the second fundamental form for an immersed surface in a hyperbolic 3-manifold.

If
$$B\colon TS\to TS$$
is a bundle endomorphism we define the covariant derivative by
$$d^\nabla B(X,Y) = \nabla_X(B(Y)) - \nabla_Y(B(Z)) - B([X,Y]).$$
Then the {\em projective Gauss-Codazzi 
equations} are:
$$\begin{array}{ll}
d^\nabla B = 0 & \qquad \mbox{(Codazzi)}\\
\Tr(B) = -2K & \qquad \mbox{(Gauss)}
\end{array}$$

We first make a preliminary calculation in a local conformal chart.
\begin{lemma}\label{dnabla_calc}
Let $g= e^{2\phi}\geu$ be a conformal metric on an open neighborhood $U$ in $\CC$.  
Let 
$$T= r(z) dz^2 + \overline{r(z)} d\bar z^2 -s(z) \cdot g$$
and $B$ the bundle endomorphism of $TU \otimes\CC$
with $T = g\cdot B$.
Then $d^\nabla B = 0$ if and only if 
$$2r_\zbar = - s_z e^{2\phi}.$$
\end{lemma}

{\bf Proof:} Note that $d^\nabla B$ is alternating so we only need to evaluate $\left(d^\nabla B\right) \left(\ddz, \ddzbar\right)$ and this is determined by its inner product with $\ddz$ and $\ddzbar$. Furthermore, as both $d^\nabla$ and $g$ are $\CC$-linear extensions of $\R$-linear tensors, we have
$$g\left(\left(d^\nabla B\right) \left(\ddz, \ddzbar\right), \ddzbar\right) = \overline{g\left(\left(d^\nabla B\right) \left(\ddzbar, \ddz\right), \ddz\right)}.$$
Therefore $d^\nabla B = 0$ if and only if the above expression is zero.

We can now compute. By the definition of $d^\nabla B$ we have
$$g\left(\left(d^\nabla B\right) \left(\ddz, \ddzbar\right), \ddz\right) = g(\nabla_{\ddz}B(\ddzbar), \ddz) - g(\nabla_{\ddzbar}B(\ddz), \ddz)$$
as $[\ddz,\ddzbar]=0$. To compute the first term we first observe that
$$\ddz g(B(\ddzbar), \ddz) = \ddz T(\ddzbar,\ddz) = -\frac12\left(s_ze^{2\phi} + s(2\phi_z e^{2\phi})\right)$$
and by the compatibility of $\nabla$ with the metric
\begin{eqnarray*}
\ddz g(B(\ddzbar), \ddz) &= &g(\nabla_{\ddz}B(\ddzbar), \ddz) + g(B(\ddzbar), \nabla_\ddz \ddz)\\ & =& g(\nabla_{\ddz}B(\ddzbar), \ddz) + 2\phi_z T(\ddzbar, \ddz) \\
&=& g(\nabla_{\ddz}B(\ddzbar), \ddz) +2\phi_z\left(-\frac{s}2e^{2\phi}\right).
\end{eqnarray*}
Setting the two equations for $\ddz g(B(\ddzbar), \ddz)$ equal to each other and solving for
$g(\nabla_{\ddz}B(\ddzbar), \ddz)$ gives
$$g(\nabla_{\ddz}B(\ddzbar), \ddz)  = -\frac{s_z}2 e^{2\phi}.$$
Applying a similar strategy we find
$$g(\nabla_{\ddzbar}B(\ddz), \ddz) = r_\zbar$$
so $d^\nabla B = 0$ is equivalent to 
$$2r_\zbar = -s_z e^{2\phi}.$$\eproof

We now restate and prove Theorem \ref{proj_GC}.

\medskip

\noindent {\bf Theorem \ref{proj_GC}}
{\em A pair $(g,B)$ on a surface $S$ satisfies the projective Gauss-Codazzi equations if and only if there exists a projective structure $\Sigma$ on $S$ such that $g$ is a conformal metric on $\Sigma$ and $B= B(\Sigma,g)$. Furthermore $\Sigma$ is unique.
}

{\bf Proof:}
First assume that $g$ is a conformal metric on a projective structure $\Sigma$ with a projective shape operator $B(\Sigma,g)$. If we write $\II_\Sigma(g)$ in a local chart as in Lemma \ref{dnabla_calc} then by the definition of $\II_\Sigma(g)$ we have
$$s = K = -4e^{-2\phi}\phi_{z\zbar}$$
(since right hand side is the formula for the curvature) and
$$r = 2Q(\Sigma, g) = 2(\phi_{zz} - \phi_z^2).$$
As $-\Tr(B) = s$ the Gauss equation holds.
It also follows that
$$2r_\zbar = -s_z e^{2\phi}$$
and therefore by Lemma \ref{dnabla_calc}, $d^\nabla B_{\Sigma}(g) = 0$.
Therefore the pair $(g, B(\Sigma,g))$ satisfies the projective Gauss-Codazzi equations.

Now assume that the pair $(g,B)$ satisfies the projective Gauss-Codazzi equations. 

 We again using the local representation of $B$ from Lemma \ref{dnabla_calc}. Then $-2\Tr(B) = s$ and by the Gauss equation
$$s = K = -4e^{-2\phi}\phi_{z\zbar}.$$
 By the Codazzi equation and Lemma \ref{dnabla_calc} we also have
$$2r_\zbar = -s_z e^{2\phi}.$$
Now choose any projective structure $\Sigma_0$ on the conformal structure $X$ given by the pair $(S,g)$. Then $(g,B_{\Sigma_0}(g))$ satisfies the projective Gauss-Codazzi equations so the $\zbar$-derivative of $2Q(\Sigma_0, g)$ is equal to the $z$-derivative of $-2Ke^{2\phi} = -2s e^{2\phi}$. In particular, it is equal to $r_\zbar$ so
$$r - 2Q(\Sigma_0, g)$$
is a holomorphic quadratic differential. Then by Theorem \ref{holomorphic} there exists a unique projective structure $\Sigma$ on $X$ with
$$2Q(\Sigma_0, \Sigma) = r - 2Q(\Sigma_0, g).$$
The cocycle property (Proposition \ref{OS_properties}) then implies
$$2Q(\Sigma, g) = 2Q(\Sigma, \Sigma_0) + 2Q(\Sigma_0, g) = r$$
so $B = B(\Sigma,g)$. 
\eproof

We also note that the following is a direct corollary of Lemma \ref{dnabla_calc}.
\begin{cor}\label{hyp_hol}
Let $g$ be a conformal metric on a projective structure $\Sigma$. Then  $Q(\Sigma,g)$ is holomorphic if and only if $g$ has constant curvature.
\end{cor}

{\bf Proof:} In a local chart $Q(\Sigma, g)$ is of the form $r dz^2$. By Lemma \ref{dnabla_calc} we have that $2r_\zbar = -K_z e^{2\phi}$ where $g = e^{2\phi} \geu$ in the chart. Therefore $r_\zbar \equiv 0$ if and only if $K_z = 0$. However, $K$ is real valued so $K_z = 0$ if and only if $K$ is constant. \eproof

We define the norm of the Schwarzian tensor by
$$\|Q(\Sigma,\hat g)\| = \frac{|Q(\Sigma,\hat g)|}{\hat g}.$$
Then we have the following description of the eigenvalues of $B(\Sigma,\hat g)$.

\begin{lemma}\label{Bevals}
The eigenvalues of $B(\Sigma,\hat g)$ are
$$-K(\hat g) \pm 4\|Q(\Sigma,\hat g)\|.$$
\end{lemma}

{\bf Proof:}  We have by definition, $\II(\Sigma,\hat g) = \hat g \cdot  B(\Sigma,\hat g)$ where
$$\II(\Sigma, \hat g) = 2\os(\Sigma,\hat g)-K(\hat g)\hat g.$$
Therefore $\Tr(B(\Sigma,\hat g)) = -2K(\hat g)$ and $\det( B(\Sigma,\hat g)) = K(g)^2-16\|Q(\Sigma,\hat g)\|^2 .$
It follows that the eigenvalues of $B_{\Sigma}(\hat g)$ are 
$$-K(\hat g) \pm 4\|Q(\Sigma,\hat g)\|.$$ 
\eproof

\section{The Gauss-Codazzi equations in $\Hs$}
Let the pair $(g, B)$ be a metric on $S$ and symmetric bundle map $B\colon TS\to TS$ that satisfy the {\em Gauss-Codazzi equations in $\Hs$}:
$$\begin{array}{ll}
d^\nabla B = 0 & \qquad \mbox{(Codazzi)}\\
\det(B) = K+1 & \qquad \mbox{(Gauss)}
\end{array}$$
where $\nabla$ is the Levi-Civita connection for $g$ and $K$ is the curvature.

If $S$ is an immersed surface in $\Hs$ (or any hyperbolic 3-manifold) with $(g,B)$ the induced metric and shape operator, then it is a classic fact that $(g,B)$ satisfies the Gauss-Codazzi equations in $\Hs$. Conversely, if $(g,B)$ is a pair that satisfies the Gauss-Codazzi equations in $\Hs$ then (locally) $(g,B)$ is the induced metric and shape operator of an immersion. We will derive both of these facts in the course of our discussion.

 We will see that there is a natural {\em dual pair} $(\hat g, \hat B)$ that satisfies the projective Gauss-Codazzi equations.

\subsection{Dual pairs}

Let $(g,B)$ be a metric and tangent bundle isomorphism and assume that $-1$ is not an eigenvalue of $B$.  We define
$$\hat g = \left(\Id +B\right)^* g$$
and we also define (if $-1$ is not an eigenvalue of $B$)
$$\hat B = (\Id + B)^{-1} (\Id - B).$$
This can be inverted. Namely
$$g = \frac14\left(\Id + \hat B\right)^* \hat g$$
and
$$B = \left(\Id +\hat B\right)^{-1} \left(\Id - \hat B\right).$$
We then have the following fact that was also proven independently by Schlenker in a recent preprint (see \cite{schlenker:dual}).
\begin{theorem}\label{duality}
The pair $(g,B)$ satisfy the Gauss-Codazzi equations in $\Hs$ if and only if the dual pair $(\hat g, \hat B)$ satisfy the projective Gauss-Codazzi equations.
\end{theorem}

In \cite{KS08}, Krasnov-Schlenker proved that if $(g,B)$ satisfy the Gauss-Codazzi equations in $\Hs$ then $(\hat g, \hat B)$ satisfies the projective Gauss-Codazzi equations. We  essentially make the same observation as in \cite{schlenker:dual} that  their proof works to give both  implications.
The following lemma will be key to the proof of this.
\begin{lemma}\label{bundle_map}
Let $(M,h)$ be a Riemannian manifold with Levi-Civita connection $\nabla$. If
$$A\colon TM \to TM$$
is a bundle isomorphism and 
$$h_A = A^*h$$
then
$$\nabla^A_X Y = A^{-1}\left( \nabla_X A(Y)\right)$$
is a connection on $M$ compatible with $h_A$. Furthermore $\nabla^A$ is torsion free (and hence the Levi-Civita connection for $h_A$) if and only if $d^\nabla A = 0$.  If $d^\nabla A = 0$ we also have $d^{\nabla^A} A^{-1} = 0$. Finally if $M$ is a surface and $A$ is symmetric then
$$K_A = \frac{K}{\det(A)}$$
where $K$ and $K_A$ are the curvature of $h$ and $h_A$.
\end{lemma}

{\bf Proof:} It is easy to check that the operator $\nabla^A$ defined by the above formula is a connection and is compatible with the metric $h_A$. We also observe that
$$A\left(\nabla^A_X Y -\nabla^A_Y X - [X,Y]\right) = d^\nabla A(X,Y)$$
so $\nabla^A$ is torsion free if and only if $d^\nabla A = 0$. If we rearrange the formula relating the two connections to 
$$\nabla_X Y = A\left(\nabla^A_X A^{-1}(Y)\right)$$
than as $\nabla$ is torsion free the previous statement implies that $d^{\nabla^A} A^{-1} = 0$. 

To find the curvature we first observe that if  $R$ and $R_A$ are the curvature tensors then $R_A = A^{-1}RA$ giving
$$h_A\left(R_A(X,Y)Z, W\right) = h(R(X,Y)(AZ), (AW)).$$

If $M$ is 2-dimensional and $A$ is symmetric then we can choose $X$ and $Y$ to be eigenvectors of $A$. Note that they will be orthogonal with respect to both $h$ and $h_A$. Then
\begin{eqnarray*}
K_A & = & \frac{h_A(R_A(X,Y)X, Y)}{h_A(X,X) h_A(Y,Y)} \\
& = & \frac{h(R(X,Y)(AX), AY)}{\det(A)^2 h(X,X) h(Y,Y)} \\
& = & \frac{\det(A)}{\det(A)^2} \cdot\frac{h(R(X,Y)X, Y)}{h(X,X)h(Y,Y)} \\
& = & \frac{K}{\det(A)}.
\end{eqnarray*}
\eproof

{\bf Proof of Theorem \ref{duality}:} Let $A = \Id + B$ and let $(g,B)$ satisfy the Codazzi equation. Then $d^\nabla A = d^\nabla \Id + d^\nabla B=0$ (since $d^\nabla \Id =0$). Applying Lemma \ref{bundle_map} this implies that $d^{\hat\nabla}A^{-1}  =0 $. But as $A^{-1} = (\Id + B)^{-1}= \left(\Id + \hat B\right)/2$ then $d^{\hat\nabla} \hat B = 0$. It follows that $(\hat g,\hat B)$ satisfies the Codazzi equation.  Letting $A= \left(\Id + \hat B\right)/2$ the same argument gives that if $\left(\hat g, \hat B\right)$ satisfies the Codazzi equation then so does $(g, B)$.

We note that for $A$ a $2\times 2$ matrix,
$$\Tr(A^{-1}) = \frac{\Tr(A)}{\det(A)}.$$
Therefore as $(\Id + B)^{-1} = \left(\Id + \hat B\right)/2$
$$2+\Tr(\hat B) = \Tr(\Id+\hat B) = \frac{2\Tr(\Id+B)}{\det(\Id +B)}.$$
Solving we get 
$$\Tr\left(\hat B\right) = \frac{2(\Tr(\Id+B)-\det(\Id+B))}{\det(\Id +B)} = \frac{2(1-\det(B))}{\det(\Id +B)}.$$
Applying the curvature formula from Lemma \ref{bundle_map} this becomes
$$\Tr\left(\hat B\right) = 2(1-\det(B))\cdot\frac{\hat K}K$$
and therefore $\Tr\left(\hat B\right) = -2\hat K$ if and only if $\det B - 1 = K$. \eproof

We observe the following.

\begin{lemma}\label{forms}
Let $(g,B)$  satisfy the Gauss-Codazzi equations in $\Hs$ and $\left(\hat g,\hat B\right)$ be its dual pair. Let $\epsilon = \sgn(\det(\Id+B)) = \sgn(\det(\Id+\hat B))$. Then  
$$K(g) dA_g = \epsilon K(\hat g) dA_{\hat g} \qquad H(g)dA_g = \frac{\epsilon}{4}\left(1-\det\left(\hat B\right)\right) dA_{\hat g},$$
where $K(g) dA_g$ and $K(\hat g) dA_{\hat g}$ are the curvature forms of the metrics and $H(g) = \Tr(B)/2$ is the mean curvature of $g$.
\end{lemma}

{\bf Proof:}
We note that the equality of the curvature forms follows from Lemma \ref{bundle_map} but we will give a proof of  both that shows they are dual.  As $(\Id+B)(\Id +\hat B) = 2\Id,$ by symmetry we have 
$$\Tr\left(\hat B\right) = \frac{2(1-\det(B))}{\det(\Id + B)} \quad \mbox{and} \quad \Tr(B) = \frac{2\left(1-\det\left(\hat B\right)\right)}{\det\left(\Id + \hat B\right)}.$$
Then as 
$$dA_{\hat g} = |\det(\Id + B)| dA_g \qquad dA_{ g} = \frac{1}{4}|\det(\Id + \hat B)| dA_{\hat g}.$$ 
Therefore
$$\Tr\left(\hat B\right) dA_{\hat g} = 2\epsilon(1-\det(B))dA_g \qquad \Tr(B)dA_{g} = \frac{\epsilon}{2}\left(1-\det\left(\hat B\right)\right)dA_{\hat g}.$$
We then apply  the Gauss equations $\det(B)-1 = K(g)$ and $\Tr\left(\hat B\right) = -2K(\hat g)$.
\eproof

\subsection{Bonnet's theorem}
We now prove Bonnet's theorem that a pair $(g, B)$ on a simply connected surface $S$ represents an immersion in $\Hs$ if and only if they satisfy the hyperbolic Gauss-Codazzi equations. This proof was told to us by F. Labourie, see \cite[Section 2]{Labourie:flat_proj} for his original  argument in the context of real projective structures and cubic differentials. Among its advantages, this approach gives an elegant description of the variation of the fundamental pair of an immersion under normal flow.

First assume that $S$ is immersed in $\Hs$. We use the Minkowski model where $\Hs$ is in $\R^{3,1}$. We then have three Levi-Civita connections: $\nabla$, $\bar\nabla$, and $\hat\nabla$ for $(S,g)$, $\Hs$ and $\R^{3,1}$. We let $N$ be the normal vector field for $S$ and define the vector field $X = \sum x_i \frac{\del}{\del x_i}$ in $\R^{3,1}$. The shape operator $B Y = \bar\nabla_Y N$ is a bundle endomorphism of $TS$ (which we implicitly identify as a sub-bundle of the restriction of $T\R^{3,1}$ to $S$). 

We let $E$ be the restriction of the tangent bundle of $T\R^{3,1}$ to $S$ and $\bar E$ the restriction of the tangent bundle $T\Hs$ to $S$. Then we can write $E = \bar E \oplus \R$ where the $\R$-bundle is spanned by $X$. Then we can write every section of $E$ as $Z + \b X$ where $Z$ is a section of $\bar E$ and $\b$ is a function. With this decomposition we have
$$\hat\nabla_Y \left(\begin{matrix} Z\\ \b \end{matrix} \right) = \left(\begin{matrix} \bar\nabla_Y Z + \b \hat\nabla_Y X \\ - \la \hat\nabla_Y Z, X \ra + Y(\b) \end{matrix} \right) = \left(\begin{matrix} \bar\nabla_Y Z + \b Y \\  \la Z, Y \ra + Y(\b) \end{matrix} \right).$$

We similarly have a decomposition $\bar E = TS \oplus \R$ where the $\R$-bundle is spanned by $N$. Then any section of $\bar E$ can be written $Z + \a N$ where $Z$ is a section of $TS$.
We have
$$\bar\nabla_Y \left(\begin{matrix} Z\\ \a \end{matrix} \right) = \left(\begin{matrix} \nabla_Y Z + \a \bar\nabla_Y N\\ \la \bar\nabla_Y Z, N \ra + Y(\a) \end{matrix} \right)= \left(\begin{matrix} \nabla_Y Z + \a B(Y)\\ -\la Z, B(Y) \ra + Y(\a) \end{matrix} \right).$$

We can then put these two computations together as $E = \bar E \oplus \R = (TS \oplus \R) \oplus \R$ and we write a section of $E$ as $Z + \a N + \b X$.
\begin{equation}\label{connection_def}
\hat\nabla_Y \left(\begin{matrix} Z\\ \a\\ \b \end{matrix} \right) = \left(\begin{matrix} \nabla_Y Z +\a B(Y)+ \b Y\\ -\la Z, B(Y) \ra + Y(\a)  \\  \la Z, Y \ra + Y(\b) \end{matrix} \right).
\end{equation}

While the inner product $\la, \ra$ is defined on all of $E$, in the above expression we only apply it to sections of $TS$. We also observe that instead of letting $E$ be the restriction of $T\R^{3,1}$ to $S$ we can define $E$ is the direct sum of $TS$ and the trivial $\R^2$-bundle over $S$.  Then if $g$ is a metric on $S$ and $B$ is a tangent bundle endomorphism of $TS$ then we can use the above expression to define $\hat\nabla$ as
$$
\hat\nabla_Y \left(\begin{matrix} Z\\ \a\\ \b \end{matrix} \right) = \left(\begin{matrix} \nabla_Y Z +\a B(Y)+ \b Y\\ -g( Z, B(Y)) + Y(\a)  \\  g( Z, Y ) + Y(\b) \end{matrix} \right).
$$

Finally we  define an inner product on $E$ by
\begin{equation}\label{inner_def}
\left\langle \left(\begin{matrix} Y\\ \a_1 \\ \b_1 \end{matrix}\right), \left(\begin{matrix} Z\\ \a_2 \\ \b_2 \end{matrix}\right) \right\rangle = g(Y,Z) + \a_1 \a_2 - \b_1\b_2.
\end{equation}

Once we have defined the connection and metric, most of what is left, is to show, by straightforward computation, that it is flat and compatible with the  metric.  The final part will be using the flat connection to define the immersion.

Denoting the curvature tensors by $R, \overline R,$ and $\hat R$, we then write  $\hat R$ as
$$\hat R(Y,Z) = \left(\begin{matrix} \hat R_{11} & \hat R_{12} & \hat R_{13} \\ \hat R_{21} & \hat R_{22} & \hat R_{23} \\ \hat R_{31} & \hat R_{32} & \hat R_{33} \end{matrix} \right)$$
where the $\hat R_{ij} =\hat R_{ij}(Y,Z)$ are bundle maps between the appropriate sub-bundles of $E$.
\begin{lemma}\label{curv_tensor}
If $B$ is symmetric (with respect to $g$) we have
\begin{itemize}
\item $\hat R_{11}(Y,Z)W =  R(Y,Z)W - g(W, B(Z)) B(Y) +g(W, B(Y)) B(Z)  + g(W,Z) Y - g(W, Y )Z$
\item $\hat R_{12} =  d^\nabla B$
\item $\hat R_{21}(Y,Z) W = g(W, d^\nabla B(Y,Z))$
\end{itemize}
and all other $\hat R_{ij} = 0$.
\end{lemma}

{\bf Proof:} The calculation is straightforward but long so we will leave out some details. Since $\hat R$ is a tensor we can assume that $[Y,Z] = 0$ when calculating $\hat R(Y,Z)$. We first have:
$$\hat\nabla_Y\hat\nabla_Z \left(\begin{matrix} W \\ 0 \\ 0 \end{matrix}\right) = \left(\begin{matrix} \nabla_Y\nabla_Z W -  g(W, B(Z)) B(Y) + g(W, Z)Y \\ -g(\nabla_Z W, B(Y)) - Yg(W, B(Z)) \\  g(\nabla_Z W, Y) + Yg(W,Z) \end{matrix} \right).$$
We obtain a similar formula after swapping $Y$ and $Z$. This directly gives $\hat R_{11}$. For $\hat R_{21}$ and $\hat R_{31}$ we also need to use the compatibility of $\nabla$ with the metric $g$ and the symmetry of $\nabla$.

The computation of the other 6 terms is similar by starting with computing $\hat\nabla_Y\hat\nabla_Z$ for the sections $N$ and $X$ of $E$. We note that the computation for $\hat R_{13}$ uses the symmetry of $\nabla$ while the computation for $\hat R_{31}$ and $\hat R_{13}$ use the symmetry of $B$. We leave the details to the reader. \eproof

\begin{cor}
The connection $\hat \nabla$ is flat if and only if $B$ is symmetric and satisfies the Gauss-Codazzi equations.
\end{cor}

{\bf Proof:} Note that $\bar R(Y,Z) W = g(W,Y)Z - g(W,Z)Y$ since $\Hs$ has constant sectional curvature $\equiv -1$. The equation
$$R(Y,Z)W - \bar R(Y,Z)W - g(W, B(Z))B(Y) + g(W, B(Y))B(Z) = 0$$
is equivalent to the Gauss equation so by Lemma \ref{curv_tensor} $\hat R_{11} = 0$ if and only if $(g,B)$ satisfies the Gauss equation.

Lemma \ref{curv_tensor} also implies that that $\hat R_{12} = 0$ and $\hat R_{21} = 0$ if and only if $(g,B)$ satisfies the Codazzi equation.
\eproof

\begin{lemma}
The inner product $\la, \ra$ is compatible with the flat connection $\hat \nabla$.
\end{lemma}

{\bf Proof:} Again the calculation is straightforward and we leave out some details.
For example we have
$$ W\left\langle \left(\begin{matrix} Y\\ 0\\ 0\end{matrix}\right), \left(\begin{matrix} Z\\ 0 \\ 0 \end{matrix}\right) \right\rangle = Wg(Y,Z) = g(\nabla_W Y, Z) + g(Y, \nabla_W Z)$$
and
\begin{eqnarray*}
& &\left\langle \hat\nabla_W \left(\begin{matrix} Y\\ 0\\ 0\end{matrix}\right), \left(\begin{matrix} Z\\ 0 \\ 0 \end{matrix}\right) \right\rangle  + \left\langle \left(\begin{matrix} Y\\ 0\\ 0\end{matrix}\right), \hat\nabla_W\left(\begin{matrix} Z\\ 0 \\ 0 \end{matrix}\right) \right\rangle\\
& =& \left\langle \left(\begin{matrix} \nabla_W Y\\ -g(Y, B(W)) \\ g(Y,W) \end{matrix}\right), \left(\begin{matrix} Z\\ 0 \\ 0 \end{matrix}\right) \right\rangle  + \left\langle \left(\begin{matrix} Y\\ 0\\ 0\end{matrix}\right), \left(\begin{matrix} \nabla_W Z\\ -g(Z, B(W)) \\ g(Z,W) \end{matrix}\right) \right\rangle\\
&=&  g(\nabla_W Y, Z) + g(Y, \nabla_W Z).
\end{eqnarray*}
This gives compatibility when we have two sections of $TS \subset E$.

We also have
$$ W\left\langle \left(\begin{matrix} Y\\ 0\\ 0\end{matrix}\right), \left(\begin{matrix} 0\\ \a \\ 0 \end{matrix}\right) \right\rangle = 0$$
and
\begin{eqnarray*}
& &\left\langle \hat\nabla_W \left(\begin{matrix} Y\\ 0\\ 0\end{matrix}\right), \left(\begin{matrix} 0\\ \a \\ 0 \end{matrix}\right) \right\rangle  + \left\langle \left(\begin{matrix} Y\\ 0\\ 0\end{matrix}\right), \hat\nabla_W\left(\begin{matrix} Z\\ 0 \\ 0 \end{matrix}\right) \right\rangle\\
&=&  \left\langle \left(\begin{matrix} \nabla_W Y\\ -g(Y, B(W)) \\ g(Y,W) \end{matrix}\right), \left(\begin{matrix} 0\\ \a \\ 0 \end{matrix}\right) \right\rangle  + \left\langle \left(\begin{matrix} Y\\ 0\\ 0\end{matrix}\right), \left(\begin{matrix} \a B(W)\\ W(\a) \\ 0 \end{matrix}\right) \right\rangle\\
& = & -\a  g(Y, B(W)) + \a g(Y, B(W)) = 0.
\end{eqnarray*}
This gives compatibility between a section of $TS$ and a multiple of the section $N$.
The other 4 terms are computed similarly. \eproof

Now that we have proven the connection is flat and compatible with the metric, we next  describe the immersion.

As $E$ is a flat connection compatible with $\la, \ra$ we have projections $\pi_p\colon E \to E_p$ such that $\la, \ra$ is the pullback via $\pi_p$ of the restriction  of $\la, \ra$ to $E_p$. For any $W \in T_q S$ we also have a decomposition of $s_*(q) W$ as a sum of $\hat\nabla_ W s(q)$ and a vector in the kernel of $(\pi_p)_*(s(q))$. It follows that
$$(\pi_p \circ s)_*(W) = (\pi_p)_*(s_*(W)) =(\pi_p)_*(\hat \nabla_W s).$$

The following Lemma gives our immersion as well as the immersions obtained under normal flow.
\begin{lemma}\label{bonnet_immersion}
Let  $A_t, B_t, C_t\colon TS\to TS$ be  tangent bundle endomorphism given by
$$A_t W = \cosh t\cdot W + \sinh t\cdot B(W)\qquad C_t W = \sinh t\cdot W + \cosh t\cdot B(W) \qquad B_t = A_t^{-1}C_t$$
and $s_t$ be the section 
$$s_t = \left(\begin{matrix} 0\\ \sinh t\\ \cosh t\end{matrix}\right).$$
Then $\pi_p\circ s_t \colon (S, A^*_t g) \to (E_p, \la, \ra)$ is an isometric immersion with shape operator  $B_t$. Furthermore $\la \pi_p\circ s_t, \pi_p \circ s_t \ra = -1$.  
In particular  $\pi_p \circ s_0\colon (S,g)\to (E_p, \la, \ra)$ is an isometric immersion with shape operator $B$.  

Finally let $N_t:S\rightarrow E$ be the section 
$$N_t = \left(\begin{matrix} 0\\ \cosh t\\ \sinh t\end{matrix}\right).$$
Then $\pi_p \circ N_t:S\rightarrow E_p$ is the normal vector to the surface $\pi_p \circ s_t:S\rightarrow E_p$.
\end{lemma}

{\bf Proof:} 
We have
$$(\pi_p \circ s_t)_*  W = (\pi_p)_* \left(\hat\nabla_W\left(\begin{matrix} 0\\ \sinh t\\ \cosh t\end{matrix}\right)\right) = (\pi_p)_*\left(\begin{matrix} A_t W\\0\\ 0\end{matrix}\right)$$
and therefore, since $\pi_p$ restricts to an isometry between fibers, we have
$$ \la (\pi_p \circ s)_*  W, (\pi_p \circ s)_* T \ra = g(A_t W,A_t T) = (A^*_t g)(W,T).$$

We now show that $B_t$ is the shape operator.  
For the section $N_t$ we have
\begin{align*}
(\pi_p\circ s_t)_*(B_tW) &= (\pi_p)_*(\hat \nabla_{B_tW} s_t)\\
& = (\pi_p)_*\left(\begin{matrix} A_t(B_tW)\\ 0\\ 0\end{matrix}\right)  = (\pi_p)_*\left(\begin{matrix} C_tW \\ 0\\ 0\end{matrix}\right)  = (\pi_p)_*\left(\hat\nabla_{W}\left(N_t\right) \right)\\
& = (\pi_p\circ N_t)_*(W).
\end{align*}
Note that for all $W$ tangent to $S$ we have $$0=\la N_t, \hat\nabla_W s_t \ra = \la \pi_p\circ N_t, \pi_p(\hat\nabla_W s_t) \ra = \la \pi_p\circ N_t, (\pi_p \circ s_t)_* W\ra$$
Thus $\pi_p \circ N_t(q) \in E_p$ is the normal vector to the surface $\pi_p \circ s_t:S\rightarrow E_p$ at $\pi_p\circ s_t(q)$. Therefore we define the vector field $n_t$  to the surface in $E_p$ given by  $n_t(\pi_p\circ s_t(q))= (\pi_p\circ N_t)(q) \in T_{\pi_p\circ s_t(q)}E_p$. 

As $(E_p, \la, \ra)$ is just another copy of the Minkowski space it also has a Levi-Civita connection $\hat\nabla^p$. In fact this flat connection is the usual Euclidean connection so $\hat\nabla^p_Y Z = Y(Z) = Z_* Y$. If $\bar\nabla^p$ is the  Levi-Civita connection for $\Hs \subset E_p$ then $\bar\nabla^p_Y Z$ is the part of $\hat\nabla^p_Y Z $ that is tangent to $\Hs$.  However, if $n$ is normal to the surface and $W$ is tangent then $\hat\nabla^p_W n$ will be tangent to $S$ so
$$\bar\nabla^p_W n = \hat\nabla^p_W n = n_*(W).$$
Therefore
\begin{align*}
(\pi_p\circ N_t)_*(W) &= (n_t\circ (\pi_p\circ s_t))_*(W)\\
& = (n_t)_* \left((\pi_p\circ s_t)_*(W)\right) \\
& = \bar\nabla^p_{(\pi_p\circ s_t)_*(W)} n_t.
\end{align*}
Thus combining we have
$$(\pi_p\circ s_t)_*(B_tW) = \bar \nabla^p_{(\pi_p\circ s_t)_*(W)} n_t$$
and $B_t$ is the shape operator to $\pi_p\circ s_t: S \rightarrow E_p$.
\eproof

We now prove Bonnet's theorem. 

\begin{theorem}\label{thm-bonnet}
Assume that $S$ is a simply connected surface and the pair $(g, B)$ is a metric and tangent bundle endomorphism of $S$. Then there is an isometric immersion
$$f\colon S \to \Hs$$
with $g = f^* g_{\Hs}$ if and only if $(g,B)$ satisfy the Gauss-Codazzi equations in $\Hs$.
Furthermore $f$ is unique up to post-composition with isometries of $\Hs$.
\end{theorem}

{\bf Proof:} We first assume that $g = f^* g_{\Hs}$ where $f$ is an immersion of $\Hs \subset \R^{3,1}$ and $B$ is the corresponding shape operator. Then the Levi-Civita connection $\hat\nabla$ on $\R^{3,1}$ is given by \eqref{connection_def}. This connection is flat ($\hat R \equiv 0$). By Lemma \ref{curv_tensor}, $\hat R \equiv 0$ if and only if $(g,B)$ satisfy the Gauss-Codazzi equations in $\Hs$.

If $(g, B)$ satisfy the Gauss-Codazzi equations then Lemma \ref{bonnet_immersion} gives an isometric immersion of $(S,g)$ in $E_p$ (for any $p \in S$) with shape operator $B$. Is $(E_p, \la, \ra)$ is a copy of Minkowski space and (again by Lemma \ref{bonnet_immersion}) the image of immersion lies in the hyperboloid this gives the desired immersion in $\Hs$.

Assume that $S$ is embedded in $\Hs \subset \R^{3,1}$.  (As our calculation is local this is sufficient.) The bundle $E$ is then the restriction of the tangent bundle of $\R^{3,1}$ to $S$. We have two flat connections on $E$. The restriction of the flat connection on $\R^{3,1}$ and the connection based on the metric and shape operator. By the discussion at the beginning of the subsection these two connections are the same. The Levi-Civita flat connection on $\R^{3,1}$ is the same as the usual flat connection on $\R^4$. In particular the projections $\pi_p\colon \R^{3,1} \to E_p$ are the identity map where $E_p$ is canonically identified with $\R^{3,1}$. Therefore the map $\pi_p \circ s_0\colon S \to E_p$ extends to an isometry from $\R^{3,1}$ to $E_p$. However, the map $\pi_p\circ s_0$ only depends on the pair $(g, B)$ which gives the uniqueness statement.
\eproof

\subsection{Normal flow}
Let  $f\colon S\to \Hs$ be an immersion with $\hat f\colon S \to T^1\Hs$  a choice of lift to the unit tangent bundle over $\Hs$. Let $\gfl \colon T^1\Hs \times \R \to T^1\Hs$ be the geodesic flow and $\pi\colon T^1\Hs  \to \Hs$ be projection to the basepoint in $\Hs$. We can then define the {\em normal flow} of  $f$ to be the  family of smooth maps $f_t = \pi \circ \gfl_t \circ \hat f:S\rightarrow \Hs$.

In the Minkowski metric,a simple calculation shows that geodesic $\gamma$ in $\Hs$  with $\gamma(0) = p$ and $\gamma'(0) = v$ is 
$$\gamma(t) = \cosh(t)p+\sinh(t)v.$$
It follows that identifying $T^1\Hs$ as a subspace of $T\R^{3,1}= \R^{3,1}\times\R^{3,1}$, that 
\begin{equation}\label{eqn-gflow}
\gfl_t(p,v) = (\cosh(t)p+\sinh(t)v, \sinh(t)p+\cosh(t)v).
\end{equation}

\begin{theorem}\label{thm-normalflow}
The metrics and shape operator for $f_t \colon S\to \Hs$ are $(g_t, B_t)$ where
$g_t = A^*_t g$ and $B_t = (A_t)^{-1} C_t$.
\end{theorem}

{\bf Proof:} 
For an immersion $f:S \rightarrow \Hs\subseteq \R^{3,1}$ let   $n: S\rightarrow \R^{3,1}$ be a choice of unit normal. Therefore by the description of geodesics in $\Hs$
$$f_t(q)  = \cosh(t) f(q) + \sinh(t) n(q).$$
By Lemma \ref{bonnet_immersion} for $f = \pi_p\circ s_0:S \rightarrow \R^{3,1}$  the normal  is $n = \pi_p \circ N_0:S\rightarrow \R^{3,1}$.  Therefore by linearity of $\pi_p$ on fibers
$$f_t(q) = \pi_p\circ(\cosh(t)s_0(q)+\sinh(t)N_0(q)) = (\pi_p\circ s_t)(q).$$
The result then follows by Lemma \ref{bonnet_immersion}.
\eproof

\subsection{A flow on metrics}
If $\hat g$ is a conformal metric on a projective structure we can define scaled metrics by $\hat g_t = e^{2t} \hat g$. Then as we will see below $\hat B_t = e^{-2t}\hat B$ are the projective shape operators for $\hat g_t$. The formulas for the dual metrics are more involved. We now  see that the dual metrics are the normal flow immersions for the original immersion.

\begin{prop}\label{dual flow}
Let $\hat g$ be a conformal metric on a complex projective structure $\Sigma$. Define the scaled metrics $\hat g_t = e^{2t} \hat g$ and with projective shape operators $\hat B_t$ and $(g_t, B_t)$ be the dual pairs. 

Then
$$\hat B_t = e^{-2t} \hat B ,\qquad  g_t = A_t^* g, \quad \mbox{and} \quad B_t = \left(A_t\right)^{-1} C_t.$$
In particular if $f:S\rightarrow \Hs$ is an immersion with fundamental pair $(g,B)$ dual to $(\hat g,\hat B)$ then the  normal flow map $f_t:S\rightarrow \Hs$ has fundamental pair $(g_t,B_t)$ dual to $(e^{2t}\hat g, e^{-2t}\hat B)$. 
 \end{prop}

{\bf Proof:} As the Osgood-Stowe differential is invariant under scaling we have $\os(\Sigma,\hat g) = \os(\Sigma,\hat g_t)$ and it follows that $\hat B_t = e^{-2t} \hat B$. 

By the definition of the dual metric we have
\begin{eqnarray*}
g_t &=& \frac14\left(\Id + \hat B_t\right)^* \hat g_t\\
& =&  \frac14\left(\Id +e^{-2t} \hat B\right)^* e^{2t}\hat g\\
& = & \frac14\left(\Id + e^{-2t}\hat B\right)^* e^{2t} \left(\Id + B\right)^* g\\
& = & \frac14 e^{2t}\left(\Id + B +e^{-2t}\left(\Id - B\right)\right)^* g \\
& = & A^*_t g.
\end{eqnarray*}

By the formula for the dual shape operator we have
\begin{eqnarray*}
B_t &=& \left(\Id + \hat{B}_t\right)^{-1} \left(\Id - \hat B_t\right) \\
&= & \left(\Id + e^{-2t} \hat B\right)^{-1} \left(\Id - e^{-2t}\hat B\right)\\
& = & \left(\Id + e^{-2t} \left(\Id + B\right)^{-1}\left(\Id - B\right)\right)^{-1} \left(\Id - e^{-2t}\left(\Id +B\right)^{-1}\left(\Id - B\right)\right)\\
& = &\left(\Id + B +e^{-2t}\left(\Id - B\right)\right)^{-1}\left(\Id +B - e^{-2t}\left(\Id - B\right)\right)\\
&=& \left(A_t\right)^{-1} C_t.
\end{eqnarray*}
The rest follows from Theorem \ref{thm-normalflow}.
\eproof

We remark that at a given point $p \in S$ the metrics $g_t$ can be singular (not positive definite) in which case the shape operator $B_t$  will not be defined at $p$. However, for any given $p$ this can happen at most two values of $t$.

We have the following description of the principal curvatures of $f_t:S\rightarrow \Hs$.

\begin{corollary}\label{cor-locconvex}
Let $f_t:S\rightarrow \Hs$ have principal curvatures $\lambda_1(t),\lambda_2(t)$. Then 
$$\lambda_i(t) = \frac{\tanh(t)+\lambda_i(0)}{1+\tanh(t)\lambda_i(0)}$$
and $f_t$ is an immersion if $\cotanh(t) \neq -\lambda_i(0)$ for $i=1,2$. Furthermore if 
$$e^{2t} \geq \max\left\{\left|\frac{\lambda_1(0)-1}{\lambda_1(0) +1}\right|, \left|\frac{\lambda_2(0)-1}{\lambda_2(0) +1}\right|\right\}$$
 then $f_t$ is locally convex. \label{Bt-evals}\end{corollary}

{\bf Proof:}
By Theorem \ref{thm-normalflow} $(f_t)^* g_{\Hs} = g_t = A_t^* g$ and the shape operator for $f_t$ is $B_t = A_t^{-1}C_t$.  The formula for $\lambda_i(t)$ follows. 

The map $f_t$ will be locally convex if $B_t$ has non-negative eigenvalues. Equivalently this holds if $\hat B_t = (\Id -B_t)^{-1}(\Id+B_t)$ has eigenvalues in $[-1,1]$. As $\hat B_t = e^{-2t}\hat B$, then it follows that $f_t$ is locally convex if $e^{2t}\geq \max\{|\hat\lambda_1(0)|,|\hat\lambda_2(0)|\}$. Rewriting in terms of $\lambda_i(0)$ we obtain the result.
\eproof

\subsection{The flow space of a pair $(g,B)$}
The one-parameter family of metrics $g_t$ can be combined to form a 3-manifold with a (possibly singular) Riemannian metric. We define the {\em flow space} $\cF_{(g,B)}(I)$ to be the product $S\times I$ with  metric $g_t \oplus dt^2$. We have the following corollary.

\begin{corollary}\label{cor-constant}
Let $(g,B)$ satisfy the Gauss-Codazzi equations  on $S$. Then $\cF_{(g,B)}(I)$ is a hyperbolic manifold wherever the metric $g_t$ is non-singular. If $S$ is simply connected then there is a smooth map $F: S\times I \rightarrow \Hs$ such that $F^*g_{\Hs} = g_t\oplus dt^2$ and $F$ is unique up to post-composition by a M\"obius transformation.
\end{corollary}

{\bf Proof:}
We first observe that for $S$ simply connected then by Theorem \ref{thm-normalflow}, we have the smooth map $F:S\times I \rightarrow \Hs$ defined by $F(q,t) = f_t(q)$ and $F^*g_{\Hs} = g_t\oplus dt^2$.  By uniqueness of immersions,  $F$ restricted to $S\times\{0\}$ is unique up to post-composition by M\"obius transformation. As $F$ is defined by normal flow on this immersion, it follows that $F$ is unique up to post-composition by a M\"obius tranformation. If $S$ is not simply connected then we obtain a hyperbolic structure on the non-singular points of $\cF_{(g,B)}(I)$ by choosing for each $U\subseteq S$ open simply connected the map $F_U:U\times I \rightarrow \Hs$. By uniqueness up to post-composition, this gives an atlas for a hyperbolic structure on the non-singular set.
\eproof

Next we give conditions for the flow space to be complete.
\begin{prop}\label{flow space}
Let $(g,B)$ be a  metric and tangent bundle endomorphism on a surface $S$ that satisfy the Gauss-Codazzi equations in $\Hs$. Assume that $g$ is complete.
\begin{itemize}
\item If the eigenvalues of $B$ are $\ge -1$ then $\cF_{(g,B)}([0,\infty))$ is a complete, hyperbolic 3-manifold with boundary. 

\item If eigenvalues of $B$ are contained in the interval $[-1,1]$ then $\cF_{(g,B)}(\R)$ is a complete, hyperbolic 3-manifold. If we further have that $S$ is simply connected then $\cF_{(g,B)}(\R)$ is isometric to $\Hs$.
\end{itemize}
\end{prop}

{\bf Proof:} For both the first and second bullets the restrictions on the eigenvalues of $B$ imply that $A_t$ (and hence $g_t$) is non-singular. Thus the flow space is a hyperbolic manifold. We now show completeness.

 Let $(p_i, t_i)$ be a Cauchy sequence in $\cF_{(g,B)}(I)$ where $I =[0,\infty)$ or $\R$. The projection
$$\cF_{(g,B)}(I)\to I$$
given by taking the last coordinate
is $1$-Lipschitz which implies that the sequence $t_i$ is convergent and lies in some bounded interval $I\cap [-T, T]$. We also have that the projection
$$\cF_{(g,B)}(I\cap [-T,T]) \to S$$
to the first coordinate is $e^T$-Lipschitz to the $g$ metric on $S$. To see this let $\lambda_i$ be the eigenvalues of $B$ then the eigenvalues of $A_t$ are $\lambda_i(t) = \cosh t+\lambda_i\sinh t$ and the Lipschitz constant is bounded by the reciprocal of the minimum of the eigenvalues for $t\in I\cap [-T, T]$. For $I =[0,\infty)$ we have $\lambda_i \geq -1$ and for $I =\R$,  $\lambda_i \in [-1,1]$. In both cases we have for $t \in I\cap [-T,T]$
$$\lambda_i(t) = \cosh t +\lambda_i\sinh t  \ge \cosh t -\sinh t \ge e^{-T}.$$
In particular, $p_i$ is convergent in $S$ so $(p_i, t_i)$ is convergent in $\cF_{(g,B)}(I)$.
Thus in both cases, $\cF_{(g,B)}(I)$ is complete.

Finally if $S$ is simply connected and  the eigenvalues are in the interval $[-1,1]$ then $\cF_{(g,B)}(\R)$ is a simply connected, complete hyperbolic manifold. Therefore by Corollary \ref{cor-constant}. we have a local isometry $F:S\times\R \rightarrow \Hs$. As a local isometry between complete hyperbolic manifolds is a covering map, then as $S\times \R$ is  simply connected, $F$ must be an isometry. \eproof

We can also recover a result of Epstein.
\begin{cor}[{Epstein \cite[Theorem 3.4]{epstein-envelopes}}] \label{embedding}
Let
$$f\colon S \to \Hs$$
be an immersion and assume that $g = f^*g_\Hs$ is a complete metric on $S$. If the absolute values of the eigenvalues of the shape operator $B$ are $\le 1$ then the map $f$ is an embedding and $S$ is homeomorphic to an open disk.
\end{cor}

{\bf Proof:} The pair $(g,B)$ will satisfy the Gauss-Codazzi equations in $\Hs$ so by Proposition \ref{flow space} there is an isometry $F$ from the flow space $\cF_{(g,B)}(\R)$ to $\Hs$. The restriction $f_0$ of $F$ to $S\times \{0\}$ is an immersion of $S$ with metric $g$ and shape operator $B$ so by Theorem \ref{thm-bonnet} we have that $f_0$ is equal to $f$ after possibly post-composing with an isometry of $\Hs$. Since $f_0$ is an embedding so is $f$. \eproof

\subsection{The hyperbolic Gauss map}
We now take the limit of the geodesic flow to obtain the hyperbolic Gauss map.
Namely, as $t\to\infty$, $\gfl_t(p)$ limits to $\partial \Hs =\chat$ so we define the {\em hyperbolic Gauss map}
$$\gfl_\infty\colon T^1\Hs\to \chat$$
by setting
$$\gfl_\infty = \lim_{t\to\infty} \pi\circ\gfl_t.$$

The map $\gfl_\infty$ is particularly simple in the Minkowski model. We let $\cL \subset \R^{3,1}$ be the light cone
$$\cL = \{ x\in \R^{3,1}\ |\ \la x,x\ra = 0,  \ x\neq 0\}.$$
If $\Hs\subseteq\R^{3,1}$ is the hyperboloid, then the Klein model is  $P(\Hs)   \subseteq \mathbb{RP}^3$, the projective hyperboloid with boundary $P(\cL)$, the projective light cone. Then by the equation for  geodesic flow in $\Hs$ (see equation \ref{eqn-gflow})
$$\gfl_\infty(q,v) = \lim_{t\rightarrow \infty}[ \cosh(t) q + \sinh(t)v] = [q+v].$$
We therefore define $\hat\gfl_\infty: T_1\Hs \rightarrow \cL$ by 
$$\hat\gfl_\infty(q,v) = q+v.$$
The identification $\partial \Hs = \chat$ is  then given by the  map $p\colon \cL \to \Sph^2 \subset \R^3$ 
$$p(x_1, x_2, x_3, x_4) = (x_1/x_4, x_2/x_4, x_3/x_4)$$
which descends to a homeomorphism on $P(\cL)$. Therefore  $\gfl_\infty: T_1\Hs \rightarrow \chat$  satisfies 
$$\gfl_\infty = p\circ \hat\gfl_\infty.$$

We will need the following lemma.
\begin{lemma}
If $v, w \in T_{(x_1, x_2, x_3, x_4)}\cL$ then
$$\la v,w \ra = \frac{1}{x_4^2}\la p_*v,p_*w \ra_{\rm{euc}}.$$
\end{lemma}
 
 {\bf Proof:} For $u = (u_1, u_2, u_3, u_4) \in \R^{3,1}$  let $\hat u = (u_1, u_2, u_3) \in \R^3$.  Thus if $u \in T_x\cL$
 $$p_*(u) = \frac{1}{x_4}\left(\hat u -\frac{u_4}{x_4}\hat x\right).$$
As $v,w \in T_x\cL$ then $\la x,x\ra =\la v,x\ra = \la w,x\ra = 0$ giving 
$$\la \hat x, \hat x\ra_{\rm{euc}}= x^2_4,\qquad \la \hat v, \hat x\ra_{\rm{euc}}= v_4x_4, \qquad \la \hat w, \hat x\ra_{\rm{euc}}= w_4x_4.$$ 
Thus
\begin{align*}
 \la p_*v,p_*w \ra_{\rm{euc}} &= \frac{1}{x_4^2}\left\langle \hat v -\frac{v_4}{x_4}\hat x, \hat w -\frac{w_4}{x_4}\hat x\right\rangle_{\rm{euc}}\\
 &= \frac{1}{x_4^2}\left(\la \hat v, \hat w\ra_{\rm{euc}}-\frac{w_4}{x_4}\la \hat v, \hat x\ra_{\rm{euc}}-\frac{v_4}{x_4}\la \hat w, \hat x\ra_{\rm{euc}}+\frac{v_4w_4}{x_4^2} \la \hat x, \hat x\ra_{\rm{euc}}\right)\\
 & = \frac{1}{x_4^2}\left(\la \hat v, \hat w\ra_{\rm{euc}}-v_4w_4\right) = \frac{1}{x^2_4}\la v, w\ra
 \end{align*}
 \eproof

For an immersion $f$ the composition $f_\infty = \gfl_\infty \circ \hat f$ gives a map from $S$ to $\chat$ called the {\em Gauss map} of $f$. We now consider the pushforward of $\hat g$  by the Gauss map $f_\infty$ and show that it is conformal on $\hat\CC$ and can be described in terms of the visual metric.
 
 For each $p \in \Hs$ the sphere $T^1_p \Hs$ has an invariant metric (under isometries of $\Hs$ that fix $p$) that is unique up to scale. This is a spherical metric and we fix the scale so that the sphere has radius $1$. The {\em visual metric} $\nu_p$ based at $p$ is the push-forward of this measure via the hyperbolic Gauss map.

\begin{lemma}\label{proj_immersion}
 Let $f\colon S\to \Hs$ be an immersion with shape operator $B$ and let $\hat f$ be the lift to $T^1\Hs$. We also let $g = f^* g_{\Hs}$. If $f_\infty = \gfl_\infty \circ \hat f$ then $f_\infty$ is an immersion at $p \in S$ if and only if $-1$ is not an eigenvalue of $B$ at $p$.
 Furthermore at $p$ we have
$$\hat g = (\Id + B)^* g = (f_\infty)^*(\nu_{f(p)})$$
so $f_\infty$ is a conformal map from $\hat g$
to $\chat$.
\end{lemma}

{\bf Proof:} By definition of $f_\infty$ we have
$$f_\infty = \gfl_\infty \circ \hat f = p\circ\hat\gfl_\infty \circ \hat f.$$
Let
$\hat f: S\rightarrow T\R^{3,1}$  be given by $\hat f(q) = (f(q), n(q))$.  Then
$$(\hat\gfl_\infty)_*(v,w) = v+w \qquad (\hat f)_*(v) = (f_*(v), f_*(Bv))$$
and therefore
$$(\hat\gfl_\infty \circ \hat f)_*(v) = f_*(v) +f_*(Bv) = f_*((\Id + B)v).$$
By post-composition by an isometry, we can assume   that $f(q) = (0,0,0,1)$. Then $\nu_{f(q)}$ is the restriction of the euclidean metric to $\Sph^2$.
Therefore
\begin{align*}
\hat g(v,w) &= g((\Id +B)v,(\Id+B)w)\\
 & = \la f_*((\Id +B)v),f_*((\Id+B)w)\ra\\
& = \la p_*(f_*((\Id +B)v)),p_*(f_*((\Id +B)w)) \ra_{\rm{euc}}\\
& = \la (f_\infty)_*v,(f_\infty)_*w \ra_{\rm{euc}}\\
& = (f_\infty)^*(\nu_{f(q)})(v,w)
\end{align*}
\eproof

We can also set
$$\gfl_{-\infty} = \lim_{t\to -\infty}\pi\circ  \gfl_t$$
and let
$$f_{-\infty} = \gfl_{-\infty} \circ \hat f.$$
By reversing the orientation of $S$ and observing that this reverses the sign of the shape operator we see that $f_{-\infty}$ is conformal for the metric
$$\check g = (\Id-B)^* g.$$
We note that as $\hat g = (\Id+B)^* g$ then
$$\check g = \left((\Id+B)^{-1}(\Id-B)\right)^* \hat g = \hat B^* \hat g.$$

\subsection{The hyperbolic Gauss map and projective structures}
We have seen that given a conformal metric $\hat g$ on a simply connected projective structure $\Sigma = (S,f)$ we get a pair $(g,B)$ that is dual to $(\hat g, \hat B)$ where $\hat B$ is the projective shape operator of $\hat g$. Then by Theorem \ref{thm-bonnet} there is an immersion $f_0\colon S\to \Hs$ with $g = f_0^* g_\Hs$ and shape operator $B$. This is the {\em dual immersion} for the pair $(\hat g, \hat B)$. We can also consider it as the dual immersion for the pair $(\hat g, \Sigma)$.

The dual immersion $f_0$ also has a Gauss map $f_\infty\colon S\to \chat$. By Lemma \ref{proj_immersion}, $f_\infty$ is a
conformal map for $\hat g$ and hence an immersion. Therefore $\Sigma_{(g,B)} = (S, f_\infty)$ is a projective structure on $S$ that is in the same conformal class as $\Sigma$. We will show that $\Sigma = \Sigma_{(g,B)}$.

\begin{theorem}\label{equivalent_proj}
Let $\hat g$ be a conformal metric on a simply connected projective structure $\Sigma = (S,f)$. If $\hat B$ is the projective shape operator, let $(g, B)$ be the pair dual to $(\hat g, \hat B)$. Then $\Sigma_{(g,B)} = \Sigma$ and there is a unique choice for isometric immersion $f_0$ with associated Gauss map $f_\infty$ such that $f_\infty = f$.
\end{theorem}

{\bf Proof:}
The proof has three steps:
\begin{enumerate}
\item Let $\hat g_0$ be another conformal metric on $\Sigma$ that agrees with $\hat g$ on a open neighborhood. If $\hat B_0$ is its projective shape operator and $(g_0, B_0)$ the pair dual to $(\hat g_0, \hat B_0)$ then
 $\Sigma_{(g,B)} = \Sigma_{(g_0, B_0)}$.

\item If there is an open neighborhood where $\hat g$ is the restriction of the hyperbolic metric on a round disk then $\Sigma = \Sigma_{g,B}$.

\item For any metric $\hat g$ on $\Sigma$ there is another conformal metric $\hat g_0$ such that $\hat g$ and $\hat g_0$ agree on open neighborhood and there is an open neighborhood where $\hat g_0$ is the restriction of a hyperbolic metric on a round disk.
\end{enumerate}
Assuming this we can prove the theorem. 
 By (3) we have a conformal metric $\hat g_0$ on $\Sigma$ that agrees with $\hat g$ on a open neighborhood and on another open neighborhood is the restriction of the hyperbolic metric on a round disk. Therefore $(\hat g, \hat B)$ and $(\hat g_0, \hat B_0)$ agree on an open neighborhood so by (1) we have $\Sigma_{(g,B)} = \Sigma_{(g_0, B_0)}$. By (2) we have $\Sigma = \Sigma_{(g_0, B_0)}$. It  follows that $\Sigma = \Sigma_{(g,B)}$.\\
\\
We now prove the individual statements.
\\
{\bf Proof of (1):} By Lemma \ref{proj_immersion}, the metrics $\hat g$ and $\hat g_0$ are conformal on $\Sigma_{(g,B)}$ and $\Sigma_{(g_0,B_0)}$.
As $\hat g= \hat g_0$ on an open neighborhood $U$ we have that $(\hat g, \hat B) = (\hat g_0, \hat B_0)$ on $U$. It follows that $(g,B)= (g_0, B_0)$ on $U$. Therefore by Theorem \ref{thm-bonnet} the isometric immersions for $(g,B)$ and $(g_0, B_0)$ can be chosen to agree on $U$ and therefore so will their Gauss maps. The Gauss maps are conformal maps from $\hat g$ to $\chat$ so if they agree on an open neighborhood they agree everywhere. This implies that $\Sigma_{(g,B)} = \Sigma_{(g_0, B_0)}$.
\\
{\bf Proof of (2):} We first observe a domain $\Omega \subset \chat$ is a projective structure and if $\hat g_0$ is a conformal metric on $\Omega$ then there is a dual immersion $r\colon \Omega \to \Hs$ for the pair $(\Omega, \hat g_0)$. Furthermore if $U$ is an open neighborhood of $\Sigma = (S, f)$ and $\hat g = f^* \hat g_0$ on $U$ then we can choose the isometric immersion $f_0$ of $(\Sigma, \hat g)$ such that $f_0 = r\circ f$ on $U$. 

We also observe that if $\Omega$ is a round disk and $\hat g_0$ is the hyperbolic metric then $r\colon \Omega\to \Hs$ can be chosen such that Gauss map for $r$ is the identity. Our assumption is that there is an open $U$ where $\hat g = f^*\hat g_0$ (with $\hat g_0$ the hyperbolic metric on a round disk). Then on $U$ the Gauss map $f_\infty$ for $f_0$ is the the composition of $f$ with the Gauss map for $r$ so $f = f_\infty$ on $U$. Two conformal maps that agree on an open set agree everywhere so $f =f_\infty$ and $\Sigma = \Sigma_{(g,B)}$.
\\
{\bf Proof of (3):} Fix a point $p \in \Sigma$ and let $\psi\colon \Sigma \to [0,1]$ be a smooth function that is 1 in a neighborhood of $p$ and has support contained in a round disk $D$. If $\hat g_D$ is the hyperbolic metric on $D$ we let $\hat g_0 = \psi \cdot \hat g_D + (1-\psi)\cdot \hat g$. \\
\\
Finally we note that in (2) we showed that $f_0$ could be chosen such that $f=f_\infty$ in the case when on an open neighborhood $\hat g$ is the restriction of a hyperbolic metric on a round disk. In (1) we showed that the two Gauss maps can be chosen to agree. Together this implies that we can chose $f_0$ so that $f=f_\infty$ in the general case. \eproof

In what follows we will always choose the dual immersion $f_0$ such that its Gauss map $f_\infty$ has the property $f = f_\infty$.

We also have an equivariant version of Theorem \ref{equivalent_proj}:
\begin{cor}\label{equivariant_gauss}
Let $\hat g$ be a conformal metric on a simply connected projective structure $\Sigma = (S,f)$. If $\hat B$ is the projective shape operator, let $(g, B)$ be the pair dual to $(\hat g, \hat B)$. Let
$$f_0\colon S\to \Hs$$
be the unique isometric immersion for $(g, B)$ such that the $f_\infty = f$. If $\Gamma$ is a deck action on $S$, $\hat g$ is $\Gamma$-invariant and 
$$\rho\colon \Gamma \to \psl$$
a homomorphism with $f\circ \gamma = \rho(\gamma) \circ f$ for all $\gamma \in \Gamma$ then $f_0\circ \gamma = \rho(\gamma) \circ f_0$.
\end{cor}

\subsection{Epstein surfaces}
Note that if $\Omega$ is a domain in $\chat$ the $\Omega$ then the identity map gives $\Omega$ a projective structure.
 If $\hat g$ is conformal metric on $\Omega$ there then the pair $(\Omega, \hat g)$ determines a dual immersion. Epstein (see \cite{epstein-envelopes}) also gave a construction of an immersed surface in $\Hs$ associated to a conformal metric on $\Omega$.
We will show that the Epstein surface is equal to the dual immersion of the pair $(\Omega,\hat g)$.

In \cite{epstein-envelopes} Epstein restricts to conformal metrics $\hat g$ on domains $\Omega \subset \chat$. However, the construction works just as well if we have a metric $\hat g$ on a surface $S$ and a map $f\colon S\to \chat$ that is conformal with respect to $\hat g$. We will work in this setting. If $S$ was a domain in $\chat$ then $f$ would  be the inclusion map.

If $\hat g$ is a conformal metric at a point $z\in \chat$ then we let $\hh_{\hat g}(z)$ be the set of points $p$ in $\Hs$ where the visual metric $\nu_p$ agrees with $\hat g$ at $z$. Then $\hh_{\hat g}(z)$ is a horosphere. Given an immersion
$$f\colon S\to \chat$$
and a metric $\hat g$ such that $f$ is locally conformal on $(S, \hat g)$  Epstein constructs an immersion
$$f_{\ep}\colon S\to\Hs$$
such that $f = \gfl_\infty \circ \hat f_{\ep}$ and $(f_{\ep})_* T_z S$ is tangent to $\hh_{f^* \hat g}(f(z))$. This is the {\em Epstein surface} for $(f, \hat g)$.

As before we can multiply $\hat g$ by $e^{2t}$ to get a family $\hat g_t = e^{2t}\hat g$ of metrics and associated Epstein surfaces.

\medskip
\noindent
{\bf Theorem \ref{epstein_equal}} {\em
Let $\Sigma = (S,f)$ be a simply connected projective structure and $\hat g$ a conformal metric on $\Sigma$. Then $f_0 = f_{\ep}$.}
\medskip

{\bf Proof:} Note that we have chosen the isometric immersion $f_0$ such that $f = f_\infty$ where $f_\infty$ is the hyperbolic Gauss map for $f_0$. Then $(f_0)_*(T_p S)$ is tangent to a horosphere $\hh$ based at $f(p) = f_\infty(p)$. By Lemma \ref{proj_immersion} at $p$ we have $\hat g = (f_\infty)^*(\nu_{f(p)})$  which implies that $\hh = \hh_{f^*{\hat g}}$ so $f_0 = f_{\ep}$. \eproof

\section{Univalence, quasiconformal extension and convexity}
\subsection{Univalence and quasidisk criteria after Epstein}
Given a locally univalent map $f\colon \Delta \to \chat$, two classical  questions are what criteria on $f$ imply univalence and what criteria imply its image is a quasidisk. In a sequence of important papers (see \cite{epstein:shur,epstein:gaussmap,epstein:univalent,epstein-envelopes}),  Epstein introduced a new approach to these problems using the theory of immersed surfaces in $\Hs$ - the {\em Epstein surfaces} that we have been discussing. Following Epstein's basic approach we reprove and strengthen Epstein's results. These type of problems have a long history which we discuss at the end of this section.

\begin{theorem}\label{thm:univalence}
 Let $\hat g$ be a complete metric on $\Delta$, and $f\colon \Delta\to \chat$ a locally univalent map such that 
 $$\|Q(\Sigma_f,\hat g)\| \leq -\frac{1}{4}K(\hat g).$$ Then the map $f$ is univalent. Furthermore if $(g,B)$ is the pair dual to $(\hat g, \hat B)$
 there is a unique isometry $F\colon \cF_{(g,B)}(\R) \to \Hs$ which extends continuously to $\Delta\times\{\infty\}$ with $F(z,\infty) = f(z)$.  
 
 \end{theorem}
{\bf Proof:}
By Lemma \ref{Bevals} the eigenvalues of $\hat B = B(\Sigma_f,\hat g)$ are $-K(\hat g) \pm 4\|Q(\Sigma_f,\hat g)\|$ and there by hypothesis are non-negative. Thus  the eigenvalues of $B$ are in the interval $(-1, 1]$. By assumption $\hat g$ is complete. Therefore as $g = \frac{1}{4}(\Id+\hat B)^*\hat g$,  and $\hat B$ has non-negative eigenvalues then $g$ is also complete. Therefore by Proposition \ref{flow space}, there is an isometry 
$$F\colon \cF_{(g,B)}(\R) \to \Hs.$$
As a space $\cF_{(g,B)}(\R)$ is homeomorphic to $\Delta \times \R$ and  $F$ maps this space to $\Hs$, we can add $\Delta \times \{+\infty\}$ to the space and then continuously extend $F$ to a map from  $\Delta \times (-\infty,+\infty]$ to $\Hs\cup \chat$. The extension is the hyperbolic Gauss map for the immersion obtained by restricting $F$ to $\Delta \times \{0\}$. Therefore by Theorem \ref{equivalent_proj} we have $F(z, +\infty) = f(z)$. As map $F$ is injective on $\Delta \times \R$ and locally injective on $\Delta \times \{+\infty\}$, this implies that $f$ is injective, proving  univalence.
\eproof

{\bf Remark:} {\em We note that in the above theorem we only need that  the dual metric $g$ is complete but stating the theorem in this generality gives a more complicated statement as it would have conditions on both metrics $\hat g, g$.} \\
\\
The following is a more refined statement and has a more involved proof.

\begin{theorem}\label{nehari2}
Let $\hat g$ be a complete metric on $\Delta$, $f\colon \Delta\to \chat$ a locally univalent map and  $\hat B = B(\Sigma_f,\hat g)$ be the projective shape operator.
\begin{enumerate}
\item If 
$$\|Q(\Sigma_f,\hat g)\| \leq -\frac{1}{4}K(\hat g)$$
then $f$ is univalent and extends continuously to the boundary of $\Delta$. 
\item If 
$$\|Q(\Sigma_f,\hat g)\| < -\frac{1}{4}K(\hat g)$$
then $f$ extends to a homeomorphism of $\hat\CC$ and $f(\Delta)$ is a Jordan domain.
\item If $K(\hat g) < 0$ and there is a $0\leq k < 1$ such that
$$-\frac{\|Q(\Sigma_f,\hat g)\|}{4K(\hat g)} \leq k < 1.$$
Then $f$ extends to a quasiconformal homeomorphism  with beltrami differential $\|\mu\|_\infty \leq k$. 
\end{enumerate}
\end{theorem}

The are a number of natural approaches  to the proof of this theorem. One   involves showing that the embedding $f$ extends continuously to the boundary by considering the limiting behavior of curves in $\Hs$ of curvature $\le 1$. Instead, we follow a more classical approach using Carath\'eodory's theory of prime ends which we now briefly describe.

\subsubsection{Extension of univalent maps, prime ends}
The following theorem due to Marie Torhorst describes precisely when a univalent map of the disk extends continuously to its boundary.

\begin{theorem}[Torhorst, \cite{uni_ext}, see also \cite{Pommerenke:book:bdry}]\label{tor}
Let $f:\Delta \rightarrow U \subseteq \hat\CC$ be a univalent map. Then $f$ extends continuously to $\partial \Delta$ if and only if $\partial U$ is locally connected.
\end{theorem}

This theorem was proven using the theory of prime ends developed by Caratheodory (see \cite{Pommerenke:book:bdry} for details). We now briefly describe this.

We let $U \subseteq \hat\CC$ be a simply connected domain in $\hat\CC$. A {\em crosscut} $C$ of $U$ is a simple arc in $U$ with $\overline C = C \cup \{a,b\}$ for some $a,b \in \partial U$ (possibly equal).

A {\em null chain} of $U$ is a sequence $\{C_n\}_{n=0}^\infty$ of  crosscuts of $U$ such that
\begin{itemize}
\item $\overline C_n \cap \overline C_{n+1} = \emptyset$
\item $C_n$ separates $C_0$ from $C_{n+1}$
\item $\diam(C_n)\rightarrow 0$
\end{itemize}
If $C_n$ is a null-chain then we define $V_n$ to be the component of $U-C_n$ not containing $C_0$. Thus $V_n$ are a decreasing collection. Two null-chains $C_n, D_n$  are equivalent if for sufficiently large $m$ there exists an $n$ such that
$$\mbox{$D_m$ separates  $C_n$ from $C_0$}\qquad \mbox{$C_m$ separates  $D_n$ from $D_0$}.$$
The equivalence class of null-chains is called a {\em prime end} of $U$. Given a prime end $p = [C_n]$ then the {\em impression} of $p$ is the compact set
$$I(p) = \bigcap_{n=1}^\infty \overline V_n.$$
This is well-defined. If $I(p)$ is a single point then $p$ is called {\em degenerate}. 

Using the theory of prime ends,  we have the following classical description of univalent maps with continuous extension due to Caratheodory and others.

\begin{theorem}[see \cite{Pommerenke:book:bdry}]\label{prime}
Let $f:\Delta \rightarrow U$ be a univalent map. Then $f$ extends continuously to $\partial \Delta$ if and only if all prime ends of $U$ are degenerate.
\end{theorem}

It follows from the above theorems that $\partial U$ being locally connected is equivalent to $U$ having only degenerate prime ends and it was this that Torhorst actually proved.

\subsubsection{Continuous extension and quasidisk criterion}
We let $\hat g$ be a complete metric on $\Delta$, $f\colon \Delta\to \chat$  a locally-univalent  map with   projective shape operator $\hat B = B(\Sigma,\hat g)$ having non-negative eigenvalues. Then we have $f$ is univalent with image $U =f(\Delta)$ and there is an isometry $F\colon \cF_{(g,B)}(\R) \to \Hs$ which extends continuously to $f$ on $\Delta\times\{\infty\}$. 

We let $G(\Hs)$ be the space of oriented geodesics in $\Hs$ and identify it with 
$$G(\Hs) =  \hat\CC^2 -\{(z,z)\ |\ z\in \hat\CC\}.$$ The space $G(\Hs)$ is naturally compactified by $\hat\CC^2$. 

Given map $F\colon \cF_{(g,B)}(\R) \to \Hs$ we define map $\hat F: \Delta \rightarrow G(\Hs)$ by letting $\hat F(z)$ be the oriented geodesic $F(\{z\}\times\R)$. In particular $\hat F(z) = (f(z), h(z))$ for some $h:\Delta\rightarrow \hat\CC$ continuous.

\begin{lemma}\label{geod}
The map $\hat F:\Delta \rightarrow G(\Hs)$ is proper and a homeomorphism onto its image. Furthermore $\hat F(z_i) \rightarrow (p,p)$ if and only if $f(z_i) \rightarrow p \in \partial U$.
\end{lemma}
 {\bf Proof:} we first prove properness. Let $\hat F(z_i) \rightarrow L$ and choose $y \in L$. Then as $F$ is a homeomorphism, there is a unique $x= (z,t)$ with $F(x) = y$ and therefore   $y \in \hat F(z)$. We choose a neighborhood $W \subseteq \Delta$ of $z$, then $F$ maps neighborhood $W\times \R$ of $x$ homeomorphically to a neighborhood of $y$. Thus geodesic $\hat F(z_i) = F(\{z_i\}\times \R)$ intersects $F(W\times \R)$ for $i$ large. As $F$ is a homeomorphism this implies $z_i \in W$. 
Thus $z_i \rightarrow z$ and by continuity $L = \hat F(z)$.  Thus $\hat F$ is proper and therefore a homeomorphism onto its image. The second item follows also.  \eproof

 We use the above to prove continuous extension.
 
 \begin{theorem}
  Let $\hat g$ be a complete metric on $\Delta$, and $f\colon \Delta\to \chat$ a locally univalent map such that projective shape operator $\hat B = B(\Sigma_f,\hat g)$ has non-negative eigenvalues. Then the map $f$ is univalent and $f$ extends continuously to its boundary.
  \end{theorem}
  
  {\bf Proof:} We have already proven $f$ is univalent and by Caratheodory, it suffices to prove that the prime ends of $\partial U$ are degenerate.
  
 We let $F\colon \cF_{(g,B)}(\R) \to \Hs$ be the isometry which extends continuously to $f$ on $\Delta\times\{\infty\}$. Let $p = [C_n]$ be a prime end of $U$. The crosscut $C_n$ can be extended to a disk $E_n$ in $\Hs$. Explicitly we let $B_n = f^{-1}(C_n)$ and $E_n = F(B_n\times\R)$. Thus $E_n$ is a disk in $\Hs$ with $\partial E_n \subseteq \hat\CC$ and $\partial E_n\cap U = C_n$. We consider $z_n \in B_n$ such that (after reducing to a subsequence) $\hat F(z_n)$ converges. If $\hat F(z_n) \rightarrow L$ then by Lemma \ref{geod} we have  $z_n\rightarrow z \in \Delta$ with $L = \hat F(z)$. Therefore by continuity $f(z_n) \rightarrow f(z) \in U$. As $f(z_n) \in C_n$ then $f(z) \in \partial U$ a contradiction. Thus $\hat F(z_n) \rightarrow (w,w)$ with $w \in \partial U$ and therefore $f(z_n) \rightarrow w$. It follows that as $\diam(C_n) \rightarrow 0$ then $\diam(E_n) \rightarrow 0$ where the diameter for $E_n$ is with respect to the Euclidean metric on the unit ball. We further have $\partial E_n$ consists of $C_n$ with endpoints $\partial C_n = \{a_n,b_n\}$ joined by a continuous arc in the complement of $U$. Thus $\partial E_n$ separates $V_n$ from $C_0$. It follows that  $\diam(V_n) \rightarrow 0$ and therefore the impression $I(p) = \cap_{n=1}^\infty \overline{V_n}$ is a single point. Thus all the prime ends of $U$ are degenerate.
 It follows from Theorem \ref{prime} that $f$ extends continuously to its boundary.
 \eproof

 We now prove our quasidisk criterion. We  will need to use a converse of the Jordan curve theorem due to Schoenflies. If $U$ is open in $\hat\CC$ then a set $K$ is {\em accessible} from $U$ if for every point $p \in K$ there is a simple arc $\alpha:[0,1]\rightarrow \hat\CC$ such that $\alpha(1) = p$ and $\alpha(t) \in U$ for $t \in [0,1)$. Clearly $K \subseteq \partial U$ but in general it may not be true that $\partial U$ is accessible from $U$ although it is true for Jordan domains. We have the following;

\begin{theorem}[{Schoenflies, \cite[Theorem VI.16.1]{newman}}]\label{conversejordan}
Let  $K \subseteq \hat\CC$ be closed with complement two disjoint open simply connected sets $U,V$ such that  $K$ is accessible from both $U$ and $V$. Then  $K$ is a Jordan curve. 
\end{theorem}

We apply this to show that when the shape operator has eigenvalues in $(-1,1)$ then the image of $\partial\Delta$ is a Jordan curve.

 \begin{theorem}\label{thm:jordan}
 Let $\hat g$ be a complete metric on $\Delta$, and $f\colon \Delta\to \chat$ a locally univalent map such that  $\hat B = B(\Sigma_f,\hat g)$ have positive eigenvalues. Then $U$ is a Jordan domain and $f$ extends  to a homeomorphism on $\overline\Delta$.
  \end{theorem}

{\bf Proof:} As $\hat B$ has positive eigenvalues then $F$ extends to a map $\Delta\times[-\infty,\infty]\rightarrow \Hs\cup\hat\CC$ with $F(z,\infty) = f(z)$ and $F(z,-\infty) = h(z)$ where $f$ is conformal with respect to $\hat g$ and $h$ is conformal with respect to $\check g$. Further $U = f(\Delta)$ and $V = h(\Delta)$ are disjoint as if $f(z) = h(w)$ then the geodesic $\hat F(z)$ and $\hat F(w)$ are asymptotic. Then we can choose sequence $z_i \rightarrow z$ such that  geodesics $\hat F(w),\hat F(z_i)$  intersect. But $\hat F(z)$ foliate $\Hs$. Thus $z = w$.

Let $p \in (U\cup V)^c$. Then we can find $w_i \in \Hs$ such that $w_i \rightarrow p$. We let $w_i \in \hat F(z_i)$. If (after reducing to a subsequence) $\hat F(z_i) \rightarrow L$ then by Lemma \ref{geod} $z_i \rightarrow z$ with $L = \hat F(z)$. As $w_i$ is on geodesic  $\hat F(z_i)$ then $p$ is an endpoint of $L$. Thus $p = f(z)$ or $p=h(z)$ contradicting $p \in (U\cup V)^c$. Thus $\hat F(z_i) \rightarrow (p,p)$ and by Lemma \ref{geod} $p \in \partial U$ and similarly $p \in\partial V$. Thus $\partial U = \partial V = (U\cup V)^c$ and $U\cup V \cup \partial U$ is a disjoint decomposition. It follows that $\partial U$ is closed and as both maps $f, h$ extend continuously to $\partial \Delta$, then $\partial U=\partial V$ is accessible from both $U$ and $V$. Thus by  Schoenflies, converse of the Jordan curve theorem (Theorem \ref{conversejordan}), $\partial U$ is a Jordan curve.
\eproof

We can now prove the more refined univalence result.
\\

{\bf Proof of Theorem \ref{nehari2}:} 
By Lemma \ref{Bevals} the eigenvalues of $\hat B$ are $-K(\hat g) \pm 4\|Q(\Sigma_f,\hat g)\|$. Therefore   $\hat B$ having non-negative eigenvalues is equivalent to item 1. Therefore $f$ is univalent as before.

If $\hat B$ satisfies items 2 or 3, then $\hat B$ has positive eigenvalues. Therefore we have two hyperbolic Gauss maps $f:\Delta\rightarrow U$ and $h:\Delta \rightarrow V$ where $f,h$ are conformal with respect to $\hat g, \check g$ respectively.  Further by Theorem \ref{thm:jordan} $U,V$ are disjoint with $\partial U = \partial V = \partial S$  a Jordan curve. Furthermore $f, h$  extend to homeomorphisms on the closed disk $\overline\Delta$. To show they agree on $\partial \Delta$ let $p \in \partial \Delta$ and $f(p) = a, h(p)=  b$. We choose $z_i \rightarrow p$ and  $\hat F(z_i)\rightarrow (a,b)$. As $f(z_i) \rightarrow a\in \partial U$, then by Lemma \ref{geod}  $a = b$ and the maps agree on $\partial \Delta$.

Thus combining the two extensions, we obtain an extension $f_{\ext}$ of $f$ to $\hat\CC$. Explicitly we let $f_{\ext}(z) = f(z)$ for $z\in\overline{\Delta}$ and $f_{\ext}(z) = h(1/\overline{z})$ for $z \in \Delta^* = \{z\in \hat\CC\  |\ |z| > 1\}$. By the above $f_{\ext}$ is a homeomorphism extending $f$. 

Now  if (3) holds we show that $f_{\ext}$ is a quasiconformal map. A standard property of quasiconformal mappings is if $U,V \subseteq \hat\CC$ are open and $q:U \rightarrow V$ is a homeomorphism which is quasiconformal on $U-\ell$ for some line $\ell$ then $q$ is a quasiconformal homeomorphism on $U$ (see \cite[Proposition 4.27]{Hubbard:book:TeichmullerI}). Thus as $\partial \Delta$ is locally a line, it suffices to show $f_{\ext}$ is a quasiconformal  on $\Delta \cup \Delta^*$. 

This follows as $F$ is conformal with respect to $\hat g$ on $\Delta\times\{\infty\}$ and conformal with respect to $\check g = \hat B^*\hat g$ on $\Delta\times\{-\infty\}$. Therefore as $\hat B$ has eigenvalues $-K(\hat g) \pm 4\|Q(\Sigma_f,\hat g)\|$ then $ f_{\ext}$ is $k$-quasiconformal with
$$K = \frac{-K(\hat g) +4\|Q(\Sigma_f,\hat g)\|}{-K(\hat g) -4\|Q(\Sigma_f,\hat g)\|} .$$
Thus $f_{\ext}$ has beltrami differential $\mu$ with
$$|\mu| = \frac{K-1}{K+1} = \frac{4\|Q(\Sigma_f,\hat g)\|}{-K(\hat g)} \leq k < 1.$$
\eproof

\subsection{Quasi-conformal reflections}
A Jordan curve $\gamma$ admits a {\em quasiconformal reflection} if there is an orientation reversing quasiconformal homeomorphism $h:\hat\CC \rightarrow \hat\CC$ fixing $\gamma$ and interchanging the two components of its complement. In \cite{Ahlfors:quasicircles}, Ahlfors gave a characterization for which Jordan curves admit a quasiconformal reflection. Epstein noted (\cite{epstein:gaussmap}) that his quasidisk criterion also gave a criterion for a Jordan curve to admit a quasiconformal reflection. We make the same observation in the following immediate corollary to Theorem \ref{nehari2}.
 
\begin{cor}\label{qc_reflection}
Let $\hat g$ be a complete metric on $\Delta$, $f\colon \Delta\to \chat$ a locally univalent map and  $\hat B = B(\Sigma_f,\hat g)$ be the projective shape operator. If $\hat g$ is negatively curved and there is a $0\leq k < 1$ such that
$$\|Q(\Sigma_f,\hat g)\| \leq -\frac{k}{4} K(\hat g)$$
then $U = f(\Delta)$ is a Jordan domain such that $\partial U$  admits a quasiconformal reflection with beltrami differential $\|\mu\|_\infty \leq k$. 
\end{cor} 
{\bf Proof:} We keep the same notation in case 3 of Theorem \ref{nehari2} above. Then we let $H:\hat\CC \rightarrow \hat\CC$ be given by $H = f_{\ext}\circ C \circ f_{\ext}^{-1}$ where $C(z) = 1/\overline{z}$.  Trivially $H(z) = z$ on $\partial U$ and as one of either $f_{\ext}$ or $f_{\ext}^{-1}$ is conformal in the composition at any point in the complement of $\partial U$, we have the beltrami differential bound is the same as for $f_{\ext}$.
\eproof

\subsubsection{History}

 The study of these two problems has a long history. Two foundational results are the Nehari univalence criterion that  if $\|Sf\|_\infty \leq \frac{1}{2}$ then $f$ is univalent (see \cite{Nehari:schwarzian}) and the Ahlfors-Weill Extension Theorem that if $\|Sf\|_\infty  \leq \frac{k}{2} < \frac{1}{2}$, then $f$ extends to a $K= \frac{1+k}{1-k}$ quasiconformal homeomorphism (see \cite{ahlforsweill}).  These results can be obtained by letting $\hat g$ be they hyperbolic metric in Theorems \ref{thm:univalence} and \ref{nehari2}.
 
 In \cite{epstein:gaussmap} Epstein proved Theorem \ref{thm:univalence} under two extra conditions. Namely he assumed that $\hat g = e^{2\phi}g_{\Hp}$ is negatively curved (rather than just non-positively curved) and that
  $$ |\phi_z(w)|(1-|w|^2) \leq c\max\left(|w|, \frac{1}{2|w|}\right).$$
In \cite{Pommerenke_Epstein}, Pommerenke replaced this later condition with the weaker assumption that  $\phi = \log |h'|$ where $h:\Delta \rightarrow \CC$ is conformal. We also note that these two authors did not state their conditions in terms of the Osgood-Stowe differential so their expressions appear different (and significantly longer).

In \cite{OS_nehari}, Osgood-Stowe were the first to state a univalence criteria in terms of their differential and their criteria also applied in higher dimensions. In our setting they showed that if $\hat g$ is a geodesically convex conformal metric on $\Delta$ with
\begin{equation}\|Q(\Sigma_f,\hat g)\|_\infty \leq -\frac{1}{4}K\label{os_cond}\end{equation}
then the map is univalent. A crucial difference is that they do not require the metric $\hat g$ to be complete (for example it could be the Euclidean metric on the disk) and they get better estimates when the diameter is finite. However, the sup-norm bound is a stronger assumption than the pointwise bound in Theorem \ref{thm:univalence}. Stronger statements of this type were later proved by Chuaqui (see \cite{Chuaqui_nehari}).

As we noted above when $\hat g$ is the hyperbolic metric (3) of Theorem \ref{nehari2} is the classical result of Ahlfors-Weill (\cite{ahlforsweill}). When $\hat g$ is is hyperbolic (2) was proven by Pommerenke-Gehring (\cite{PGnehari}). Under the extra technical assumptions on the metrics described above, Theorem \ref{nehari2} was proved by Epstein in \cite{epstein:gaussmap}.

 As well as  generalizing Epstein's criteria, our approach is short (sections 4 to 6 of this paper), proving the main results of \cite{epstein:gaussmap,epstein:univalent,epstein-envelopes} and does not require the use of Epstein's unpublished  paper \cite{epstein:shur}  where he proved a version of Shur's Lemma for hyperbolic space. Also as  Epstein's  first paper \cite{epstein-envelopes} (where he introduced  Epstein surfaces and proved many of their properties) remains unpublished, we hope that this paper will bring the results in this important paper to a wider audience.

\subsection{Example $f(z)= z^c$} 
\begin{figure}[htbp] 
   \centering
   \includegraphics[width=2.5in]{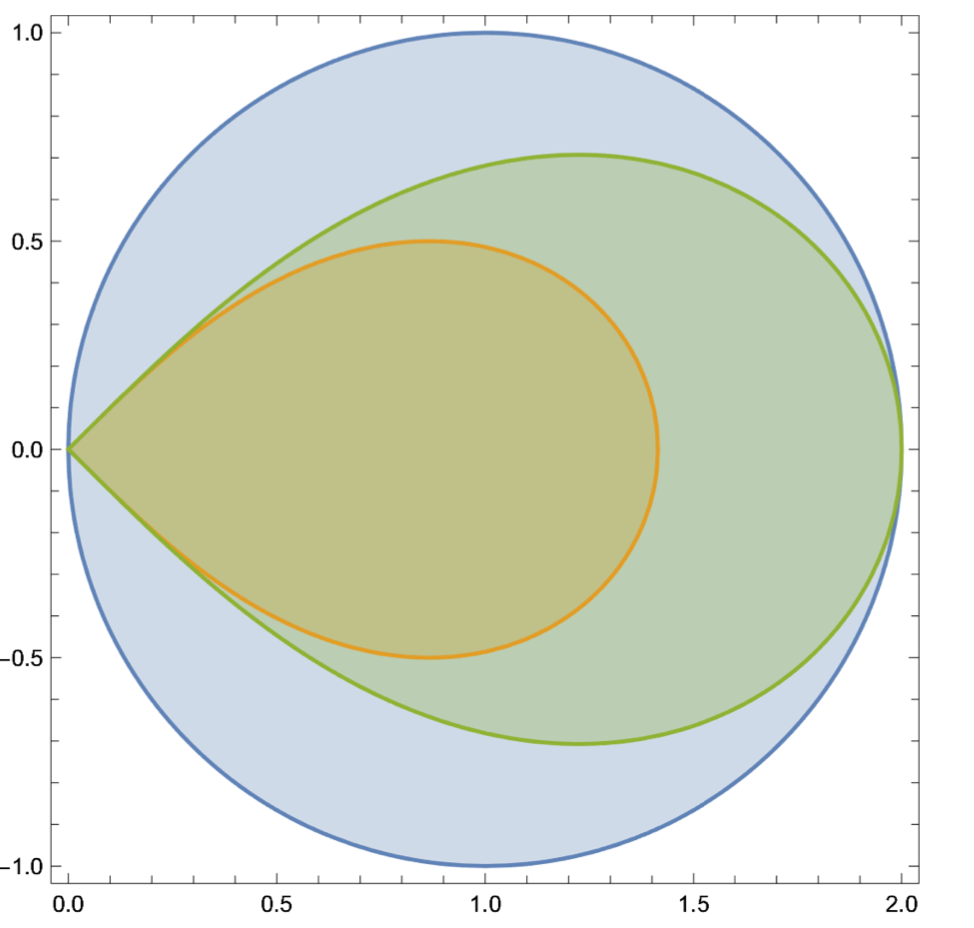} 
   \caption{Comparing univalence conditions for $f(z)=z^c$}
   \label{univalent}
\end{figure}
We now consider the above univalence condition for the projective $\Sigma = (f,\Hp)$  given by $f(z) = z^c$ on the upper half plane. We choose a general conformal metric $\hat g = e^{2\phi}\geu$ which is invariant under the action of $\R$  of $t\cdot z = e^tz$. Then in polar coordinates 
$$\phi  = -\log(r) + h(\theta)$$
for some function $h(\theta)$. In local coordinates the univalence condition is
$$|Q(\Sigma,\hat g)| \leq \frac{1}{4}\Delta \phi. $$ 
After computing we get
$$Q(\Sigma,\hat g)  = \frac{1}{4z^2}(h'(\theta)^2 -  h''(\theta) +c^2)$$
and the univalence criterion becomes
$$|h'(\theta)^2 +c^2-  h''(\theta)| \leq h''(\theta).$$
Letting  $h(\theta) = -\log(\sin(\theta))$ gives the hyperbolic metric and therefore we obtain Nehari's condition $||S(f)||_\infty < 1/2$ or equivalently $|c^2-1| \leq 1$. 

We now let $h(\theta) = -2\log(\sin(\theta))$. Then we have
$$|4\mbox{cotan}^2(\theta) +c^2 -2\mbox{cosec}^2(\theta)| \leq 2\mbox{cosec}^2(\theta) .$$
After simplification, this reduces to 
$$|c^2-2| \leq 2.$$ To compare, this gives a better univalence condition than the Nehari inequality   $|c^2-1| \leq 1$ and is a closer approximation to  the full set of univalent maps of the form $z^c$ which is the set $|c-1| \leq 1$ (see Figure \ref{univalent}).

\section{The hyperbolic extension of a projective structure}
If a projective structure $\Sigma$ is not simply connected the dual immersions will not map to $\Hs$. Instead we replace $\Hs$ with $\ext(\Sigma)$, the {\em hyperbolic extension} of $\Sigma$.

We first define this when $\Sigma$ is simply connected. In the general case we will take a quotient.

Let the pair $\Sigma = (S, f)$ define a projective structure with $f\colon S\to \chat$ and immersion. Then $U \subset S$ is a {\em round disk} in $\Sigma$ if $f(U)$ is a round disk. For each round disk $U$ let $\bar{H}_U = U \times [0,1]$ and $H_U = U \times [0,1)$. We then fix a continuous map $F_U\colon \bar{U}_H \to \to \Hs\cup\chat$ such that $F_U(z,1) = f(z)$ and $F_U$ restricted to $H_U$ is a homeomorphism to a closed hyperbolic half-space in $\Hs$. If $U_0$ and $U_1$ are round disks in $\Sigma$ with $U_0 \cap U_1 \neq \emptyset$ then we define $(z_0, t_0) \sim (z_1, t_1)$ if $F_{U_0}(z_0, t_0)  = F_{U_1}(z_1,t_1)$. This defines an equivalence relation of the disjoint union of $\bar{H}_U$ and the {\em projective extension} $\ext(\Sigma)$ of $\Sigma$ is the quotient of this union by the equivalence relation $\sim$. Similarly we define $\ext^{\rm o}(\Sigma)$ as the subspace of $\ext(\Sigma)$ given by the equivalence relation on the disjoint union of $H_U$. Then $\ext(\Sigma)$ is homeomorphic to $S \times [0,1]$ and there is a  map $F\colon \ext(\Sigma) \to \Hs \cup \chat$ such that on each $\bar{U}_H$ we have $F = F_U$. In particular $F(z,1) = f(z)$. Note that $F$ pulls back a hyperbolic metric to
 $\ext^{\rm o}(\Sigma)$ but the boundary is not smooth. Instead $S \times \{0\}$ is a locally convex {\em pleated surface}. We also note that while if $(z_0, 1) \sim (z_1, t_1)$ then $t_1 = 1$ in general the equivalence relation does not respect the product structure on the $\bar{U}_H$. For further details see \cite{Kamishima:Tan}.

If $\Sigma$ is a projective structure on a surface $S$ that does not come from a pair $(S,f)$ then we can lift $\Sigma$ to a projective structure $\tilde\Sigma$ on the universal cover $\tilde S$ of $S$. Then $\tilde\Sigma$ does come from a pair $(\tilde S, f)$ and furthermore the deck action on $\tilde\Sigma$ will extend to a deck action on $\ext(\tilde\Sigma)$ where the action on the hyperbolic part of $\ext(\tilde\Sigma)$ will be by isometries. Then the quotient of $\ext(\tilde\Sigma)$ by this deck action will be $\ext(\Sigma)$.

\subsection{Embedding the flow space in $\ext(\Sigma)$}
We consider conformal metrics $\hat g$  on $\Sigma$ such that the dual pair $(g,B)$ locally convex surfaces. We show that then the associated flow space $\cF_{(g,B)}( [0,\infty))$ isometrically embeds in $\ext^{\rm o}(\Sigma)$.  Further using this, we will show that for two such metrics $\hat g_1,\hat g_2$ with $\hat g_1 \leq \hat g_2$ then   $\cF_{(g_2,B_2)}( [0,\infty))$ isometrically embeds in $\cF_{(g_1,B_1)}( [0,\infty))$. 

\begin{prop}\label{flow_embeds}
Let $\hat g$ be a  conformal metric on a projective structure $\Sigma$ on a surface $S$ with projective shape operator $\hat B$ and dual pair  $(g,B)$. If $g$ is complete and the eigenvalues $\hat B$ are between $(-1,1]$ the flow space $\cF_{(g,B)}( [0,\infty))$ isometrically embeds in $\ext^{\rm o}(\Sigma)$.

The intersection of $\Sigma \times \{0\}$ with the internal boundary of $\ext(\Sigma)$ is a union of plaques. 
\end{prop}

{\bf Proof:} First assume that $\Sigma$ is simply connected. By Lemma \ref{flow space} since $g$ is complete the flow space $\cF_{(g,B)}([0,\infty))$ is complete and any piece of a hyperbolic plane tangent to $\Sigma\times \{0\}$ will extend to an immersed copy of $\htwo$. The local isometry from the flow space to $\Hs$ this plane will be embedded and hence it must be embedded in the flow space. In particular it will extend to a round disk $U$ in $\Sigma$. One can then construct the flow space similarly to the construction of $\ext(\Sigma)$ however rather than take all round disks in $\Sigma$ we just take those that arise from tangent planes to $\Sigma \times \{0\}$. Then the resulting quotient space will be a subspace of $\ext(\Sigma)$. On the other hand this quotient space will be isometric to the original flow space $\cF_{(g,B)}([0,\infty))$ so the flow space isometrically embeds in $\ext^{\rm o}(\Sigma)$.

Note that the construction of $\ext(\Sigma)$ and the flow space is canonical so any deck action that preserves both the projective structure and the metric $\hat g$ will extends to isometries of $\ext(\Sigma)$ and the flow space. In particular the quotient of the flow space (which will be the flow space of the quotient) will isometrically embed in the quotient of $\ext(\Sigma)$ (which will be the extension of the quotient). Therefore to prove the proposition when $\Sigma$ isn't simply connected we lift to the universal cover, construct the embedding there and then take the quotient.
 \eproof

There is a natural metric on $\Sigma$ closely associated to $\ext(\Sigma)$. As before assume that $\Sigma$ is simply connected. If $U$ is a round disk let $\hat g_U$ be the hyperbolic metric. Then Thurston defines the {\em projective metric} $\hat g_\Sigma$  to be the infimum, at each point $z\in \Sigma$ of $\hat g_U$ at $z$ where $U$ varies over all round disks that contain $z$. If $\Sigma$ is hyperbolic then the infimum is realized and $\hat g_\Sigma$ is everywhere positive. We note that while this metric will be continuous, it will not be smooth. If $\Sigma$ is not simply connected we observe that the projective metric on the universal cover will be equivariant and descend to a metric on $\Sigma$. This is the projective metric $\hat g_{\Sigma}$ on $\Sigma$.

An important special cases is when $\Sigma$ is an open domain in $\chat$. That is $\Omega = \chat \smallsetminus \Lambda$ for some closed set $\Lambda$ in $\chat$. The convex hull of $\Lambda$ in $\Hs$ is the smallest, closed convex set whose closure in $\chat  = \del\Hs$ is $\Lambda$. The convex hull is the intersection of all half spaces in $\Hs$ whose closure in $\chat$ contains $\Lambda$. Then by definition $\ext(\Sigma)$ is the union of the complementary half-spaces. It follows that $\ext(\Sigma)$ is the complement of the interior of the convex hull.

Let $\Sigma$ be a projective structure on a simply connected surface $S$ determined by an immersion $f_\infty\colon S\to \chat$. We can assume that the restriction of $F\colon \ext(\Sigma) \to \Hs\cup \chat$ to $\Sigma$ is $f_\infty$. If $\hh$ is a horoball in $\ext(\Sigma)$ based at $z \in \Sigma$ then $F(\hh)$ is a horoball in $\Hs$ based at $w = f_\infty(z)$. Then $F(\hh) = \hh_{g_w}(w)$ where $g_w$ is a conformal metric for the tangent space $T_w\chat$. If $\hat g$ is conformal metric on $\Sigma$ and $\hat g = (f_\infty)^* g_w$ at $z$ then we write $\hh_{\hat g}(z) = \hh$.

The next lemma gives a condition on the conformal metric $\hat g$ for the horoballs $\hh_{\hat g}(z)$ to embed in $\ext(\Sigma)$.
\begin{lemma}
If $\hat g$ is a conformal metric on a simply connected projective structure $\Sigma$ with $\hat g \ge g_\Sigma$ then the horoball $\hh_{\hat g}(z)$ is contained in $\ext(\Sigma)$ for all $z \in \Sigma$.
\end{lemma}

{\bf Proof:} By the definition of the projective metric for any $z \in \Sigma$ there is a round disk $U$ in $\Sigma$ such that $g_U \le \hat g$ at $z$. Then the half space $H_U$ will isometrically embed in $\ext(\Sigma)$ and the horoball $\hh_{g_U}(z)$ will be contained in $H_U$ and therefore will be contained in $\ext(\Sigma)$. As $g_U \le \hat g$ at $z$ we also have that $\hh_{\hat g}(z)$ is contained in $\hh_{g_U}(z)$ so $\hh_{\hat g}(z)$ is contained in $\ext(\Sigma)$. \eproof

\begin{prop}\label{mono-flow}
Let $\hat g_1$ and $\hat g_2$ be conformal metrics on a simply connected projective structure and assume that both projective shape operators $\hat B_1, \hat B_2$ have eigenvalues in the interval $(-1,1]$. If $\hat g_1 \le \hat g_2$ then the flow space $\cF_{(g_2, B_2)}([0,\infty))$ isometrically embeds in the flow space $\cF_{(g_1, B_1)}([0,\infty))$.
\end{prop}

{\bf Proof:} By Proposition \ref{flow_embeds} the flow spaces $\cF_{(g_i,B_i)}([0,\infty))$ embed in $\ext(\Sigma)$ and  is the union of horoballs $\hh_{\hat g_i}(z)$. Since for each $z \in \Sigma$ the horoballs $\hh_{\hat g_2}(z)$ are contained in the horoballs $\hh_{\hat g_1}(z)$ this implies that the flow space $\cF_{(g_2, B_2)}([0,\infty))$ isometrically embeds in the flow space $\cF_{(g_1, B_1)}([0,\infty))$. \eproof

\subsection{Convexity}
We derive a convexity criterion for the dual surface of a pair $(\hat g,\hat B)$. 
\begin{theorem}\label{loc_conv}
Let $(\hat g, \hat B)$ be  dual to $(g, B)$ with $\hat B = B(\Sigma,\hat g)$. Let $f$ be an  immersion into $\Hs$ with fundamental pair $(g,B)$ and $f_t = \pi\circ \gfl_t \circ \hat f$. Then for $$e^{2t_0} > \sup\{|K(\hat g)| + 4\|Q(\Sigma,\hat g)\|\}$$  $f_{t_0}$ is a locally convex immersion. In particular if $g_\Sigma$ is the projective metric on $\Sigma$ and $\hat g$ is complete then
$$g_\Sigma \leq \sup\{|K(\hat g)| + 4\|Q(\Sigma,\hat g)\|\} \hat g.$$
\end{theorem}

{\bf Proof:}  As $\hat B_t = e^{-2t}\hat B$ and $\hat B = B(\Sigma,\hat g)$ has eigenvalues $-K(\hat g) \pm 4\|Q(\Sigma,\hat g)\|$ then  for 
$$e^{2t} > \sup\{ |K(\hat g)| + 4\|Q(\Sigma,\hat g)\|\}$$
we have $\hat B_t$ has eigenvalues in $(-1,1)$. Thus by duality, the eigenvalues of $B$ are in $(0,\infty)$. It follows that $f_t$ has positive principal curvatures and therefore is a locally convex immersion.

 By \cite{BBB}, if $\hat g$ is  complete and pair $(\hat g, \hat B)$  has locally convex immersed dual surface, then $g_\Sigma \leq \hat g$ where 
 $g_\Sigma$ is the projective metric on $\Sigma$.  Therefore for $t > \sup\{|K(\hat g)| + 4\|Q(\Sigma,\hat g)\|\}$ we have
$g_\Sigma \leq e^{2t}\hat g .$ Taking the infimum over all such $t$ we get
$$g_\Sigma \leq \sup\{|K(\hat g)| + 4\|Q(\Sigma,\hat g)\|\} \hat g$$
\eproof

In \cite[Theorem 2.8]{BBB} we prove that for $g_h$ the hyperbolic metric then
$$g_{\Sigma} \leq (1+2\|\phi_\Sigma\|_\infty) g_h.$$
As $Q(\Sigma, g_h) = -\phi_\Sigma/2$  we see that the above theorem is a generalization for arbitrary conformal metrics.

\section{Hyperbolic 3-manifolds}
Let $M$ be a hyperbolic 3-manifold. Then $M = \Hs/\Gamma$ where the deck group $\Gamma$ is a discrete group of isometries of $\Hs$. The group $\Gamma$ will also act on $\partial\Hs = \chat$ with domain of discontinuity $\Omega$. Then the quotient $\bar M = (\Hs \cup \Omega)/\Gamma$ is a 3-manifold with boundary. If $\bar M$ is compact then $M$ is a {\em conformally compact} hyperbolic 3-manifold.  The boundary $\Sigma = \Omega/\Gamma$ is a projective structure and is the {\em projective boundary} of $M$. While conformally compact is the usual terminology for Einstein manifolds, for hyperbolic 3-manifolds one usually refers to $\bar M$ as {\em convex co-compact} because, as we will see shortly, there are convex $\Gamma$-invariant subsets of $\Hs$ with compact $\Gamma$-quotient. We will define a correspondence between convex, compact submanifolds of $M$ and conformal metrics on $\Sigma$ whose projective shape operator has eigenvalues in the interval $(-1,1]$.

We begin with the classical classification of Kleinian groups by their conformal boundary. We'll state this theorem in a somewhat unusual way that will be convenient for our purposes. 

\begin{theorem}\label{AB_param}
Let $\bar X$ be a smooth, compact, 3-manifold with boundary with $X$ the interior of $\bar X$. Assume that $\bar X$ is homeomorphic to a convex co-compact hyperbolic 3-manifold and let $\hat g$ be a smooth Riemannian metric on $\del \bar X$. Then there exists a convex, co-compact hyperbolic 3-manifold $\bar M_{\hat g}$ and a diffeomorphism $\phi_{\hat g}\colon \bar X\to \bar M_{\hat g}$ such that the restriction of $\phi$ to $\del \bar X$ is a conformal map from the conformal class of $\hat g$ to $\del_c M$.

Furthermore if $\hat g_0$ is another Riemannian metric on $\del \bar X$ and there is a diffeomorphism of $\bar X$ to itself such the restriction to $\del \bar X$ is a conformal map from the conformal class of $\hat g$ to the conformal class of $\hat g_0$ then there is an isometry from $\bar M_{\hat g}$ to $\bar M_{\hat g_0}$ whose restriction to the boundary is the given conformal map.
\end{theorem}

\newcommand{\GG}{{\mathcal G}}
\newcommand{\PP}{{\mathcal P}}
Let $\GG(\del\bar X)$ be the space of smooth Riemannian metrics on $\del\bar X$. Then for $\hat g \in \GG(\del\bar X)$ we can use the diffeomorphism $\phi_{\hat g}\colon \bar X \to \bar M_{\hat g}$ to pull back the projective structure on $\del \bar M_{\hat g}$ to a projective structure $\Sigma_{\hat g}$ on $\del\bar X$. Let ${\hat B}_{\hat g}$ be the projective shape operator and $\GG_c(\del \bar X) \subset \GG(\del \bar X)$ the subspace of metrics where ${\hat B}_{\hat g}$ has eigenvalues in the interval $(-1,1]$.

The pair $(\hat g, {\hat B}_{\hat g})$ has a dual pair $(g, B_g)$. The dual pair determines an isometric immersion $f_{\hat g} \colon \del \bar X \to M_{\hat g}$. To construct $f_{\hat g}$ we observe that the metric $g$  and shape operator lift to the universal cover of each component of $\del \bar X$. We can then construct the isometric immersion on the universal cover and observe that it will descend to the map $f_{\hat g}$  by the naturality of the constructions. In fact the same argument gives a local isometry $F \colon \cF_{( g,  B_g)}([0,\infty]) \to M_{\hat g}$.

The hyperbolic manifold $M = \Hs/\Gamma$ is {\em Fuchsian} if there is a $\Gamma$-invariant copy of $\Hp$ in $\Hs$. This equivalent to the limit set of $\Gamma$ being contained in a round circle. The quotient $\Hp/\Gamma$ is the unique totally geodesic surface in $M$ whose inclusion in $M$ is a homotopy equivalence.

A convex co-compact hyperbolic 3-manifold is Fuchsian if its limit set is a round circle.
\begin{lemma}\label{convex_embed}
Given $\hat g\in \GG_c(\del \bar X)$ the dual immersion
$$f_{\hat g}\colon \del \bar X \to M_{\hat g}$$
is an embedding and the image bounds a convex manifold homeomorphic to $\bar X$ or $M_{\hat g}$ is Fuchsian and $f_{\hat g}$ is 2-to-1 with image a totally geodesic surface. In this last case ${\hat B}_{\hat g}$ is the identity.
\end{lemma}

{\bf Proof:} 
We will show that the map $F\colon \cF\to \bar M_{\hat g}$ is an embedding unless $M_{\hat g}$ is Fuchsian and $\hat g$ is the hyperbolic metric.
The restriction of $F$ to $\del\bar X\times \{\infty\}$ is a diffeomorphism to $\del \bar M_{\hat g}$ and $F$ is an immersion everywhere. Since $\del \bar X$ is compact this implies that $F$ is an embedding on $\cF_{(g, B_g)}([T, \infty])$ for large $T$. Let $T_0$ be the infimum of $T$ where the map is an embedding. Then there must be $p\neq q$ in $S$ such that $F(p, T_0) = F(q, T_0)$ and the $F$-image of $S\times \{T_0\}$ will be tangent at this point of intersection. Furthermore their normals will be pointing in opposite directions. Convexity then implies that the surfaces are totally geodesic and 2-to-1 on a neighborhood of $p$ and $q$. This open condition is also a closed condition so $F$ will be totally geodesic and 2-to-1 on the components of $\del \bar X \times \{T_0\}$ that contain $(p, T_0)$ and $(q, T_0)$. Since the surfaces $\del\bar X\times \{t\}$ are strictly convex when $t>0$ this can only happen when $T_0 = 0$.
\eproof

We can also start with a smooth, convex submanifold $N$ of a convex co-compact hyperbolic 3-manifold $M$ and produce a conformal metric $\hat g$ on $\del_c M$ such that $N = N_{\hat g}$.

\begin{lemma}\label{convex_subman}
Let $N \subset M$ be a  convex, compact submanifold of a convex co-compact hyperbolic manifold $M$. If $\Sigma$ is the projective boundary of $M$ then there is a unique conformal metric $\hat g$ on $\Sigma$ such that $N_{\hat g} = N$.
\end{lemma}

{\bf Proof:} We'll assume that $M = \Hs$ and $N$ is closed convex subset of $\Hs$  with smooth boundary and we will find a conformal metric on $\Omega=\chat\smallsetminus \bar N$. Then for each $z\in \Omega$ there is a unique horosphere $\hh_z$ that intersects $N$ in a single point $p$.  We then let $\hat g$ be the conformal metric on $\Omega$ with $\hat g = \nu_p$ at $z$. Then $\del N$ is the Epstein surface for $\hat g$. By Theorem \ref{epstein_equal} we have that $\del N$ is (the image of) the dual immersion for the pair $(\Omega, \hat g)$.

When $N$ is a smooth, convex submanifold of a convex co-compact manifold take the pre-image of $N$ in the universal cover $\Hs$ and apply the construction from the previous paragraph to get an equivariant conformal metric on the domain of discontinuity. This will descend to a metric $\hat g$ on $\del\bar M$ with $N = N_{\hat g}$. \eproof

We note that the construction of the metric $\hat g$ works even when $N$ is convex but not smooth. However, in this case $\hat g$ will not be smooth so we cannot construct the dual metric $g$ as in the smooth case.  This will be discussed further below.

\section{W-volume}
As in the previous section we fix compact, hyperbolizable 3-manifold $\bar X$. We first define the {\em W-volume} on $\cG_c(\del\bar X)$. We will then extend it to $\cG(\del\bar X)$. Recall that associated to $\hat g \in \cG_c(\del\bar X)$ we have a convex, hyperbolic 3-manifold $N_{\hat g}$ that is homeomorphic $\bar X$ 
\begin{itemize}
\item 
 $M_{\hat g}$ is a convex co-compact hyperbolic structure on $X$ with $\hat g$ a conformal metric on $\del_c M_{\hat g}$;
 
 \item ${\hat B}_{\hat g}$ is the projective shape operator for the metric $\hat g$ and the projective boundary of $M_{\hat g}$;

\item $(g, B_g)$ is the pair dual to $(\hat g, {\hat B}_{\hat g})$;

\item $N_{\hat g}$ is a smooth convex submanifold of $M_{\hat g}$;

\item the induced metric on $\del N_{\hat g}$ is $g$ and the shape operator is $B_g$.
\end{itemize}
The mean curvature $H_g$ of $\del N_{\hat g}$ is one half the trace of $B_g$ and is a function on $\del N_{\hat g} = \del \bar X$. We define the boundary term
$$\cB(\hat g) = \int_{\del\bar X} H_g dA_g$$
where $dA_g$ is the area form for $g$. Then the W-volume is the function $W\colon \cG_c(\del\bar X)\to \R$ given by
$$W(\hat g) = \vol(N_{\hat g}) - \frac12\cB(\hat g).$$

The space of metric $\cG_c(\del\bar X)$ is an open subspace of the space of symmetric 2-tensors on $\del X$ with the $C^\infty$-topology. We are interested in calculating the variation of W-volume in this topology. More explicitly if $\hat g_t$ is a smooth family of metrics in $\cG_c(\del\bar X)$ then we want to find
$$\left.\frac{d}{dt} \right|_{t=0} W(\hat g_t)$$
The smooth family $\hat g_t$ determines smooth families $B_{\hat g_t}$, $g_t$, $B_{g_t}$, etc. Furthermore if $\delta \hat g$ is the time zero derivative of $\hat g_t$ then $\delta \hat g$ determines the time zero derivatives $\delta \hat B$, $\delta g$, $\delta B$, etc. for the associated objects.

We begin with the generalized Sch\"afli formula for the change in volume of $N_{\hat g}$ due to Rivin-Schlenker. They give a formula for the variation of the volume of a compact hyperbolic $n$-manifold, or more generally, a compact Einstein manifold. We will only state it in the setting we will use and include a proof due to Souam (\cite{Souam04}) for completeness. While we only need the statement for smooth manifolds with boundary we state it for manifolds with corners as this will come up in the proof. For a manifold with corners at each point $x$ in the co-dimension two face $E$ there is {\em exterior dihedral angle} $\theta_E(x)$.

\begin{theorem}[Rivin-Schlenker, \cite{rivin:schlafli} and Souam, \cite{Souam04}]\label{schlafli}
Let $M$ be a compact hyperbolic 3-manifold with corners and $\delta g$ a variation of the hyperbolic metric. Then
$$\delta V = \frac12\sum_E \int_E \delta \theta_E dE+ \int_{\del M}  \langle g, \delta B\rangle_g + \frac{1}{2}\langle \delta g, B\rangle_g$$
\end{theorem}

{\bf Proof:} We can decompose $M$ into finitely many manifolds with corners each of which embeds in $\Hs$. If we prove the formula for each piece then the terms on the faces that are paired will cancel. The corners will come in two types: some will come from corners of the original manifold while others will be new. In the latter case, after the corners in the pieces are glued the will become a smooth part of the interior, in which cases the total angle will be $0$, or a smooth part of the boundary, in which case the total angle will be $\pi$. In both cases, this implies the total change in dihedral angle is zero and and they will not contribute to the last term of the  formula.

We now assume that $M$ is a manifold with corners that embeds in $\Hs$. It will be useful to think of $M$ as both a fixed smooth manifold with corners and as a subspace of $\Hs$. As $M$ embeds in $\Hs$ there is a smooth vector field $\xi$, defined on a neighborhood of $M$, with flow $\phi_t$, $g_t = \phi_t^* g_{\Hs}$ and $\delta g$ is the time zero derivative of $g_t$. 
We also assume that $\xi$ is non-zero on $M$. If it isn't, by compactness, the norm of $\xi$ is bounded on $\del M$ so we can add an infinitesimal isometry to $\xi$ that is larger than this bound. Then the sum will be non-zero.

The family of metrics $g_t$ on $M$ also determines a family of shape operators $B_t$ and second fundamental forms ${\rm II}_t$ on $\del M$.  When there is a subscript we view the object as being on the fixed manifold $M$. 

Fix $u \in T_x M$ and define a parameterized surface $\sigma\colon (-\epsilon, \epsilon)^2 \to \Hs$ such that:
$$\sigma(0,0) = x \qquad\quad \sigma(s,t) = \phi_t(\sigma(s,0)) \qquad\quad \sigma_*(0,0) \frac{\del}{\del s} = u.$$
This surface will be singular if $u$ and $\xi$ are parallel. We define the vector field $n$ along $\sigma$ such that $n(s,t)$ is the normal vector to $\phi_t(\del M)$ at $\sigma(s,t)$.  If we fix $t$ the covariant derivative of $n$ along the path $s \mapsto \sigma(s,t)$ defines another vector field along $\sigma$ which we label $\hat B u$ as we have
$$\hat B u(0,t) = (\phi_t)_*(B_t u).$$

We can also evaluate the curvature tensor in terms of covariant derivatives along $\sigma$:
\begin{eqnarray*}
R(\xi, u) n & = & \frac{D}{dt} \left(\frac{Dn}{ds}\right) - \frac{D}{ds}\left(\frac{Dn}{dt}\right)\\
& = & \nabla_\xi \hat B u - \nabla_u\left(\nabla_\xi n\right).
\end{eqnarray*}
Note that while the parameterized surface may be singular, the vector field $\hat B u$ is a well defined vector field along the integral curve of $\xi$ through $x$ and $\nabla_\xi$ is a smooth vector field on $\del M$.

We let $\la, \ra$ be the usual inner product on $\Hs$.
Since we are in constant curvature $=-1$ we have
$$R(\xi, u) n = \la n, \xi \ra u - \la n, u\ra \xi = \la n, \xi \ra u.$$
Combining this equation with the one above gives
$$\nabla_\xi \hat Bu =  \la n, \xi\ra u + \nabla_u \nabla_\xi n.$$

Given any vector $v \in T_x M$ we define $\tilde v$ to be the parallel vector field on $(-\epsilon, \epsilon)^2$ such that $\sigma_*(0,0) \tilde v = v$. On the integral curve of $\xi$ through $x$ this determines a vector field which we also label $\tilde v$. We then have
$$g_t(u,v) =  \la\tilde u(\phi_t(x)), \tilde v(\phi_t(x))\ra$$
and
$$\delta g(u,v) = \xi \la\tilde u, \tilde v\ra  = \la\nabla_\xi \tilde u, v\ra + \la u, \nabla_\xi \tilde v\ra.$$
As ${\rm II}_t(u, v) = \la\hat B u, \tilde v\ra$ at $\phi_t(x)$ then
$$\delta {\rm II}(u,v) = \la\nabla_\xi \hat B u,  v\ra + \la B u, \nabla_\xi \tilde v\ra.$$
On the other hand
$$\delta{\rm II}(u,v) = \delta(g(Bu, v)) =  \delta g(Bu, v) + \la\delta B u, v\ra.$$
Combining we get
\begin{eqnarray*}
\la\delta B u, v\ra& = &\delta{\rm II}(u,v) - \delta g(Bu, v) \\
&= & \la \nabla_\xi \hat B u,  v \ra + \la  B u, \nabla_\xi  \tilde v\ra - \la \nabla_\xi \widetilde{Bu},  v \ra - \la Bu ,\nabla_\xi \tilde v\ra\\
&= & \la \nabla_\xi \hat B u,  v \ra -\la \nabla_\xi \widetilde{Bu},  v \ra .
\end{eqnarray*}
Substituting for $\nabla_\xi \hat Bu$ we get
$$\la n, \xi \ra \la u, v\ra = - \la \nabla_{\tilde u} \nabla_\xi n, v \ra + \la \nabla_\xi \widetilde{Bu}, v \ra + \la\delta B u, v\ra.$$
This expression gives a relationship between four bilinear maps in $u$ and $v$. We will take the trace of each of them. For the term on the left the trace is $2\la n, \xi \ra$. The trace of $\la \nabla_u \left(\nabla_\xi n\right), v\ra$ is the divergence of $\nabla_\xi n$ with respect to the induced metric on $\del M$. For the middle term on the right we assume that $u$ is an eigenvector of $B$ with eigenvalue of $\lambda$ and noting that $\widetilde{Bu} = \lambda \tilde u$ we get
$$\delta g(Bu, u)  = \la \nabla_\xi \widetilde{Bu}, u \ra + \la Bu, \nabla_\xi \tilde u \ra = 2\lambda \la \nabla_\xi \tilde u , u \ra = 2\la \nabla_\xi \widetilde{Bu}, u \ra$$ 
and therefore the trace of $\frac12 \delta g(Bu,v)$ is equal to the trace of $\la \nabla_\xi \widetilde{Bu}, v \ra$.
In our notation from section \ref{pair_notation} this gives
$$\la \xi, n \ra dA = -\frac12 \div_{\del M}\left(\nabla_\xi n\right) dA + \la g, \delta B\ra_g + \frac12\la \delta g, B\ra_g$$
By the divergence theorem we have
$$\delta V = \int_{\del M} \langle \xi, n\rangle dA$$
where $n$ is the outward normal on $\del M$ and $dA$ is the area form. We note that the divergence theorem holds for manifolds with corners.
Thus
$$\delta V = \int_{\del M} \left(  -\frac12 \div_{\del M}\left(\nabla_\xi n\right) dA + \la g, \delta B\ra_g + \frac12\la \delta g, B\ra_g\right).$$

Finally we  consider the divergence term $\div_{\del M}\left(\nabla_\xi n\right) dA $.
 We note that as $\langle n, n\rangle = 1$ then differentiating we have $\langle \nabla_\xi n, n\rangle =0.$ Thus $\nabla_\xi n$ is a tangent vector on $\partial M$ and smooth except on the edges $E$. Each edge $E$  belongs to exactly two faces $F^\pm_E$  whose normal $n^\pm_E$  extends to $E$. Further we define the conormal $v_E^\pm$ to $E$ to be the normal to $E$ in the face $F^\pm_E$. Then the divergence theorem, applied to each face, gives
$$\int_{\partial M}\div_{\partial M}(\nabla_\xi n)dA =  \sum_{E} \int_{E} \langle \nabla_\xi n^+_E, v^+_E\rangle + \langle \nabla_\xi n^-_E, v^-_E\rangle{d\ell}.$$
As
$\langle n^+_E, n^-_E\rangle = \cos\theta_E $, differentiating, we get
$$\langle \nabla_\xi n^+_E, n^-_E\rangle +\langle n^+_E, \nabla_\xi n^-_E\rangle  = -\sin\theta_E \delta\theta_E.$$
Also as $n^+_E$ and $v^+_E$ are an orthonormal basis for the 2-plane orthogonal to $E$ and $\langle \nabla_\xi n^+_E, n^+_E\rangle = 0$ we have 
$$\nabla_\xi n^+_E - \la \nabla_\xi n^+_E, v^+_E\ra v^+_E$$
is tangent to $E$. Therefore the inner produce with $n^-_E$ is zero giving
$$\langle \nabla_\xi  n^+_E, n^-_E\rangle =  \langle \nabla_\xi  n^+_E, v^+_E\rangle\langle v^+_E,n^-_E\rangle.$$
 As $\langle v^+_E,n^-_E\rangle =  \sin(\theta_E)$ (and similarly  $\langle v^-_E,n^+_E\rangle =  \sin(\theta_E)$) we have
$$ -\sin(\theta_E)\delta\theta_E= (\langle \nabla_\xi n^+_E, v^+ _E\rangle+\langle \nabla_\xi n^-_E, v^-_E\rangle)\sin(\theta_E) .$$
Canceling, we obtain the equation
$$   \langle \nabla_\xi n^+_E, v^+_E\rangle + \langle \nabla_\xi n^-_E, v^-_E\rangle = -\delta\theta_E.$$
The result follows. \eproof

Next we state Krasnov and Schlenker's variation formula for W-volume in terms of the interior data.
\begin{theorem}[{Krasnov-Schlenker, \cite{KS08}}]
Let $M$ be a compact hyperbolic 3-manifold with boundary and $\delta g$ a variation of the hyperbolic metric. If $\delta W$ is the induced variation of the W-volume then
$$\delta W = \frac{1}{2}\int_{\del M} \langle g, \delta B\rangle_g+ \langle \delta g, B_0\rangle_g$$
where $(g, B)$ is the dual pair, $B_0$ is the traceless part of $B$ and $\delta g$ and $\delta B$ are the induced variations.
\label{var-in}\end{theorem}

{\bf Proof:}
Let $\eta$ be the tangent bundle isomorphism such that $\delta g = 2g\cdot \eta$. Then 
$$\delta(dA_g) = \Tr(\eta) dA_g =  \frac{1}{2}\Tr_g(\delta g) dA_g = \langle \delta g, \Id \rangle_g.$$
Applying the formula for $\delta V$ from Theorem \ref{schlafli} we have
\begin{eqnarray*}\delta W &=& \delta V -\delta\cB\\
 &=& \int  \langle g, \delta B\rangle_g + \frac{1}{2}\langle \delta g, B\rangle_g- \frac{1}{2}\int \delta H dA_g + H\delta(dA_g)\\ &= & \int  \langle g, \delta B\rangle_g + \frac{1}{2}\langle \delta g, B\rangle_g - \frac{1}{2}\int \langle g, \delta B\rangle_g+ H \langle \delta g,\Id \rangle_g\\
&=& \frac{1}{2}  \int  \langle g,\delta B\rangle_g  +  \langle\delta g, B-H\Id\rangle_g \\ &=& \frac{1}{2}\int \langle g, \delta B\rangle_g+ \langle \delta g, B_0\rangle_g.
\end{eqnarray*}
\eproof

Krasnov-Schlenker also gave a variational formula for W-volume in terms of data at infinity. We show that this formula  follows  directly from Lemma \ref{invert} comparing the variation of a pair of metrics.  
\begin{theorem}[{Krasnov-Schlenker, \cite{KS08}}]
Given $\hat g \in \cG_c(\del\bar X)$ and $\delta \hat g$ a variation of $\hat g$. If $\delta W$ is the induced variation of the W-volume then
 $$\delta W = -\frac{1}{4}\int_{\del \bar X} \langle \hat g, \delta \hat B\rangle_{\hat g}+ \langle \delta \hat g, \hat B_0\rangle_{\hat g}.$$
where $\hat B_0$ is the traceless part of the shape operator at infinity $\hat B$. 
\label{var-out}\end{theorem}
{\bf Proof:}
We have $\hat g = P^*g$ and $g = \hat P^*\hat g$  where 
   $$P =\Id + B \qquad \mbox{and}\qquad \hat P = \frac{1}{2}(\Id +\hat B).$$ 
   Further as $\hat g \in \cG_c(\del\bar X)$ then $\det(\Id+B) > 0$.
Then by    Lemma \ref{invert}
    $$ \langle g, \delta P\rangle_g+  \langle \delta g,  P_0\rangle_g  = -\langle \hat g, \delta \hat P\rangle_{\hat g}-  \langle \delta \hat g,  \hat P_0\rangle_{\hat g}.$$   
   As 
   $$P_0 = B_0\qquad \delta P = \delta B,
    \qquad \mbox{and}\qquad P_0 = \frac{1}{2}\hat B_0 \qquad \delta \hat P = \frac{1}{2}\delta \hat B,$$ 
    then by linearity
$$ \langle g, \delta B\rangle_g+  \langle \delta g,  B_0\rangle_g  = -\frac{1}{2}\left( \langle \hat g, \delta \hat B\rangle_{\hat g}+  \langle \delta \hat g,  \hat B_0\rangle_{\hat g}\right).$$
 \eproof

We finally give the variational formula in terms of the Osgood-Stowe differential. 

\begin{theorem}\label{Wvar}
The variational formula for W-volume on $\mathcal G_c(\partial \bar X)$ at a metric $\hat g$ is given by
$$\delta W = \frac{1}{4}\int_{\partial \bar X} \delta KdA_{\hat g} - 2\Re \int_{\partial\bar X} \langle Q(\hat g),\mu(\delta \hat g)\rangle_{\hat g} .$$
\end{theorem}
{\bf Proof:}
By  \cite{KS08} (see also Theorem \ref{var-out} above) 
 $$\delta W = -\frac{1}{4}\int \langle \hat g, \delta \hat B\rangle_{\hat g}+ \langle \delta \hat g, \hat B_0\rangle_{\hat g}.$$
 By the projective Gauss equation we have $\Tr(\hat B) = -2K(\hat g)$ which gives
 $$\langle \hat g, \delta \hat B\rangle_{\hat g} = \frac{1}{2} \Tr(\delta\hat B)dA_{\hat g}  = -\delta K(\hat g)dA_{\hat g}.$$ 
 By Theorem \ref{equivalent_proj} we have 
$$\hat{\II}_0 = 2\os(\Sigma,\hat g) = 2Q(\Sigma,\hat g) + 2\overline{Q(\Sigma,\hat g)} .$$
We have  $\delta \hat g = 2\hat g\cdot \eta$ where  the strain $\eta_0$  is 
$$\eta_0 = \nu(\delta \hat g) \ddz\otimes d\zbar + \overline{\nu(\delta \hat g)} \ddzbar\otimes dz$$
and $\mu(\delta \hat g) = \nu(\delta \hat g) \ddz\otimes d\zbar$.
Thus  as $\hat{\II}_0 = \hat g\cdot \hat B_0$ we have by Lemma \ref{pairing}
$$\langle \delta \hat g, \hat B_0\rangle_{\hat g} = \langle 2\hat g\cdot \eta , \hat B_0\rangle_{\hat g}  = \langle \hat g\cdot \hat B_0, 2 \eta \rangle_{\hat g}  = \langle \hat{\II}_0 ,2\eta_0\rangle_{\hat g} = 8 \Re \langle Q(\Sigma,\hat g), \mu(\delta \hat g) \rangle_{\hat g} .$$
 Substituting, the variational formula  follows.  \eproof

\subsection{Conformal variations}
If we fix the conformal class of the metric $\hat g$ then we are fixing the convex co-compact hyperbolic 3-manifold $M_{\hat g}$ so we are studying the W-volume of convex submanifolds of $M_{\hat g}$. We note that the second term in Theorem \ref{Wvar} vanishes. This will allow us to make explicit computations which will enable to us to give a general definition of W-volume for any smooth metric on $\del \bar X$.

We now prove Theorem \ref{Wdiff} but only for metrics in $\cG_c(\del \bar X)$. Once we have done this we will be able to define W-volume for general metrics.
\medskip

\begin{prop}\label{WdiffGc}
Let $\hat g \in \cG_c(\del \bar X)$ and $u$ a smooth function on $\del\bar X$ such that $\hat g_t = e^{2tu}\hat g$ is in $\cG_c(\del\bar X)$ for $t \in [0,1]$. Then
$$W(\hat g_1) -W(\hat g_0) = -\frac{1}{4} \int_{\del \bar X} u\left( \Omega_{\hat g_0}+ \Omega_{\hat g_1}\right).$$
\end{prop}

{\bf Proof:}
Let $K_t = K(\hat g_t)$ be the Gaussian curvature of $\hat g_t$ and $dA_t$ be the area forms. For this proof we need more than the time zero derivative so we let $\dot K_t$, $\dot{dA}_t$, etc. represent time $t$ derivatives. By Gauss-Bonnet $$\int K_t dA_t = 2\pi\chi(\partial M)$$
so
$$\int \dot K_t dA_t = -\int K_t\dot{dA}_t.$$
As $dA_t = e^{2tu}dA$ we also have $\dt{dA_t}\  = 2udA_t$.
As the family $\hat g_t$ fixes the conformal structure there is only the curvature term in the formula from Theorem \ref{Wvar} and we have
$$\dot W_t = -\frac14 \int K_t \dot{dA}_t = -\frac12 \int uK_t dA_t.$$
The formula for curvature under conformal change (see \cite[Theorem 1.159]{Besse}) gives
$$K_t = e^{-2tu}(K(\hat g)-t\Delta u),$$ where $\Delta$ is the Laplacian of $\hat g$. Therefore
$$\dot W_t = -\frac{1}{2}\int u(K(\hat g)-t\Delta u) dA_{\hat g} = -\frac{1}{2}\int uK(\hat g)dA_{\hat g} +\frac{1}{2}t\int u\Delta u dA_{\hat g}$$
Integrating we have
$$W(\hat g_1) = W(\hat g) -\frac{1}{2}\int u\cdot K(\hat g) dA_g + \frac{1}{4}\int u\Delta(u)dA_g .$$
Then
$K(h) = e^{-2u}(K(g) -\Delta u)$ giving
$$W(h) = W(g) -\frac{1}{4}\int u\cdot K(g) dA_g - \frac{1}{4}\int uK(h) dA_h  $$
Therefore
$$W(h) - W(g) = -\frac{1}{4}(\Omega_{g}+ \Omega_{h})(u).$$ 
\eproof

We then get Krasnov-Schlenker's scaling formula for metrics in $\cG_c(\del \bar X)$.
\begin{cor}[{Krasnov-Schlenker, \cite{KS08}}]\label{scaling}
If $g \in \cG_c(\del\bar X)$ and $s \ge 0$ then
$$W\left(e^{2s}\hat g\right) = W(\hat g) - s\pi\chi(\del\bar X).$$
\end{cor}

{\bf Proof:} When $s \ge 0$ the metrics $e^{2t\cdot s} \hat g$ are in $\cG_c(\del \bar X)$ for $t\in [0,1]$ so we  can apply Proposition  \ref{WdiffGc} with $u \equiv s$. The result then follows by observing that
$$\int_{\del \bar X} s \left(\Omega_{\hat g_0} + \Omega_{\hat g_1}\right) = 2s(2\pi\chi(\del \bar X))$$
by Gauss-Bonnet. \eproof

We can now define W-volume for an arbitrary $\hat g \in \cG(\del \bar X)$ as follows. By compactness the absolute value of  eigenvalues of $\hat B_{\hat g}$ are bounded by some constant $\Lambda>1$. Then for and $s > \frac12 \log\Lambda$ we have that the projective shape operator of the metric $e^{2s}\hat g$ has eigenvalues in the interval $(-1,1)$ by Proposition \ref{dual flow}. We then set
$$W(\hat g) = W(e^{2s}\hat g) + s\pi\chi(\del\bar X).$$
and observe that  $W(\hat g)$  doesn't depend on the choice of $s$ by Corollary \ref{scaling}.

We now observe that the results on W-volume on $\cG_c(\del \bar X)$ extend to $\cG(\del \bar X)$. As already noted the scaling property  is due to Krasnov-Schlenker (see \cite{KS08}).

\begin{theorem}\label{Wvol_omnibus}
The following hold for $W$ on $\cG(\del \bar X)$.
\begin{itemize}
\item If $\hat g \in \cG(\del\bar X)$ and $s \in \R$ then
$$W\left(e^{2s}\hat g\right) = W(\hat g) - s\pi\chi(\del\bar X).\qquad \mbox{\bf  [Krasnov-Schlenker, \cite{KS08}]}$$
\item The variational formula for W-volume on $\cG(\partial \bar X)$ at a metric $\hat g$ is given by
$$\delta W = \frac{1}{4}\int_{\partial \bar X} \delta KdA_{\hat g} - 2\Re \int_{\partial\bar X} \langle Q(\hat g),\mu(\delta \hat g)\rangle_{\hat g} .$$
\item Let $\hat g_0, \hat g_1 \in \cG(\del \bar X)$ and $u$ a smooth function on $\del\bar X$ such that $\hat g_1 = e^{2u}\hat g_0$. Then
$$W(\hat g_1) -W(\hat g_0) = -\frac{1}{4} \int_{\del \bar X} u\left( \Omega_{\hat g_0}+ \Omega_{\hat g_1}\right).$$
\end{itemize}
\end{theorem}

{\bf Proof:}
The first item follows trivially from the definition.

For the variational formula we let $\hat g_t$ be a variation with $\hat g_0 = \hat g$ and $\delta W, \delta K , \delta \hat g$ be the derivatives at $t=0$. We fix $s$ such that $e^{2s}\hat g_t \in \cG_c(\del\bar X)$ for $t \in (-\epsilon,\epsilon)$. Then by definition
$$W(\hat g_t) = W(e^{2s}\hat g_t) - s\pi\chi(\del\bar X).$$
It follows that 
$$\delta W = \frac{1}{4}\int_{\partial \bar X} \delta K(e^{2s}\hat g_t) dA_{e^{2s}\hat g} - 2\Re \int_{\partial\bar X} \langle Q(e^{2s}\hat g),\mu(\delta (e^{2s}\hat g))\rangle_{e^{2s}\hat g}.$$
We have $K(e^{2s}\hat g_t) = e^{-2s} K(\hat g_t)$ giving $\delta K(e^{2s}\hat g_t) = e^{-2s}\delta K$. Also $dA_{e^{2s}\hat g} = e^{2s} dA_{\hat g}$.
Further  $Q(e^{2s}\hat g) = Q(\hat g)$, and $\mu(\delta(e^{2s} \hat g)) = \mu(\delta\hat g)$.  Finally  $\langle \cdot,\cdot \rangle_{\hat g} = \langle \cdot,\cdot \rangle_{e^{2s}\hat g} $ as the pairing only depends on the conformal class of the metric. Thus
$$\delta W = \frac{1}{4}\int_{\partial \bar X} \delta KdA_{\hat g} - 2\Re \int_{\partial\bar X} \langle Q(\hat g),\mu(\delta \hat g)\rangle_{\hat g} .$$
For the last item, let $\hat g_t = e^{2tu}\hat g_0$. By compactness of $\partial \bar X$ we can choose a $T$ such that $e^{2T}\hat g_t \in  \cG_c(\del\bar X)$ for all $t\in [0,1]$. Also the curvature form in invariant under scaling the metric giving $\Omega_{e^{2T}\hat g_i} = \Omega_{\hat g_i}$. Thus 
$$W(\hat g_1)-W(\hat g_0) = W(e^{2T}\hat g_1)-W(e^{2T}\hat g_0)=  -\frac{1}{4} \int_{\del \bar X} u\left( \Omega_{\hat g_0}+ \Omega_{\hat g_1}\right).$$
\eproof 

We can now prove our generalization of the Krasnov-Schlenker maximality property.
\medskip

\noindent
{\bf Theorem \ref{wmax}}
{\em
Let $\hat g$ be a constant curvature metric on $\del \bar X$ and 
let $\hat h = e^{2u}\hat g$ be a smooth metric on $\del\bar X$ with the same area as $\hat g$. Then
$$W(\hat g) - W\left(\hat h\right) \geq \frac{1}{4}\|\nabla u\|_2^2$$
where $\nabla$ is the gradient of $\hat g$. }
\vspace{.1in}

{\bf Proof:}
Assume $\hat g$ has constant curvature $K < 0$. 
If $\Delta$ is the Laplacian of $\hat g$, then $K(\hat h) = e^{-2u}(K-\Delta u)$. Then by  Theorem  \ref{Wdiff}
$$W(\hat g) - W\left(\hat h\right) =\frac{1}{4} \int u(\Omega_{\hat g}+\Omega_{\hat h}) = -\frac{1}{4}\int u(\Delta u-2K) dA_{\hat g}.$$
Then by the divergence theorem applied to $u\nabla u$ we have Green's identity
$$0 = \int \nabla\cdot(u\nabla u) dA_{\hat g} = \int u\Delta u dA_{\hat g} + \int  \nabla u\cdot \nabla u \ dA_{\hat g} = \int u\Delta u dA_{\hat g} +\|\nabla u\|_2^2.$$
Thus as $t \leq e^t-1$ for all $t \in \R$ and $K$ is non-positive,
$$W(\hat g)-W(\hat h) -\frac{1}{4}\|\nabla u\|_2^2= \frac{K}{4}\int 2u dA_{\hat g} \geq \frac{K}{4}\int (e^{2u}-1)dA_{\hat g} =  \frac{K}{4}\int dA_{\hat h}-dA_{\hat g} = 0$$
as $\hat g, \hat h$ have the same area.
\eproof

\subsection{W-volume for projective structures}
We now define a relative version of W-volume and show that it naturally extends to a pair of conformal metrics on an arbitrary  projective structure. In  earlier work (see \cite{KSprojRvol}) Krasnov-Schlenker define this for two conformal metrics on a projective structure where one is the projective metric.

We will see that much of the theory for W-volume for convex co-compact $3$-manifolds easily extends to this W-volume for pairs. We will also show that we can use this W-volume for pairs to obtain a surprising bound on the $L^2$-norm of a projective structure in terms of the length of the associated bending measured lamination.
 
Let $\hat g_0, \hat g_1 \in \cG(\del\bar X)$ be in the same conformal class and define
$$W(\hat g_0, \hat g_1) = W(\hat g_1) - W(\hat g_0).$$
Of course this definition makes sense for any pair of metrics in $\cG(\del \bar X)$. However, when the metrics are in the same conformal class $N_{\hat g_0}$ and $N_{\hat g_1}$ are both submanifolds of the same convex co-compact hyperbolic 3-manifold. In fact if $N_{\hat g_0} \subset N_{\hat g_1}$ then
$$W(\hat g_0, \hat g_1) = \vol(N_{\hat g_1}\smallsetminus N_{\hat g_0}) - \frac12\left(\cB(\hat g_1) - \cB(\hat g_0)\right).$$

With this as motivation we now define a relative W-volume for pairs of conformal metrics on a projective structure. Let $\Sigma$ be a projective structure and let $\cG(\Sigma)$ be smooth conformal metrics on $\Sigma$ and $\cG_c(\Sigma)$ the subspace of metrics such that the projective shape operator has eigenvalues in the interval $(-1,1]$. 
For  $\hat g \in \cG_c(\Sigma)$ by Proposition \ref{flow_embeds}  the flow space $\cF_{(g,B)}([0,\infty))$ embeds in $\ext^{\rm o}(\Sigma)$ with compact complement. We define $N_{\hat g} = \ext^{\rm o}(\Sigma)\backslash\cF_{(g,B)}((0,\infty))$. Then for $\hat g,\hat h \in \cG_c(\Sigma)$ we define
$$W(\hat g, \hat h) = \vol(N_{\hat h})- \vol(N_{\hat g})-\frac{1}{2}(\cB(\hat h) -\cB(\hat g)).$$

We would like to extend this definition to arbitrary pairs of conformal metrics. Our strategy will be the same as when we were considering convex hyperbolic 3-manifolds. In particular, the variational formulas will continue to hold with the only change being that  in the definition of $W(\hat g_0, \hat g_1)$ the boundary $\cB(\hat g_0)$ is positive so in the variational formulas the term coming from the variation of $\hat g_0$ will have opposite sign. It follows that given $\hat g_0, \hat g_1 \in \cG_c(\Sigma)$  and $t\ge s \ge 0$ we have
$$W(\hat g_0, \hat g_1)  = W(e^{2t}\hat g_0, e^{2s}\hat g_1) -(t-s)\pi\chi(\Sigma).$$
Then for two arbitrary smooth metrics $\hat g_0, \hat g_1 \in \cG(\Sigma)$ we can choose $s$ and $t$ such that $e^{2t} \hat g_0, e^{2s}\hat g_1 \in \cG_c(\Sigma)$  and define
$$W(\hat g_0, \hat g_1)  = W(e^{2t}\hat g_0, e^{2s}\hat g_1) -(t-s)\pi\chi(\Sigma).$$
As before the choice of $s$ and $t$ aren't unique but the definition of $W(\hat g_0, \hat g_1)$ is independent of this choice.

The following properties of W-volume for pairs of metrics in $\cG(\Sigma)$ follows easily from Theorem \ref{Wvol_omnibus}.

\begin{theorem}\label{proj_omnibus}
The following hold for $W$ on $\cG(\Sigma)$.
\begin{itemize}
\item If $\hat g_0, \hat g_1 \in \cG(\Sigma)$ and $s,t \in \R$ then
$$W(\hat g_0, \hat g_1)  = W(e^{2t}\hat g_0, e^{2s}\hat g_1) -(t-s)\pi\chi(\Sigma).$$
\item The variational formula for W-volume at $(\hat g_0,\hat g_1) \in \cG(\Sigma)\times\cG(\Sigma)$ is
$$\delta W = \frac{1}{4}\int_{\Sigma} \delta K(\hat g_1)dA_{\hat g_1}-\delta K(\hat g_0)dA_{\hat g_0}$$
\item Let $\hat g_0, \hat g_1 \in \cG(\Sigma)$ and $u$ a smooth function such that $\hat g_1 = e^{2u}\hat g_0$. Then
$$W(\hat g_0,\hat g_1) = -\frac{1}{4} \int_{\Sigma} u\left( \Omega_{\hat g_0}+ \Omega_{\hat g_1}\right).$$
In particular if $\hat g_0,\hat g_1$ are non-positively curved with $\hat g_1 \leq \hat g_0$ then
$$W(\hat g_0,\hat g_1) \leq 0.$$
\item If $\hat g_0$ is constant curvature and $\hat g_1 = e^{2u}\hat g_0$  has the same area as $\hat  g_0$ then
 $$W(\hat g_0,\hat g_1) \leq -\frac{1}{4}\|\nabla u\|_2^2$$
 where $\nabla$ is the gradient for $\hat g_0$. 
\end{itemize}
\end{theorem}

\begin{prop}\label{boundary bound}
If $\hat g_0, \hat g_1 \in \cG_c(\Sigma)$ with $\hat g_0\le \hat g_1$ then 
$$W(\hat g_0, \hat g_1) \ge \frac12\left(\cB(\hat g_0) - \cB(\hat g_1)\right).$$
\end{prop}

{\bf Proof:}
Since $\hat g_0, \hat g_1 \in \cG_c(\Sigma)$ and $\hat g_0 \le \hat g_1$ then by Proposition \ref{mono-flow} $N_{\hat g_0} \subseteq N_{\hat g_1}$. Therefore
\begin{eqnarray*}
W(\hat g_0, \hat g_1) 
&= & \vol\left(N_{\hat g_1})-\vol(N_{\hat g_0}\right) -\frac12\left(\cB\left(\hat g_1\right) - \cB\left(\hat g_0 \right)\right)\\
&\geq& \frac12\left(\cB\left(\hat g_0\right) - \cB\left(\hat g_1\right)\right)
\end{eqnarray*}
\eproof

\subsection{$W(\hat g_h,\hat g_\Sigma)$}
To prove Theorem \ref{newbound}, we will obtain upper and lower bounds on $W(\hat g_h, \hat g_\Sigma)$. As the projective metric $\hat g_\Sigma$ is not  smooth, we first need to extend $W$. We define $\cG^0_c(\Sigma) =\overline{ \cG_c(\Sigma)} \subseteq \cG^0(\Sigma)$ where the closure is with respect to the sup-norm topology on the space of continuous conformal metrics $\cG^0(\Sigma)$. In the appendix we will show that $W$ extends continuously to $\cG^0_c(\Sigma)\times \cG^0_c(\Sigma)$ and that $\hat g_\Sigma \in  \cG^0_c(\Sigma)$. This will suffice for our purposes.

 To get an upper bound we will rescale the metrics so that they have the same area and the use the maximality property for metrics of constant curvature. For the lower bound we scale the hyperbolic metric so that its Epstein surface is convex and then we'll apply Proposition \ref{boundary bound}. For the first estimate we need a formula for the area of the projective metric $\hat g_\Sigma$. For the second we will need formulas for the areas of the scaled metrics $(\hat g_h)_t$ and $(\hat g_\Sigma)_t$ and their Epstein surfaces. We begin with hyperbolic metric:

\begin{lemma}\label{intmean}
Let $\hat g_h$ be the hyperbolic metric  on $\Sigma$ and $(g_h)_t$ the dual metrics to $e^{2t}\hat g_h$.
 Then  for $t >\frac12\log \left(1+2\left\|\Phi_\Sigma\right\|_\infty\right)$ we have
$$\area((g_h)_t) = -2\pi\chi(\Sigma) \cosh^2(t)-e^{-2t}\left\|\Phi_\Sigma\right\|^2_2 .$$

\end{lemma}

{\bf Proof:}
If $d{\hat A}_t$ is the area form for $(\hat g_h)_t$ and $dA_t$ is the area form for $g_t$ then by Proposition \ref{dual flow} we have
$$d A_t = \frac14\left|\det\left(\Id + \hat B_t\right)\right| d\hat A_t= \frac14\left|\det\left(\Id + e^{-2t} \hat B_0\right)\right| e^{2t}d\hat A_0.$$
The determinant of $\Id + \hat B_t$ is non-negative when the eigenvalues of $\hat B_t$ are $\ge -1$ and by Proposition \ref{loc_conv} this is exactly when $e^{2t} \ge 1 + 2\left\|\Phi_\Sigma\right\|_\infty$. Therefore when  $t >\frac12\log \left(1+2\left\|\Phi_\Sigma\right\|_\infty\right)$ we can drop the absolute value signs and integrate to get
$$\area((g_h)_t) =\int dA_t = \frac{1}{4}\int \det\left(\Id+e^{-2t}\hat B_0\right) e^{2t} d\hat A_0.$$
By  Proposition \ref{Bevals}, $\hat B$ has eigenvalues $1\pm 2\|\phi_\Sigma(z)\|$.  Thus
\begin{eqnarray*}
\det\left(\Id+e^{-2t}\hat B_0\right) &=& \left(1+e^{-2t}(1+ 2\left\|\Phi_\Sigma\right\|)\right)\left(1+e^{-2t}(1- 2\left\|\Phi_\Sigma\right\|)\right)\\&=&( 1+e^{-2t})^2 - 4\left\|\Phi_\Sigma\right\|^2e^{-4t}.
\end{eqnarray*}
As $\area(\hat g_h) = -2\pi\chi(\Sigma)$ we have
\begin{eqnarray*}
\area((g_h)_t) &=&  \frac{1}{4}\left(e^{t}+e^{-t})^2 (-2\pi\chi(\Sigma)\right) -e^{-2t}\left\|\Phi_\Sigma\right\|^2_2\\ &=& -2\pi\chi(\Sigma)\cosh^2(t) -e^{-2t}\left\|\Phi_\Sigma\right\|^2_2.
\end{eqnarray*}
\eproof

For the projective metric $\hat g_\Sigma$ we have the following;

\begin{lemma}\label{area_pr}
Let $\Sigma$ be a projective structure with Thurston parameterization $\left(Y_\Sigma, \lambda_\Sigma\right)$. Then
\begin{itemize}
\item
$$\area(\hat g_\Sigma) = -2\pi\chi(\Sigma) + L\left(\lambda_\Sigma\right)$$
\item
$$\area((g_\Sigma)_t) = -2\pi\cosh^2(t)\chi(\Sigma) + L\left(\lambda_\Sigma\right)\sinh(t)\cosh(t).$$
\item
$$\area((g_\Sigma)_t)-\area((g_h)_t) = L\left(\lambda_\Sigma\right)\sinh(t)\cosh(t)+e^{-2t}\left\|\Phi_\Sigma\right\|^2_2.$$
\item
$$\area((\hat g_\Sigma)_t)-\area((\hat g_h)_t) = e^{2t}L\left(\lambda_\Sigma\right).$$
\end{itemize}
\end{lemma}

{\bf Proof:} In the Thurston parameterization of projective structures $\Sigma$ is given by a pair $\left(Y_\Sigma, \lambda_\Sigma\right)$ (see \cite{Kamishima:Tan} for a complete description of Thurston's parametrization). We can assume that the support of $\lambda_\Sigma$ is a single curve $\gamma$ with weight $\theta$ so that $L(\lambda_\Sigma) = \theta\ell_{Y_\Sigma}(\gamma)$. The Epstein surface at time $t$ is the union of a Euclidean cylinder of circumference $\ell_{Y_\Sigma}(\gamma)\cosh(t)$ and width $\theta\sinh(t)$ with the complement a surface of constant negative curvature equal to $\tanh^2(t)-1$. The area of the cylinder is $\theta\ell_{Y_\Sigma}(\gamma) \sinh(t)\cosh(t) = L(\lambda_\Sigma)\sinh(t)\cosh(t)$. By the Gauss-Bonnet formula the area of the remainder of the surface is $-2\pi\cosh^2(t)\chi(\Sigma)$. This gives the formula for $\area((g_\Sigma)_t)$.

A similar computation gives the $\area(\hat g_\Sigma)$. In this metric the annulus has circumference $\ell_{Y_\Sigma}(\gamma)$ and width $\theta$ and the complementary surface is hyperbolic. The reasoning of the previous paragraph implies the formula for $\area(\hat g_\Sigma)$.

Laminations supported on single curve are dense in the space of all measure laminations and the lengths and areas will vary continuously. This gives the general case.  The formulas for the difference of  two areas then follows easily from the prior lemma. \eproof

We now give an elementary bound for $W(\hat g_h,\hat g_\Sigma)$ in terms of the length of the bending lamination.

\begin{lemma}\label{upper}
$$W(\hat g_h, \hat g_\Sigma) \le \frac14 L\left(\lambda_\Sigma\right)$$
\end{lemma}

{\bf Proof:}
By Lemma \ref{area_pr}
$$\area\left(\hat g_\Sigma\right) = 2\pi|\chi(\Sigma)|+L(\lambda_\Sigma).$$
The maximality property of Theorem \ref{proj_omnibus}  gives that if $\hat g$ is smooth of the same area as $e^{2t}\hat g_h$ then
$$W(e^{2t}\hat g_h, \hat g) \leq 0.$$
 By Proposition \ref{W_proj_extension}, W-volume is continuous on $\cG^0_c(\Sigma)\times\cG^0_c(\Sigma)$. Further by Proposition \ref{gw_proj}  we have $\hat g_\Sigma \in \cG^0_c(\Sigma)$  giving, 
$$W\left(\frac{2\pi|\chi(\Sigma)|+L(\lambda_\Sigma)}{2\pi|\chi(\Sigma)|}\hat g_h,\hat g_\Sigma\right) \leq 0.$$
Combining this with the scaling property we get
$$W(\hat g_h,\hat g_\Sigma) -\frac{1}{2}\log\left(\frac{2\pi|\chi(\Sigma)|+L(\lambda_\Sigma)}{2\pi|\chi(\Sigma)|}\right)\pi|\chi(\Sigma)| \leq 0.$$
As  $\log(1+x) \le x$ we have
$$W(\hat g_h,\hat g_\Sigma) \leq \frac{1}{2}\log\left(1+\frac{L(\lambda_\Sigma)}{2\pi|\chi(\Sigma)|}\right)\pi|\chi(\Sigma)| \leq \frac{1}{4}L(\lambda_\Sigma).$$
\eproof

Next we get a lower bound. The key here is that volumes are non-negative.

\begin{lemma}\label{lower}
$$W(\hat g_h, \hat g_\Sigma) \ge \frac{e^{-2T}\left\|\Phi_\Sigma\right\|_2^2}2 - \frac{L\left(\lambda_\Sigma\right)\cosh(2T)}4$$
where $e^{2T} = 1+2\left\|\Phi_\Sigma\right\|_\infty$.
\end{lemma}

{\bf Proof:} 
By Proposition \ref{loc_conv} the eigenvalues for the shape operator for the Epstein surface of $(\hat g_h)_T$ are non-negative so $(\hat g_h)_T $ is the metric at infinity for a  locally convex Epstein surface.  For all $t\ge 0$, $(\hat g_\Sigma)_t$ is the metric at infinity for a locally convex Epstein surface. Since $\hat g_h \le \hat g_\Sigma$ therefore $e^{2T}\hat g_h \le e^{2T}\hat g_\Sigma$ and by Proposition \ref{boundary bound}
\begin{eqnarray*}
W(\hat g_h, \hat g_\Sigma) & = & W\left((\hat g_h)_T, (\hat g_\Sigma)_T\right) \\
&\geq& \frac12\left(\cB\left((g_h)_T\right)-\cB\left((g_\Sigma)_T\right) \right).
\end{eqnarray*}
Applying Lemma \ref{area_pr} we have
\begin{eqnarray*}
\cB\left((g_h)_T\right)-\cB\left((g_\Sigma)_T\right) &= &\frac{1}{2}\left(\area((\hat g_h)_T)-\area((\hat g_\Sigma)_T)\right)-\left(\area((g_h)_T)-\area((g_\Sigma)_T)\right) \\
&=& -\frac{1}{2} L(\lambda_\Sigma)e^{2T}+L(\lambda_\Sigma)\sinh(T)\cosh(T)+e^{-2T}\left\|\Phi_\Sigma\right\|_2^2\\
&=& \frac{1}{2}L(\lambda_\Sigma)(2\sinh(T)\cosh(T)-e^{2T})+e^{-2T}\left\|\Phi_\Sigma\right\|_2^2 \\
&=& e^{-2T}\left\|\Phi_\Sigma\right\|_2^2-\frac{1}{2}L(\lambda_\Sigma)\cosh(2T).\\
\end{eqnarray*}
Therefore it follows that
\begin{eqnarray*}
W(\hat g_h, \hat g_\Sigma)& \ge &\frac12 \left(e^{-2T}\left\|\Phi_\Sigma\right\|_2^2 -\frac{L(\lambda_\Sigma)}2 \cosh(2T)\right).
\end{eqnarray*}
\eproof

We now prove Theorem \ref{newbound}. \newline
\vspace{.001in}

\noindent
{\bf Theorem \ref{newbound}}
{\em Let $\Sigma$ be a projective structure with Thurston parameterization $\left(Y_\Sigma, \lambda_\Sigma\right)$. Then

$$\|\Phi_\Sigma\|_2 \leq (1+\|\Phi_\Sigma\|_\infty)\sqrt{L(\lambda_\Sigma)}.$$
}
{\bf Proof:} Combining Lemmas \ref{upper} and \ref{lower} we get
$$\frac{e^{-2T}\left\|\Phi_\Sigma\right\|_2^2}2 - \frac{L\left(\lambda_\Sigma\right)\cosh(2T)}4\leq \frac{1}{4}L(\lambda_\Sigma)$$
giving
$$\left\|\Phi_\Sigma\right\|_2^2 \leq \frac{1}{2}L(\lambda_\Sigma)e^{2T}(1+\cosh(2T)) =   \frac{1}{4}L(\lambda_\Sigma)\left(e^{2T}+1\right)^2.$$
Applying $e^{2T} = 1+2\left\|\Phi_\Sigma\right\|_\infty$ we obtain 
$$\left\|\Phi_\Sigma\right\|_2 \leq   (1+\left\|\Phi_\Sigma\right\|_\infty)\sqrt{L(\lambda_\Sigma)}.$$
\eproof

\appendix
\section{W-volume for non-smooth metrics} 
\subsection{Hyperbolic manifolds}
As we noted above, a convex submanifold  $N$ of a convex co-compact hyperbolic 3-manifold $M$ determines a conformal metric $\hat g$ on $\del\bar M$. We would like to define the W-volume for $\hat g$ (and $N$) even when $\hat g$ is not smooth. While the volume term (the volume of $N$) is still well defined, the integral of mean curvature is not.

For these non-smooth convex manifolds it will be convenient to fix the ambient convex co-compact hyperbolic 3-manifold $M$. In Lemma \ref{convex_subman} we saw that convex submanifolds of $M$ determine conformal metrics on $\del_c M$. We let $\cG^0_c(\del\bar X; M)$ be the space of such metrics.  This will contain the space of smooth metrics denoted $\cG_c(\del\bar X;M)$. As all metrics are conformally equivalent, given $\hat g$ and $\hat h$ in $\cG^0_c(\del\bar X; M)$, we have $\hat h = e^{2u} \hat g$ for some continuous function $u$. We then define the sup-norm by $\|\hat g - \hat h\|_\infty = \|u\|_\infty$. We also have $\hat h \ge \hat g$ if $u\ge 0$ and $\hat h \le \hat g$ if $u \le 0$.

When $\del N$ is not smooth the integral of mean curvature is not defined but the area of $\del N$ is. The following proposition will then be useful for defining W-volume for non-smooth metrics.
\begin{prop}[{{\cite[Lemma 3.1]{BBB}}}]\label{mean_integral}
Let $\hat g$ be a smooth conformal metric on a projective structure $\Sigma$ and assume that the eigenvalues of the shape operator $B$ for Epstein surface are non-negative. If $g$ is the Epstein metric then
$$\int_{\Sigma} H dA = \frac12\area(\hat g) -\area(g) - \pi\chi(\Sigma)$$
where $H = \frac 12 \tr B$ is the mean curvature.
\end{prop}
With this proposition as a  motivation we can define
$$W(\hat g) = \vol(N_{\hat g}) -\frac14\area(\hat g) +\frac12\area(g) + \frac\pi2\chi(\del\bar X)$$
for $\hat g \in \cG^0_c(\del\bar X)$. When $\hat g \in \cG_c(\del\bar X)$ this agrees with our original definition. For this definition to be useful we need it to be continuous. We prove this next.

\begin{prop}\label{cont_ext}
Let $M$ be a convex co-compact hyperbolic structure on $X$. Then the W-volume is continuous on $\cG^0_c(\del \bar X;M)$ with the sup-norm topology.
\end{prop}

{\bf Proof:} Let $\hat g$ be a metric in $\cG^0_c(\del \bar X;M)$ and $N = N_{\hat g}$ the corresponding convex submanifold of $M_{\hat g}$. Then for $t \ge 0$, we have $\hat g_t \in \cG^0_c(\del\bar X)$ as $N_t = N_{\hat g_t}$ is the $t$-neighborhood of $N$. We also have that $\area(\hat g_t)$, $\area(g_t)$ (where $g_t$ is the induced metric on $\del N_t$) and $\vol(N_t)$ are continuous in $t$. Therefore, for $\epsilon>0$ we can choose a $\delta>0$ such that $\delta\le \epsilon$, $|\area(\hat g) -\area(\hat g_{2\delta})| < \epsilon$, $|\area(g) - \area(g_{2\delta})| < \epsilon$ and $|\vol(N) - \vol(N_{2\delta})| < \epsilon$.

Now assume that $\hat h \in \cG_c(\del \bar X)$ with $\hat g < \hat h < \hat g_{2\delta}$. Then $\area(\hat g) \le \area(\hat h) \le \area(\hat g_{2\delta})$. We also have that $N \subset N_{\hat h} \subset N_{2\delta}$ which implies that $\vol(N)\le V(N_{\hat h}) \le \vol(N_{2\delta})$ and, as the nearest point projections of $\del N_{\hat h}$ to $\del N$ and $\del N_{\hat g_{2\delta}}$ to $\del N_{\hat h}$ are area decreasing,  $\area(g) \le \area(h) \le \area(g_{2\delta})$. It then follows that
$$\left|\cB(\hat g) - \cB(\hat h)\right| \le \frac\epsilon2 + \epsilon = \frac{3\epsilon}2$$
and therefore
$$\left| W(\hat g) - W(\hat h) \right| \le \epsilon + \frac{3\epsilon}4 = \frac{7\epsilon}4.$$

Next we assume that $\hat h \in \cG_c(\del\bar X)$ with $\|\hat g -\hat h\|_\infty \le \delta$. Let $\hat h_{\delta} = e^{2\delta}\hat h$ and we have $\hat g \le \hat h_\delta \le \hat g_{2\delta}$. Then, by the above,
$$\left| W(\hat g) - W(\hat h_{\delta})\right| \le \frac{7\epsilon}{4}.$$
Since $W(\hat h_{\delta}) = W(\hat h) +\delta\pi|\chi(\del\bar X)|$ and $\delta < \epsilon$ this implies that
$$\left| W(\hat g) - W(\hat h) \right| \le \frac{7\epsilon}4 + \epsilon\pi|\chi(\del \bar X)| = \left(\frac74 + \pi|\chi(\del\bar X)|\right) \epsilon.$$
Continuity follows. \eproof

To extend our results for smooth metrics to continuous metrics we further need to know that smooth metrics are dense in $\cG_c^0(\del\bar X; M)$.
\begin{prop}\label{gw}
$\cG_c(\del\bar X; M)$ is dense in $\cG^0_c(\del \bar X; M)$.
\end{prop}
{\bf Proof:} Fix $\hat g \in \cG_c^0(\del\bar X; M)$ and 
 let $\rho\colon M\to [0,\infty)$ be the distance from $N_{\hat g}\subset M$. This function is convex along geodesics and is nowhere-constant on geodesics in $M\smallsetminus N_{\hat g}$. Therefore $\rho^2$ is strictly convex on $M \smallsetminus N_{\hat g}$ and by Theorem 2 in \cite{GreeneWu_Acta} we have that $\rho^2$ can be approximated by smooth, strictly  convex functions. In particular, for all $\epsilon>0$ we can find a smooth convex function $f\colon M\smallsetminus N\to \R$ such that $\|f -\rho^2\|_\infty < \epsilon^2/4$. Note that $N_{e^{2\epsilon}\hat g}$ is the $\epsilon$-neighborhood of $N_{\hat g}$ and we have that $f \ge 3\epsilon^2/4$ on $\del N_{e^{2\epsilon}\hat g}$ and $f \le \epsilon^2/4$ on $\del N_{\hat g}$. Therefore  $S = f^{-1}\left(\epsilon^2/2\right)$ will be a smooth convex surface in $N_{e^{2\epsilon}\hat g} \smallsetminus N_{\hat g}$ and it will separate $\del N_{\hat g}$ from $\del N_{e^{2\epsilon}\hat g}$. This implies that $S$ bounds a smooth,  convex submanifold $N'$ in $M$. Finally we note that if $\hat h$ is the metric for $N'$ then $\hat g < \hat h < e^{2\epsilon}\hat g$. To see this (after lifting to the universal cover) any horoball $\hh$ that is tangent to $N'$ will be disjoint from $N_{\hat g}$ and intersect will $N_{e^{2\epsilon}\hat g}$.  Therefore $\hh$ will be nested between the horoballs, based at the same point as $\hh$, tangent to $N_{\hat g}$ and tangent to $N_{e^{2\epsilon}\hat g}$. This gives the estimate and that $\|\hat g - \hat h\|_\infty < \epsilon$. As $\epsilon>0$ was arbitrary and $\hat h$ is smooth the proposition follows. \eproof

\subsection{Projective structures}
As above, we will now extend $W$ on $\cG_c(\Sigma) \times \cG_c (\Sigma)$ to general locally convex metrics. By Proposition \ref{mean_integral} we  have for $\hat g,\hat h \in \cG_c(\Sigma)$
$$W(\hat g, \hat h) = \vol(N_{\hat h})- \vol(N_{\hat g}) +\frac14\area(\hat g) -\frac12\area(g) -\frac14\area(\hat h) +\frac12\area(h).$$
We define $\cG_c^0(\Sigma)$ to be the closure of $\cG_c(\Sigma)$ in the sup-norm topology on $\mathcal \cG^0(\Sigma)$. Then repeating the proof of Proposition \ref{cont_ext} we have 
\begin{prop}\label{W_proj_extension} Let $\Sigma$ be a projective structure on a closed surface $S$. Then $W$ extends continuously to $\cG_c^0(\Sigma)\times \cG_c^0(\Sigma)$ in the sup-norm topology.
\end{prop}
Finally we show that $\hat g_\Sigma \in \cG_c^0(\Sigma)$. 

\begin{prop} If $\Sigma$ is a projective structure then $\hat g_\Sigma \in \cG_c^0(\Sigma)$.\label{gw_proj}
\end{prop}

{\bf Proof:}
The proof  follows that of Proposition \ref{gw} using Greene-Wu, except now we do our analysis on $\ext(\Sigma)$. Let $\rho\colon \ext(\Sigma) \to [0,\infty)$ be the distance to $\Sigma\times\{0\} \subseteq \ext(\Sigma)$. This function is convex along geodesics and is nowhere-constant on geodesics in $\Sigma\times (0,\infty)$. Therefore $\rho^2$ is strictly convex on $\Sigma\times (0,\infty)$ and by Theorem 2 in \cite{GreeneWu_Acta} we have that $\rho^2$ can be approximated by smooth, strictly  convex functions. We choose a smooth convex function $\Sigma\times (0,\infty) \to \R$ such that $\|f -\rho^2\|_\infty < \epsilon^2/4$. As in Proposition \ref{gw}  the surface $S = f^{-1}\left(\epsilon^2/2\right)$ will be a smooth convex surface approximating $\Sigma\times\{0\}$ in the $\epsilon$ neighborhood of $\Sigma\times\{0\}$. Specifically, the metric $g$ on $S$ will have dual metric $\hat g$ with $\|\hat g - \hat g_\Sigma\| < \epsilon$.\eproof

\bibliography{bib,math}
\bibliographystyle{math}

\end{document}

%% file: macros.tex
\newcommand{\bb}{\mathbb}

\newcommand{\cx}{{\bb C}}

\newcommand{\htwo}{{\bb H}^2}

\newcommand{\ddz}{\frac{\partial}{\partial z}}

\renewcommand{\bold}[1]{\medskip \noindent {\bf \boldmath #1
                        }\nopagebreak[4]}

\newcommand{\del}{\partial}

\newcommand{\zbar}{{\overline{z}}}

\newcommand{\chat}{\widehat{\cx}}

\newcommand{\area}{\operatorname{area}}

\newcommand{\diam}{\operatorname{diam}}

\renewcommand{\div}{\operatorname{div}}

\renewcommand{\Re}{\operatorname{Re}}

\newcommand{\sgn}{\operatorname{sgn}}

\newcommand{\Tr}{\operatorname{Tr}}
\newcommand{\tr}{\operatorname{tr}}
\newcommand{\vol}{\operatorname{vol}}

\newtheorem{theorem}{Theorem}[section]
\newtheorem{prop}[theorem]{Proposition}
\newtheorem{lemma}[theorem]{Lemma}
\newtheorem{cor}[theorem]{Corollary}

\newcommand{\cB}{{\cal B}}

\newcommand{\cF}{{\cal F}}
\newcommand{\cG}{{\cal G}}

\newcommand{\cL}{{\cal L}}

